\pdfoutput=1

\documentclass[preprint]{article}
\usepackage{amsmath}
\usepackage{amssymb}
\usepackage{amsthm}
\usepackage{graphicx}
\usepackage{epic,eepic,epsfig}
\usepackage{color}
\usepackage{caption}
\usepackage{subcaption}
\usepackage{placeins}	
\usepackage{multirow}
\usepackage{epstopdf}
\usepackage{varwidth}
\usepackage[table]{xcolor}

\newtheorem{theorem}{Theorem}[section]
\newtheorem{lemma}[theorem]{Lemma}
\newtheorem{corollary}[theorem]{Corollary}

\theoremstyle{definition}
\newtheorem{definition}{Definition}[section]

\newcommand{\revb}[1]{#1}
\newcommand{\revr}[1]{#1}

\newcommand{\bsym}[1]{\mathbf{#1}}

\definecolor{tabclr}{cmyk}{0,0,1,0}

\begin{document}

\title{A mass conservative scheme for fluid-structure interaction problems
by the staggered discontinuous Galerkin method} 

\author{
Siu Wun Cheung\thanks{Department of Mathematics, The Chinese University of Hong Kong, Hong Kong SAR(\texttt{tonycsw2905@hotmail.com})}
\and
Eric T. Chung\thanks{Department of Mathematics, The Chinese University of Hong Kong, Shatin, New Territories, Hong Kong SAR, China (\texttt{tschung@math.cuhk.edu.hk}) }
\and
Hyea Hyun Kim\thanks{Department of Applied Mathematics and Institute of Natural Sciences, Kyung Hee University,
Korea (\texttt{hhkim@khu.ac.kr})}
}

\maketitle

\begin{abstract}
In this paper, we develop a new mass conservative numerical scheme for the simulations
of a class of fluid-structure interaction problems.
We will use the immersed boundary method to model the fluid-structure interaction, while the fluid flow is governed by the incompressible Navier-Stokes equations.
The immersed boundary method is proven to be a successful scheme to model fluid-structure interactions.
To ensure mass conservation, we will use the staggered discontinuous Galerkin method to discretize the incompressible Navier-Stokes equations.
The staggered discontinuous Galerkin method is able to preserve the skew-symmetry of the convection term.
In addition, by using a local postprocessing technique, the weakly divergence free velocity can be used to
compute a new postprocessed velocity, which is exactly divergence free and has a superconvergence property.
This strongly divergence free velocity field is the key to the mass conservation.
Furthermore, energy stability is improved by the skew-symmetric discretization of the convection term.
We will present several numerical results to show the performance of the method.
\end{abstract}

\section{Introduction}
\label{sec:intro}
Fluid-structure interaction, which models the interaction of movable structures and the surrounding fluid flow, is the key to the design of many engineering problems.
There are in literature a variety of methods to model fluid-structure interactions, and among them the immersed boundary (IB) method and the immersed interface method (IIM)
are proven to be very successful.
The immersed interface method \cite{iim} was first introduced by Li, and a detailed discussion can be found in \cite{iim-book}.
The immersed boundary method was first introduced by Peskin \cite{peskin77} for the numerical approximation of blood flow around the heart valves,
and a detailed discussion on the applications of IB method is given in \cite{peskin02}.
These methods have been successfully extended to other applications.
In this paper, we will focus on the development of our scheme using the immersed boundary approach,
since it can be combined with the staggered discontinuous Galerkin method and gives a mass conservative scheme.

One key feature of immersed boundary method is that the Eulerian mesh in the Cartesian coordinate system is fixed, and the configuration of the immersed structure does not necessarily adapt to the Eulerian mesh. This avoids the high cost of mesh updating. The source term which represents the effects of the force exerted by the immersed structure on the fluid is modelled by a Dirac delta function. In the original formulation of Immersed boundary method, finite difference methods are used in spatial discretization for the governing equations of the fluid flows. Since the material points of the immersed boundary may not adapt to the Eulerian grid, the Dirac delta function needs to be approximated. The construction of approximations of the Dirac Delta function is discussed in \cite{peskin02}.

On the contrary, in finite element and other Galerkin methods, the Dirac Delta functions can be handled directly by the variational formulation and therefore approximations of the Dirac Delta functions are not needed.
In \cite{boffi03}, a finite element approach for immersed boundary method (FE-IBM) was proposed. More recent researches on FE-IBM can be found in \cite{boffi07} and \cite{boffi16}.

In this paper, we present a staggered discontinuous Galerkin immersed boundary method (SDG-IBM). IB method is used for modelling the fluid-structure interaction, and the fluid flow is modelled by incompressible Navier-Stokes equations which would be solved numerically by a discontinuous Galerkin method based on staggered meshes.
Discontinuous Galerkin methods have been applied to problems in fluid dynamics and wave propagations with great success, see for example \cite{carrero05,hdg,cockburn05,ldg,conv-diff,houston09,Liu-Shu,hpdg,shahbazi07,nguyen11}.
On the other hand, staggered meshes bring the advantages of reducing numerical dissipation in computational fluid dynamics \cite{sfvm,sdm,sdm1},
and numerical dispersion in computational wave propagation \cite{semi,fully,newdg,newdg1,meta,jcp-max,geo}.
Combining the ideas of DG methods and staggered meshes, a new class of staggered discontinuous Galerkin (SDG) methods for approximations of the incompressible Navier-Stokes equations was proposed \cite{sdg-ns1}. The new class of SDG methods inherits many good properties, including local and global conservations, optimal convergence, and superconvergence through the use of a local postprocessing technique in \cite{hdg,hdg1}. Furthermore, energy stability is achieved by spectro-consistent discretizations with a novel splitting of the diffusion and the convection term. An analysis of the SDG method for incompressible Navier-Stokes equations is given in \cite{sdg-ns2}. For a more complete discussion on the SDG method, see also \cite{newdg,newdg1,meta,jcp-max,curlcurl,convdiff,kcl-2013-Stokes-SDG,kcl-2013-FETI-DP-Stokes} and the references therein.
We remark that another class of discontinuous Galerkin methods based on space-time staggered meshes is proposed in \cite{Dumbser1,Dumbser2,Dumbser3}.

In the finite element formulation of IB method in \cite{boffi03}, the convection term was neglected and linearized Navier-Stokes equations was considered. In our proposed method, by an iterative approach and a skew-symmetric discretization of the convection term, we can also handle the convection term without losing any stability in terms of energy.
\revb{Our stability result is subject to CFL type restriction on time step since our scheme treats
the fluid structure interaction explicitly. Otherwise the implementation is not feasible due to
the presence of nonlinear term in the fluid model, 
which also requires iteration. 
We note that a stability result
without time step restriction was proven for a simple linear fluid model when fluid structure interaction was
treated implicitly using an iterative method, see \cite{causin-implicit}.
}

Another important issue of IB method is that the loss in volume enclosed by the immersed structure in the numerical approximation, which can be resolved by improving the divergence-free property of the interpolated velocity field which drives the Lagrangian markers, see \cite{peskin02} for a detailed discussion.
A key component of our method is the use of postprocessing techniques to obtain a pointwise divergence-free velocity field approximation at each time level, which is used to drive the Lagrangian markers of the immersed boundary and acts as a convection velocity in the iterative approach of solving the incompressible Navier-Stokes equations. In particular, by using the pointwise divergence-free postprocessed velocity to drive the Lagrangian markers of the immersed boundary, our method significantly resolves the numerical error of lack of volume conservation. In these regards, our method has advantages over the FE-IBM and other discontinuous Galerkin methods.

The paper is organized as follows. In Section \ref{sec:ibm}, we will have a brief discussion on the problem formulation of the IB method. Next, in Section \ref{sec:sdgibm}, we will present the derivation of SDG-IBM. In Section \ref{sec:stab}, we will provide a stability analysis of SDG-IBM. Then, in Section \ref{sec:num}, we will present extensive numerical examples to show the performance of SDG-IBM. Finally, a conclusion is given.

\section{Problem description}
\label{sec:ibm}

Suppose, for $t \in [0,T]$, in a two-dimensional domain $\Omega \subset \mathbb{R}^2$, the immersed boundary is an elastic incompressible fibre, modeled by a simple closed curve $\Gamma_t$ contained in $\Omega$. The Eulerian coordinates of $\Gamma_t$ are denoted by $\bsym{X}(s,t)$, where $0 \leq s \leq L$ is the Lagrangian coordinates labeling material points along the curve, and
\begin{equation}
\bsym{X}(0,t) = \bsym{X}(L,t) \text{ for } t \in [0,T].
\label{eq:IB_def}
\end{equation}
The motion of the fluid is described by the incompressible Navier-Stokes equations
\begin{equation}
\begin{split}
\rho \bsym{u}_t - \mu \bsym{\Delta} \bsym{u} + \rho \bsym{u}\cdot \nabla \bsym{u} + \nabla p & = \bsym{F} \; \mbox{ in }\Omega \times (0,T),\\
\text{div} \, \bsym{u} & = 0 \; \mbox{ in }\Omega \times (0,T),\\
\bsym{u} & = 0 \; \mbox { on } \partial \Omega \times (0,T),\\
\bsym{u} & = \bsym{u}_0 \; \mbox { in } \Omega \times \{0\},
\end{split}
\label{eq:IB_NS}
\end{equation}
where $p$ is the pressure with $\int_\Omega p \, dx = 0$, $\bsym{u} = (u_1, u_2)$ is the velocity and $\bsym{F} = (F_1, F_2)$ is the source term. Here $\rho$ and $\mu$ are the density and the viscosity of the fluid, respectively. 
Let $\bsym{f}(s,t)$ denote the elastic force density resulted from the deformation of the immersed boundary.
In the IB method, the force $\bsym{F}(\bsym{x},t)$ exerted on the fluid by the immersed boundary is given by
\begin{equation}
\bsym{F}(\bsym{x},t) = \int_0^L \bsym{f}(s,t) \delta(\bsym{x} - \bsym{X}(s,t)) \, ds \text{ in } \Omega \times (0,T).
\label{eq:IB_source}
\end{equation}
Finally, a no-slip condition is imposed between the immersed boundary and the fluid. The motion of the immersed boundary is described by the Euler-Lagrange equation
\begin{equation}
\begin{split}
\frac{\partial}{\partial t} \bsym{X}(s,t) & = \bsym{u} (\bsym{X}(s,t)) \; \mbox { in } [0,L] \times [0,T],\\
\bsym{X}(s,0) & = \bsym{X}_0(s) \; \mbox { in } [0,L].
\label{eq:IB_EL}
\end{split}
\end{equation}
\revb{In the current work, we only consider the case when both $\rho$ and $\mu$ are uniform. Extension to a more general model with varying $\rho$ and $\mu$ across the interface will be addressed in our future work.}

We consider a simple model with a massless closed curve $\Gamma_t$ immersed in an incompressible fluid. Suppose $\gamma$ is the tension in $\Gamma_t$ and $\tau$ is the unit tangent to $\Gamma_t$. Then the local force density $\bsym{f}$ acting on the fluid by $\Gamma_t$ is given by
\begin{equation}
\bsym{f} = \frac{\partial}{\partial s} (\gamma \tau).
\end{equation}
We assume $\gamma$ is proportional to $\left \vert \dfrac{\partial \bsym{X}}{\partial s} \right \vert$. Then we have
\begin{equation}
\gamma \tau = \kappa \dfrac{\partial \bsym{X}}{\partial s} \implies \bsym{f} = \kappa \frac{\partial^2 \bsym{X}}{\partial s^2},
\label{eq:sim_mod}
\end{equation}
where $\kappa$ is the elasticity constant of the material along the immersed boundary.

\section{Derivation of SDG-IBM}
\label{sec:sdgibm}
In this section, we will give a detailed derivation of SDG-IBM.
We will start with the temporal discretization, and then discuss the details of full discretization.
We will discuss an iterative approach of linearizing the nonlinear convection term of Navier-Stokes equations \eqref{eq:IB_NS}.
Next, we will give the construction of the staggered mesh and the construction of finite element spaces with staggered continuity property.
After that, we will explain the derivation of the SDG method and the resultant system of linear equations in each iteration.
We will also present the postprocessing technique (c.f. \cite{hdg}) to obtain a pointwise divergence-free velocity field and discuss the significance of the post-processed velocity in our method.
Then, we will move on to discuss the discretization of the source term \eqref{eq:IB_source} in the simple model \eqref{eq:sim_mod}.
Finally, we will discuss the full discretization of the Euler-Lagrange equation \eqref{eq:IB_EL}.

\subsection{BE/FE temporal discretization}
\label{sec:IB_time}
We will first discretize the continuous problem in time, and obtain a temporally discrete and spatially continuous system.
We will use backward-Euler method for the temporal discretization of Navier-Stokes equations. In order to avoid a fully implicit system of equations at each time-step, we use forward-Euler method in time discretization for Euler-Lagrange equation~\eqref{eq:IB_EL} and the fibre force~\eqref{eq:IB_source}. A similar approach was employed by \cite{boffi03}, and such an approach is regarded as the BE/FE scheme \cite{stockie99}.
\revb{We note that fully implicit scheme was considered in \cite{causin-implicit} for a simple linear fluid
model. For our nonlinear fluid model, that approach is not feasible.}

Let $K$ be the number of divisions in $[0,T]$ in the temporal domain, $\Delta t = T/K$ be the time step size and $t_n = n \Delta t$. From now on, a function with a superscript $n$ stands for evaluation of the function at time $t = t_n$. For $n = 1, 2, \ldots, K$, given $\bsym{u}^{n-1}$, our goal is to solve for $(\bsym{u}^n, p^n)$ in the following system of nonlinear PDEs:
\begin{equation}
\begin{split}
\frac{\rho}{\Delta t} \bsym{u}^n - \mu \bsym{\Delta} \bsym{u}^n + \rho \bsym{u}^n\cdot \nabla \bsym{u}^n + \nabla p^n  & = \frac{\rho}{\Delta t} \bsym{u}^{n-1} + \bsym{F}^n \; \mbox{ in }\Omega,\\
\text{div} \, \bsym{u}^n & = 0 \; \mbox{ in }\Omega,\\
\bsym{u}^n & = 0 \; \mbox { on } \partial \Omega,\\
\bsym{u}^0 & = \bsym{u}_0 \; \mbox { in } \Omega,
\end{split}
\label{eq:IB_NS2}
\end{equation}
where the source term $\bsym{F}^n$ is given by
\begin{equation}
\bsym{F}^n(\bsym{x}) = \int_0^L \bsym{f}^{n-1}(s) \delta(\bsym{x} - \bsym{X}^{n-1}(s)) \; ds \; \mbox{ in } \Omega.
\label{eq:IB_source2}
\end{equation}
On the other hand, the immersed boundary $\bsym{X}^n$ is evolved by
\begin{equation}
\begin{split}
\bsym{X}^{n} & = \bsym{X}^{n-1} + \Delta t \, \bsym{u}^{n}\left(\bsym{X}^{n-1}\right) \; \mbox{ in } [0,L],\\
\bsym{X}^{0} & = \bsym{X}_0 \; \mbox{ in } [0,L].
\end{split}
\label{eq:IB_EL2}
\end{equation}

\subsection{Linearization of Navier-Stokes equations by iterative approach}
\label{sec:NS_iter}
In our method, for solving the system \eqref{eq:IB_NS2} of nonlinear PDE at $t = t_n$, the nonlinear convection term is linearized by a sequence of Picard fixed-point iterations:
\begin{equation}
\begin{split}
\frac{\rho}{\Delta t} \bsym{u}^n_m - \mu \bsym{\Delta} \bsym{u}^n_m + \rho \bsym{V}^n_m \cdot \nabla \bsym{u}^n_m + \nabla p^n_m & = \frac{\rho}{\Delta t} \bsym{u}^{n-1} + \bsym{F}^n \; \mbox{ in }\Omega,\\
\text{div} \, \bsym{u}^n_m & = 0 \; \mbox{ in }\Omega,\\
\bsym{u}^n_m & = 0 \; \mbox { on } \partial \Omega,\\
\bsym{u}^0 & = \bsym{u}_0 \; \mbox { in } \Omega,
\end{split}
\label{eq:IB_NS3}
\end{equation}
where $\bsym{V}^n_m$ is a given pointwise divergence-free velocity field depending on $\bsym{u}^{n}_{m-1}$.

The choice of the velocity field $\bsym{V}^n_m$ in the formulation of \eqref{eq:IB_NS3} will be discussed in Section \ref{sec:post}. The SDG method for solving \eqref{eq:IB_NS3} in a particular iteration will be discussed in Sections \ref{sec:mesh}--\ref{sec:sys}. The fixed point $(\bsym{u}^n, p^n)$ of the sequence $\{(\bsym{u}_m^n, p_m^n)\}_{m=1}^\infty$ is then our solution for \eqref{eq:IB_NS2}. In practice, we set a suitable stopping criterion for the Picard fixed-point iterations when the number of iterations done is sufficient or when the successive difference of the elements in a particular iteration is small enough.

\subsection{Staggered meshes}
\label{sec:mesh}

Let $\mathcal{T}_u$ be a triangulation of the two-dimensional domain $\Omega$ by a
set of triangles without hanging nodes.
We introduce the notation $\mathcal{F}_u$ to denote the
set of all edges in the triangulation $\mathcal{T}_u$ and
$\mathcal{F}_u^0$ to denote the subset of all interior edges in
$\mathcal{F}_u$ excluding those on the boundary of $\Omega$.
For each triangle in $\mathcal{T}_u$, we take an interior point $\nu$, denote the initial triangle by $\mathcal{S}(\nu)$, and divide $\mathcal{S}(\nu)$ into three triangles by
joining the point $\nu$ and the three vertices of $\mathcal{S}(\nu)$.
We also denote the set of all interior points $\nu$ by $\mathcal{N}$, the set of all new edges generated by the subdivision of triangles by $\mathcal{F}_p$, and the triangulation after subdivision by $\mathcal{T}$.
Note that the interior point $\nu$ of each triangle in $\mathcal{T}_u$ should be chosen such that the new triangulation $\mathcal{T}$ observes the shape regularity criterion.
In practice, we can simply choose $\nu$ as the centroid of the triangle.
Also,
$\mathcal{F} = \mathcal{F}_u\cup\mathcal{F}_p$
denotes the set of all edges of triangles in $\mathcal{T}$
and $\mathcal{F}^0 = \mathcal{F}^0_u \cup \mathcal{F}_p$ denotes the set of all interior edges of triangles in $\mathcal{T}$.
For each edge $e \in \mathcal{F}_u$, we let $\mathcal{R}(e)$ be the
union of the all triangles in the new triangulation $\mathcal{T}$ sharing the edge $e$.
{Figure}~\ref{fig_mesh} demonstrates these definitions. The edges $e \in \mathcal{F}_u$ are represented in solid lines and the $e \in \mathcal{F}_p$ are represented in dotted lines.

\setlength{\unitlength}{1mm}
\begin{figure}[ht!]
\vspace*{-1pc}
\begin{center}
\begin{picture}(80,50)(0,10)
\path(0,20)(30,50)(80,40) \path(30,50)(40,10) \path(0,20)(40,10)
\path(80,40)(40,10) \put(19,24){$\bullet$} \put(49,34){$\bullet$}
\dottedline(20,25)(30,50) \dottedline(20,25)(40,10)
\dottedline(20,25)(0,20) \dottedline(50,35)(30,50)
\dottedline(50,35)(40,10) \dottedline(50,35)(80,40)
\put(6,40){$\mathcal{S}(\nu_{1})$}
\put(60,47){$\mathcal{S}(\nu_{2})$} \put(28,23){$\mathcal{R}(e)$}
\put(36,34){$e$} \put(15,27){$\nu_1$} \put(52,38){$\nu_2$}
\end{picture}
\end{center}
\caption{An illustration of the staggered mesh in two dimensions.}\label{fig_mesh}
\end{figure}
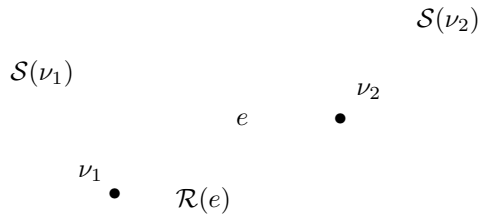

For each edge $e \in \mathcal{F}$, we will also define a unit normal vector $\bsym{n}_{e}$ in the following way.
If $e \in \mathcal{F} \setminus \mathcal{F}^0$ is a boundary edge, then we define $\bsym{n}_{e}$ as the outward unit normal vector of $e$ from $\Omega$.
If $e\in\mathcal{F}^0$ is an interior edge, then $\bsym{n}_{e}$ is fixed as one of the two possible unit normal vectors on $e$.
When it is clear that which edge we are considering, we omit the index $e$ and write the unit normal vector as $\bsym{n}$.

To end this section, we define the jumps in the following way:
for any edge $e \in \mathcal{F}$, denote one of the triangles in the refined triangulation $\mathcal{T}$, which contains $e$ by $\tau^+$, and denote the other triangle, if exists, by $\tau^-$.
The outward unit normal vectors on $e$ in $\tau^+$ and $\tau^-$ are denoted by $\bsym{n}^+$ and $\bsym{n}^-$, respectively.
Also, for any quantity $\phi$, the notations $\phi^\pm$ are defined on the edge $e$ by the values of $\phi \vert_{\tau^{\pm}}$ restricted on $e$.
Then, if $\phi$ is a scalar quantity, the notation $[\phi]$ over an edge $e$ defined as
\begin{equation}
[\phi]\vert_e:= (\bsym{n} \cdot \bsym{n}^+) \phi^+ + (\bsym{n} \cdot \bsym{n}^-) \phi^-.
\end{equation}
If $\bsym{\Phi}$ is a vector quantity, then the notation $[\bsym{\Phi} \cdot \bsym{n}]$ is similarly defined as
\begin{equation}
[\bsym{\Phi} \cdot \bsym{n}]|_e:=(\bsym{n}\cdot \bsym{n}^+) (\bsym{\Phi}^+ \cdot \bsym{n})
+ (\bsym{n}\cdot \bsym{n}^-) (\bsym{\Phi}^- \cdot \bsym{n}).
\end{equation}

\subsection{SDG finite element spaces}
\label{sec:SDGspace}

We will define the SDG finite element spaces.
Let $k \geq 0$ be a non-negative integer.
Let $\tau\in\mathcal{T}$ and $e\in\mathcal{F}$.
We define $P^{k}(\tau)$ and $P^k(e)$ as the space of polynomials whose order is not greater than $k$ on $\tau$ and $e$, respectively.
We will also define norms on the spaces.
We use the standard notations $ \| \cdot \|_{0, \Omega}$ to denote the standard $L^2$ norm on $\Omega$ and $ \| \cdot \|_{0, e}$ to denote the $L^2$ norm on an edge $e$.

First, we define the following locally $H^{1}(\Omega)$-conforming finite element space for velocity:
\begin{equation}
U^{h}
= \{ v \: : \: v|_{\tau} \in P^{k}(\tau); \; \tau\in\mathcal{T}; \;
v \;\text{is
continuous over} \; e \in \mathcal{F}_{u}^{0}; \;
v|_{\partial \Omega}=0 \}.
\label{eq:Uh}
\end{equation}
Note that for any $v\in U^h$, we have $v|_{\mathcal{R}(e)} \in
H^1(\mathcal{R}(e))$ for each edge $e \in \mathcal{F}_u$.
We define the following discrete $L^2$-norm $\| \cdot \|_X$ and discrete $H^1$-norm $\| \cdot \|_Z$ on the space $U^h$:
\begin{equation}
\begin{split}
\| v \|_X & = \left( \| v \|_{0, \Omega}^2 + \sum_{e \in \mathcal{F}_u^0} h_e \| v \|_{0,e}^2 \right)^\frac{1}{2}, \\
\| v \|_Z & = \left( \| \nabla_h v \|_{0, \Omega}^2 + \sum_{e \in \mathcal{F}_p} h_e^{-1} \| [v] \|_{0,e}^2 \right)^\frac{1}{2},
\end{split}
\label{eq:Unorm}
\end{equation}
\revr{where $\nabla_h$ denotes the gradient operator applied piecewise on the given triangulation $\mathcal{T}$.}
For $\bsym{v} = (v_1, v_2) \in [U^h]^2$, we also define an energy norm
\begin{equation}
\| \bsym{v} \|_h = ( \| v_1 \|_Z^2 + \| v_2 \|_Z^2)^\frac{1}{2}.
\end{equation}

Next, we define the following locally $H(\text{div};\Omega)$-conforming finite element space for velocity gradients:
\begin{equation}
W^{h}
= \{ \bsym{\Psi} \: : \: \bsym{\Psi}|_{\tau} \in P^{k}(\tau)^{2}; \: \tau\in\mathcal{T}; \:
\bsym{\Psi} \cdot \bsym{n} \: \text{is continuous over} \:
e\in \mathcal{F}_{p} \}.
\label{eq:Wh}
\end{equation}
Note that for any $\bsym{\Psi} \in W^h$, we have
$\bsym{\Psi} \vert_{\mathcal{S}(\nu)} \in H(\text{div};\mathcal{S}(\nu))$
for each $\nu \in \mathcal{N}$.
We define the following discrete $L^2$-norm $\| \cdot \|_{X'}$ and discrete $H(\text{div};\Omega)$-norm $\| \cdot \|_{Z'}$ on the space $W^h$:
\begin{equation}
\begin{split}
\| \bsym{\Psi} \|_{X'} & = \left( \| \bsym{\Psi} \|_{0, \Omega}^2 + \sum_{e \in \mathcal{F}_p} h_e \| \bsym{\Psi} \cdot \bsym{n} \|_{0,e}^2 \right)^\frac{1}{2}, \\
\| \bsym{\Psi} \|_{Z'} & = \left( \| \text{div}_h \bsym{\Psi} \|_{0, \Omega}^2 + \sum_{e \in \mathcal{F}_u^0} h_e^{-1} \| [\bsym{\Psi} \cdot \bsym{n}] \|_{0,e}^2 \right)^\frac{1}{2}.
\end{split}
\label{eq:Wnorm}
\end{equation}
\revr{Here $\text{div}_h$ denotes the divergence operator applied piecewise on the given triangulation $\mathcal{T}$.}

We also define the following locally $H^{1}(\Omega)$-conforming finite element space for pressure:
\begin{equation}
P^h
=\{q \: : \: q|_{\tau} \in P^{k}(\tau); \; \tau\in\mathcal{T}; \;
q \;\text{is
continuous over} \; e \in \mathcal{F}_{p};\; \int_\Omega q \;dx=0\}.
\label{eq:Ph}
\end{equation}
We define the following discrete $L^2$-norm $\| \cdot \|_P$ on the space $P^h$:
\begin{equation}
\| q \|_P = \left( \| q \|_{0, \Omega}^2 + \sum_{e \in \mathcal{F}_p} h_e \| q \|_{0,e}^2 \right)^\frac{1}{2}. \\
\label{eq:Pnorm}
\end{equation}

Finally, we define a finite element space for the Eulerian coordinates of the immersed boundary. Suppose we have a partition of the interval $D = [0,L]$ in the Lagrangian coordinate system:
\begin{equation}
0 = s_0 < s_1 < s_2 < \ldots < s_m = L.
\label{eq:Lpart}
\end{equation}
We denote the subintervals by $J_i = (s_{i-1}, s_i)$ and define the following space:
\begin{equation}
S^h = \{ \bsym{Y} \: : \: \bsym{Y}|_{J_i} \in P^1(J_i); \; 1 \leq i \leq m; \; \bsym{Y} \text{ is continuous at } s_i; \; \bsym{Y}(0) = \bsym{Y}(L)\}.
\label{eq:Sh}
\end{equation}
For any $\bsym{Y} \in S^h$, $\bsym{Y}$ is an $m$-sided polygon with vertices $\bsym{Y}(s_i)$.

\subsection{SDG spatial discretization}
\label{sec:sdg}

In view of \eqref{eq:IB_NS3}, at each time step $n$ and each iteration $m$, one needs to solve the system of linear PDEs:
\begin{equation}
\begin{split}
\alpha \bsym{u} - \mu \bsym{\Delta} \bsym{u} + \rho \bsym{V} \cdot \nabla \bsym{u} + \nabla p & = \bsym{F} \; \mbox{ in }\Omega,\\
\text{div} \, \bsym{u} & = 0 \; \mbox{ in }\Omega,\\
\bsym{u} & = 0 \; \mbox { on } \partial \Omega.\\
\end{split}
\label{eq:IB_NS4}
\end{equation}
We introduce the auxiliary variables
\begin{equation}
\begin{split}
\bsym{w} & = \sqrt{\mu}  \, \nabla u_1 - \frac{\rho}{2\sqrt{\mu}} u_1 \bsym{V},\\
\bsym{z} & = \sqrt{\mu} \, \nabla u_2 - \frac{\rho}{2\sqrt{\mu}} u_2 \bsym{V},\\
\widetilde{\bsym{w}} & = u_1 \bsym{V},\\
\widetilde{\bsym{z}} & = u_2 \bsym{V}.
\end{split}
\label{eq:newvar}
\end{equation}
Then \eqref{eq:IB_NS4} can be reformulated as a system of first-order linear PDEs:
\begin{equation}
\begin{split}
\alpha u_1 -\sqrt{\mu} \, \text{div}\, \bsym{w} + \frac{\rho}{2\sqrt{\mu}} \bsym{V}\cdot\bsym{w} + \frac{\rho^2}{4 \mu} \bsym{V}\cdot\widetilde{\bsym{w}} + p_x & = F_1 \; \mbox{ in }\Omega, \\
\alpha u_2 - \sqrt{\mu} \, \text{div}\, \bsym{z} + \frac{\rho}{2\sqrt{\mu}} \bsym{V}\cdot\bsym{z} + \frac{\rho^2}{4 \mu} \bsym{V}\cdot\widetilde{\bsym{z}} + p_y & = F_2 \; \mbox{ in }\Omega, \\
\text{div} \, \bsym{u} & = 0 \; \mbox{ in }\Omega,\\
\bsym{u} & = 0 \; \mbox{ on }\partial \Omega.
\end{split}
\label{eq:IB_NS5}
\end{equation}

We will derive the discrete problem in our SDG formulation
starting from the system of first order equations in \eqref{eq:newvar} and \eqref{eq:IB_NS5}.

Multiplying the first equation of \eqref{eq:newvar} by $\bsym{\Psi}_1 \in W^h$
and integrating over $\mathcal{S}(\nu)$ for $\nu \in \mathcal{N}$, we obtain
\begin{equation}
\int_{\mathcal{S}(\nu)} \bsym{w} \cdot \bsym{\Psi}_1 \; dx = - \sqrt{\mu} \int_{\mathcal{S}(\nu)} u_1 \text{div } \bsym{\Psi}_1 \; dx
+ \sqrt{\mu} \int_{\partial \mathcal{S}(\nu)} u_1 \bsym{\Psi}_1 \cdot \bsym{n} \; d\sigma
- \frac{\rho}{2\sqrt{\mu}} \int_{\mathcal{S}(\nu)} \widetilde{\bsym{w}} \cdot \bsym{\Psi}_1 \; dx.
\label{eq:w1}
\end{equation}
Similarly, multiplying the second equation of \eqref{eq:newvar} by $\bsym{\Psi}_2 \in W^h$
and integrating over $\mathcal{S}(\nu)$ for $\nu \in \mathcal{N}$, we have
\begin{equation}
\int_{\mathcal{S}(\nu)} \bsym{z} \cdot \bsym{\Psi}_2 \; dx = - \sqrt{\mu} \int_{\mathcal{S}(\nu)} u_2 \text{div } \bsym{\Psi}_2 \; dx
+ \sqrt{\mu} \int_{\partial \mathcal{S}(\nu)} u_2 \bsym{\Psi}_2 \cdot \bsym{n} \; d\sigma
- \frac{\rho}{2\sqrt{\mu}} \int_{\mathcal{S}(\nu)} \widetilde{\bsym{z}} \cdot \bsym{\Psi}_2 \; dx.
\label{eq:z1}
\end{equation}

Multiplying the third equation of \eqref{eq:newvar} by $\bsym{\Psi}_3 \in W^h$
and integrating over $\mathcal{S}(\nu)$ for $\nu \in \mathcal{N}$, we have
\begin{equation}
\int_{\mathcal{S}(\nu)} \widetilde{\bsym{w}} \cdot \bsym{\Psi}_3 \; dx = \int_{\mathcal{S}(\nu)} u_1 \bsym{V} \cdot \bsym{\Psi}_3 \; dx.
\label{eq:tw1}
\end{equation}
Similarly, multiplying the fourth equation of \eqref{eq:newvar} by $\bsym{\Psi}_4 \in W^h$
and integrating over $\mathcal{S}(\nu)$ for $\nu \in \mathcal{N}$, we have
\begin{equation}
\int_{\mathcal{S}(\nu)} \widetilde{\bsym{z}} \cdot \bsym{\Psi}_4 \; dx = \int_{\mathcal{S}(\nu)} u_2 \bsym{V} \cdot \bsym{\Psi}_4 \; dx.
\label{eq:tz1}
\end{equation}

Multiplying the first equation of \eqref{eq:IB_NS5} by $v_1 \in U^h$
and integrating over $\mathcal{R}(e)$ for $e \in \mathcal{F}_u^0$, we have
\begin{equation}
\begin{split}
& \alpha \int_{\mathcal{R}(e)} u_1 \, v_1 \; dx + \sqrt{\mu} \int_{\mathcal{R}(e)} \bsym{w} \cdot \nabla v_1\; dx - \sqrt{\mu} \int_{\partial \mathcal{R}(e)} ( \bsym{w}\cdot \bsym{n}) v_1 \;d\sigma
+ \frac{\rho}{2\sqrt{\mu}} \int_{\mathcal{R}(e)} \bsym{V}\cdot \bsym{w} \, v_1 \; dx \\
& + \frac{\rho^2}{4 \mu} \int_{\mathcal{R}(e)} \bsym{V}\cdot \widetilde{\bsym{w}} \, v_1 \; dx
- \int_{\mathcal{R}(e)} p (v_1)_x
+ \int_{\partial \mathcal{R}(e)} p v_1 n_1 \;d\sigma = \int_{\mathcal{R}(e)} F_1 v_1\; dx.
\end{split}
\label{eq:div_w}
\end{equation}
Similarly, multiplying the second equation of \eqref{eq:IB_NS5} by $v_2 \in U^h$
and integrating over $\mathcal{R}(e)$ for $e \in \mathcal{F}_u^0$, we have
\begin{equation}
\begin{split}
& \alpha \int_{\mathcal{R}(e)} u_2 \, v_2 \; dx + \sqrt{\mu} \int_{\mathcal{R}(e)} \bsym{z} \cdot \nabla v_2 \; dx - \sqrt{\mu} \int_{\partial \mathcal{R}(e)} (\bsym{z}\cdot \bsym{n}) v_2\;d\sigma
+ \frac{\rho}{2\sqrt{\mu}} \int_{\mathcal{R}(e)} \bsym{V}\cdot \bsym{z} \, v_2 \; dx \\
& + \frac{\rho^2}{4 \mu} \int_{\mathcal{R}(e)} \bsym{V}\cdot \widetilde{\bsym{z}} \, v_2 \; dx
- \int_{\mathcal{R}(e)} p (v_2)_y \; dx
+ \int_{\partial \mathcal{R}(e)} p v_2 n_2 \;d\sigma  = \int_{\mathcal{R}(e)} F_2 v_2\; dx.
\end{split}
\label{eq:div_z}
\end{equation}

Finally, multiplying the third equation of \eqref{eq:IB_NS5} by $q\in P^h$,
and integrating over $\mathcal{S}(\nu)$ for $\nu \in \mathcal{N}$, we have
\begin{equation}
-\int_{\mathcal{S}(\nu)} \bsym{u} \cdot \nabla q \; dx+ \int_{\partial\mathcal{S}(\nu)} (\bsym{u}\cdot \bsym{n}) q \; d\sigma= 0.
\label{eq:div_0}
\end{equation}

Summing those equations in \eqref{eq:w1}--\eqref{eq:div_0} over all $\mathcal{R}(e)$ and $\mathcal{S}(\nu)$,
our staggered discontinuous Galerkin method for
\eqref{eq:IB_NS4} is obtained: find
$(\bsym{u}_h,\bsym{w}_h,\bsym{z}_h,\widetilde{\bsym{w}}_h,\widetilde{\bsym{z}}_h,p_h)\in [U^h]^2\times [W^h]^4 \times P^h$ such that
for any $\bsym{v} = (v_1,v_2) \in [U^h]^2, \bsym{\Psi}_1, \bsym{\Psi}_2, \bsym{\Psi}_3, \bsym{\Psi}_4 \in W^h, q \in P^h$, we have
\begin{equation}
\begin{split}
\alpha (\bsym{u}_h,\bsym{v})_{0,\Omega}  + \sqrt{\mu} B_h(\bsym{w}_h,v_1)+ \sqrt{\mu} B_h(\bsym{z}_h,v_2)+ & \\
\frac{\rho}{2\sqrt{\mu}}R_h\left(\bsym{w}_h + \frac{\rho}{2\sqrt{\mu}}\widetilde{\bsym{w}}_h,v_1\right) + \frac{\rho}{2\sqrt{\mu}}R_h\left(\bsym{z}_h + \frac{\rho}{2\sqrt{\mu}}\widetilde{\bsym{z}}_h,v_2\right)
+  b_h^*(p_h, \bsym{v})& = (\bsym{F},\bsym{v})_{0,\Omega},\\
\sqrt{\mu} B_h^*(u_{h,1},\bsym{\Psi}_1) - \frac{\rho}{2\sqrt{\mu}} (\widetilde{\bsym{w}}_h,\bsym{\Psi}_1)_{0,\Omega} & = (\bsym{w}_h,\bsym{\Psi}_1)_{0,\Omega},\\
\sqrt{\mu} B_h^*(u_{h,2},\bsym{\Psi}_2) - \frac{\rho}{2\sqrt{\mu}} (\widetilde{\bsym{z}}_h,\bsym{\Psi}_2)_{0,\Omega} & = (\bsym{z}_h,\bsym{\Psi}_2)_{0,\Omega},\\
R_h^*(u_{h,1},\bsym{\Psi}_3) &=(\widetilde{\bsym{w}}_h,\bsym{\Psi}_3)_{0,\Omega},\\
R_h^*(u_{h,2},\bsym{\Psi}_4) &=(\widetilde{\bsym{z}}_h,\bsym{\Psi}_4)_{0,\Omega},\\
b_h(\bsym{u}_h,q)&=0,
\end{split}
\label{eq:NS_discrete}
\end{equation}
where bilinear forms $B_h(\bsym{\Psi},v)$ and $B^*_h(v,\bsym{\Psi})$ are defined as
\begin{equation}
\begin{split}
B_{h}(\bsym{\Psi},v) &=  \int_{\Omega} \bsym{\Psi} \cdot \nabla_h v
\; dx
- \sum_{e\in\mathcal{F}_{p}}\int_{e} \bsym{\Psi} \cdot \bsym{n} \: [v] \; d\sigma,  \\
B^{*}_{h}(v,\bsym{\Psi}) &= -\int_{\Omega}  v \: \text{div}_h\,
\bsym{\Psi} \;dx + \sum_{e\in\mathcal{F}_{u}^{0}}\int_{e} v \:
[\bsym{\Psi} \cdot \bsym{n}] \; d\sigma,
\end{split}
\label{eq:bilinear1}
\end{equation}
and the bilinear forms $b_h^*(q,\bsym{v})$ and $b_h(\bsym{v},q)$
as 
\begin{equation}
\begin{split}
b_h^*(q,\bsym{v}) &= -\int_\Omega q \,\text{div}_h\, \bsym{v} \;dx +
\sum_{e\in\mathcal{F}_p}\int_e q [\bsym{v}\cdot\bsym{n}] \;d\sigma,\\
b_h(\bsym{v},q)&=\int_{\Omega} \bsym{v} \cdot \nabla_h q \;dx
-\sum_{e \in \mathcal{F}_u^0} \int_e \bsym{v} \cdot \bsym{n} [q] \; d\sigma.
\end{split}
\label{eq:bilinear2}
\end{equation}
The bilinear forms $R_h(\bsym{\Psi},v)$ and $R_h^*(v,\bsym{\Psi})$ are also defined as
\begin{equation}
\begin{split}
R_h(\bsym{\Psi},v) & = \int_\Omega (\bsym{V} \cdot \bsym{\Psi}) \, v \; dx,\\
R_h^*(v,\bsym{\Psi}) & = \int_\Omega v \, (\bsym{V} \cdot \bsym{\Psi}) \; dx.
\end{split}
\label{eq:bilinear3}
\end{equation}
Moreover, $(\cdot,\cdot)_{0,\Omega}$ denotes the standard $L^2(\Omega)$ inner product.

By \cite{newdg1}, the two bilinear forms in \eqref{eq:bilinear1} satisfy the adjoint relation
\begin{equation}
B_h(\bsym{\Psi},v) = B_h^*(v,\bsym{\Psi})
\label{eq:adjoint-B}
\end{equation}
for all $v\in U_h$ and $\bsym{\Psi} \in W^h$.
The bilinear forms $B_h$ and $B_h^*$ are also continuous with respect to suitable discrete norms
\begin{equation}
\begin{split}
\vert B_h(\bsym{\Psi}, v) \vert & \leq \| \bsym{\Psi} \|_{X'} \| v \|_Z, \\
\vert B_h^*(v, \bsym{\Psi}) \vert & \leq  \| v \|_X \| \bsym{\Psi} \|_{Z'}, \\
\end{split}
\label{eq:cont-B}
\end{equation}
for all $v\in U_h$ and $\bsym{\Psi} \in W^h$.
Moreover, the bilinear forms $B_h$ and $B_h^*$ satisfy a pair of inf-sup conditions: there exists constants $\beta_1$ and $\beta_2$, independent of $h$, such that
\begin{equation}
\begin{split}
\inf_{v \in U^h \setminus \{0\}} \sup_{\bsym{\Psi} \in W^h \setminus \{\bsym{0}\}}
\frac{B_h(\bsym{\Psi}, v)}{ \| \bsym{\Psi} \|_{X'} \| v \|_Z} & \geq \beta_1, \\
\inf_{\bsym{\Psi} \in W^h \setminus \{\bsym{0}\}} \sup_{v \in U^h \setminus \{0\}}
\frac{B_h^*(v, \bsym{\Psi})}{\| v \|_X \| \bsym{\Psi} \|_{Z'}} & \geq \beta_2. \\
\end{split}
\label{eq:inf-sup-B}
\end{equation}

By \cite{kcl-2013-Stokes-SDG}, the two bilinear forms in \eqref{eq:bilinear2} satisfy the adjoint relation
\begin{equation}\label{eq:adjoint-b}
b_h^*(q,\bsym{v})=b_h(\bsym{v},q)
\end{equation}
for all $q \in P_h$ and $\bsym{v}\in [U_h]^2$.
The bilinear form $b_h$ is also continuous: there exists a constant $C_b$ such that
\begin{equation}
\vert b_h(\bsym{v},q) \vert \leq C_b \| \bsym{v} \|_h \| q \|_{0, \Omega},
\label{eq:cont-b}
\end{equation}
for all $q \in P_h$ and $\bsym{v}\in [U_h]^2$.
Moreover, the bilinear form $b_h$ satisfies an inf-sup condition: there exists a constant $\gamma$, independent of $h$, such that
\begin{equation}
\inf_{q \in P^h \setminus \{0\}} \sup_{\bsym{v} \in [U^h]^2 \setminus \{\bsym{0}\}}
\frac{b_h(\bsym{v},q)}{\| \bsym{v} \|_h \| q \|_{0,\Omega}} \geq \gamma.
\label{eq:inf-sup-b}
\end{equation}

Finally, the two bilinear forms in (\ref{eq:bilinear3}) satisfy
\begin{equation}\label{eq:adjoint-r}
R_h^*(v,\bsym{\Psi})=R_h(\bsym{\Psi},v)
\end{equation}
for all $v \in U_h$ and $\bsym{\Psi} \in W^h$.

\subsection{Linear system}
\label{sec:sys}

In this section,
we derive the linear system resulting from \eqref{eq:NS_discrete}. We denote the corresponding matrix representation of the bilinear forms $B_h$, $b_h$ and $R_h$ by $B$, $C$ and $R$, respectively. Then by the adjoint properties, the matrix representation of the bilinear forms $B^*_h$, $b^*_h$ and $R^*_h$ are given by $B^\top$, $C^\top$ and $R^\top$, respectively. Also, the notations for the finite element solutions would be abused to denote their corresponding vector representations.

The second and the third equations of \eqref{eq:NS_discrete} can be written as
\begin{equation}
\begin{split}
\sqrt{\mu} B^\top u_{h,1} - \frac{\rho}{2\sqrt{\mu}} M\bsym{\widetilde{w}}_h & = M\bsym{w}_h,\\
\sqrt{\mu} B^\top u_{h,2} - \frac{\rho}{2\sqrt{\mu}} M\bsym{\widetilde{z}}_h & = M\bsym{z}_h,
\end{split}
\label{eq:matrix_eq1}
\end{equation}
where $M$ is the mass matrix for the space $W^h$.
Similarly, the fourth and the fifth equations of \eqref{eq:NS_discrete} can be written as
\begin{equation}
\begin{split}
R^\top u_{h,1} & = M\bsym{\widetilde{w}}_h,\\
R^\top u_{h,2} & = M\bsym{\widetilde{z}}_h.
\end{split}
\label{eq:matrix_eq2}
\end{equation}
Lastly, the first and the last equations of \eqref{eq:NS_discrete} can be written as
\begin{equation}
\begin{split}
\alpha
\begin{pmatrix}
\widetilde{M} u_{h,1}\\
\widetilde{M} u_{h,2}
\end{pmatrix}
+ \sqrt{\mu}
\begin{pmatrix}
B\bsym{w}_h\\
B\bsym{z}_h
\end{pmatrix}
 + \frac{\rho}{2\sqrt{\mu}}
\begin{pmatrix}
R(\bsym{w}_h + \frac{\rho}{2\sqrt{\mu}}\bsym{\widetilde{w}}_h)\\
R(\bsym{z}_h + \frac{\rho}{2\sqrt{\mu}}\bsym{\widetilde{z}}_h)
\end{pmatrix}
 + C^\top p_h &=
\begin{pmatrix}
F_{h,1} \\
F_{h,2}
\end{pmatrix},\\
C\bsym{u}_h & = 0,
\end{split}
\label{eq:matrix_eq3}
\end{equation}
where $\widetilde{M}$ is the mass matrix for the space $U^h$.
We can now obtain a linear system with the unknowns $\bsym{w}_h, \bsym{z}_h, \bsym{\widetilde{w}}_h, \bsym{\widetilde{z}}_h$ eliminated. Combining \eqref{eq:matrix_eq1} and \eqref{eq:matrix_eq2}, we have
\begin{equation}
\begin{split}
\bsym{w}_h & = M^{-1} (\sqrt{\mu} B^\top u_{h,1} - \frac{\rho}{2\sqrt{\mu}} R^\top u_{h,1}),\\
\bsym{z}_h & = M^{-1} (\sqrt{\mu} B^\top u_{h,2} - \frac{\rho}{2\sqrt{\mu}} R^\top u_{h,2}),\\
\bsym{\widetilde{w}}_h & = M^{-1} R^\top u_{h,1},\\
\bsym{\widetilde{z}}_h & = M^{-1} R^\top u_{h,2}.
\end{split}
\label{eq:matrix_eq4}
\end{equation}
We note that the elimination can be done by solving small problems in each $\mathcal{S}(\nu)$ since $M$ is a block diagonal matrix with each block corresponding to the mass matrix
of $W^h|_{\mathcal{S}(\nu)}$.

We further introduce the notations
\begin{equation}
\begin{split}
\bsym{\Delta}_h &= -BM^{-1}B^\top,\\
\bsym{V} \cdot \nabla_h &= -\frac{1}{2} BM^{-1}R^\top + \frac{1}{2} RM^{-1}B^\top, \\
A &= \alpha \widetilde{M} -\mu \bsym{\Delta}_h + \rho \bsym{V} \cdot \nabla_h.
\end{split}
\label{eq:mat_def}
\end{equation}
We note that the negative of the discrete Laplacian operator $-\Delta_h$ is symmetric and positive-definite, and the discrete convection operator $\bsym{V} \cdot \nabla_h$ is skew-symmetric.
Combining \eqref{eq:matrix_eq3} and \eqref{eq:matrix_eq4}, the algebraic system of the discrete problem \eqref{eq:NS_discrete} can then be reduced to
\begin{equation}
\left(\begin{array}{ccc}
A & 0 & \multirow{2}{*}{$C^\top$}\\
0 & A & \\
\multicolumn{2}{c}{C} & 0
\end{array}\right)
\begin{pmatrix}
u_{h,1}\\
u_{h,2}\\
p_h
\end{pmatrix}
=
\begin{pmatrix}
F_{h,1}\\
F_{h,2}\\
0
\end{pmatrix}
\label{eq:NS_system}
\end{equation}
and the above system is solved for the unknowns $(u_{h,1},\,u_{h,2},\,p_h)$.

\subsection{Postprocessing}
\label{sec:post}

In this section,
we present a postprocessing technique for the velocity, which was
introduced in \cite{hdg}.
In our case, we perform the postprocessing on each $\mathcal{S}(\nu)$
to obtain a divergence-free velocity with a higher convergence rate.

Let $(\bsym{u}_h,\bsym{w}_h,\bsym{z}_h,\widetilde{\bsym{w}}_h,\widetilde{\bsym{z}}_h,p_h)\in [U^h]^2\times [W^h]^4 \times P^h$ be the solution of \eqref{eq:NS_discrete}. We introduce the notations
\begin{equation}
\begin{split}
\widehat{\bsym{w}}_h & = \frac{1}{\sqrt{\mu}} \left(\bsym{w}_h + \frac{\rho}{2\sqrt{\mu}} \widetilde{\bsym{w}}_h\right),\\
\widehat{\bsym{z}}_h & = \frac{1}{\sqrt{\mu}} \left(\bsym{z}_h + \frac{\rho}{2\sqrt{\mu}} \widetilde{\bsym{z}}_h\right),\\
L_h & = \begin{pmatrix} \widehat{\bsym{w}}_h^T \\ \widehat{\bsym{z}}_h^T \end{pmatrix}.
\end{split}
\label{eq:newvar2}
\end{equation}
Then $L_h$ is an approximation for the  matrix $L$ of $\nabla \bsym{u}$.

Let $\bsym{u}_h^{\star} \in P^{k+1}(\mathcal{S}(\nu))^2$ be the post-processed velocity.
For every edge $e\in\partial\mathcal{S}(\nu)$, $\bsym{u}_h^\star$ satisfies
\begin{equation}
\label{eq:post1}
\int_e (\bsym{u}_h^{\star} - \bsym{u}_h) \cdot \bsym{n} \, v \; d\sigma = 0, \quad\forall v\in P^k(e)
\end{equation}
and
\begin{equation}
\label{eq:post2}
\int_e \Big(  (\bsym{n}\times \nabla) (\bsym{u}_h^{\star}) - \bsym{n}\times (\{ L_h^T \} \bsym{n}) \Big)  \, (\bsym{n}\times \nabla) v \; d\sigma = 0, \quad\forall v\in P^k(e).
\end{equation}
In the two-dimensional case, we have $\bsym{n}\times \nabla = n_2 \partial_1 - n_1 \partial_2$,
$\bsym{n}\times \bsym{a} = n_1 a_2 - n_2 a_1$, and $\nabla \times \bsym{a} = \partial_1 a_2 - \partial_2 a_1$.
In addition, $\bsym{u}_h^\star$ satisfies
\begin{equation}
\label{eq:post3}
\int_{\mathcal{S}(\nu)} ( \bsym{u}_h^{\star} - \bsym{u}_h) \cdot \nabla v \; dx = 0, \quad \forall v\in P^k(\mathcal{S}(\nu))
\end{equation}
and
\begin{equation}
\label{eq:post4}
\int_{\mathcal{S}(\nu)} (\nabla\times \bsym{u}_h^{\star} - \mathcal{L}_h ) \, v \mathcal{B} \; dx = 0, \quad \forall v\in P^{k-1}(\mathcal{S}(\nu))
\end{equation}
where $\mathcal{L}_h = (L_h)_{21} - (L_h)_{12}$ and $\mathcal{B}$ is the bubble function,
defined by the product of barycentric coordinates of vertices of $\mathcal{S}(\nu)$.

We solve \eqref{eq:post1}-\eqref{eq:post4} to obtain the post-processed velocity $\bsym{u}_h^{\star}$.
In \cite{sdg-ns2}, it is shown that $\bsym{u}_h^{\star}$ is exactly divergence-free.

The pointwise divergence-free property of post-processed velocity is vital in SDG-IBM. First, in the sequence of Picard fixed point iterations in \eqref{eq:IB_NS3}, the velocity field $\bsym{V}_m^n$ is chosen to be the post-processed velocity from $\bsym{u}_{m-1,h}^{n}$. Second, in the full discretization of \eqref{eq:IB_EL2}, the Lagrangian markers are driven by the post-processed velocity of the fixed point velocity field at a certain time step. More details will be explained
in Section~\ref{sec:EL}.

\subsection{Discretization of source term}
\label{sec:source}
In Section~\ref{sec:post}, we have discussed the linearization of the convection term by the post-processed velocity in each iteration. In Sections \ref{sec:mesh} -- \ref{sec:sys}, we have discussed the SDG method for solving the linearized equation in each iteration, given a particular source term. To complete the discussion on our method for solving \eqref{eq:IB_NS3}, it remains to discuss the spatial discretization of the source term given by \eqref{eq:IB_source2} at each time level.

\begin{figure}[ht!]
\centering
\includegraphics[width=3.0in]{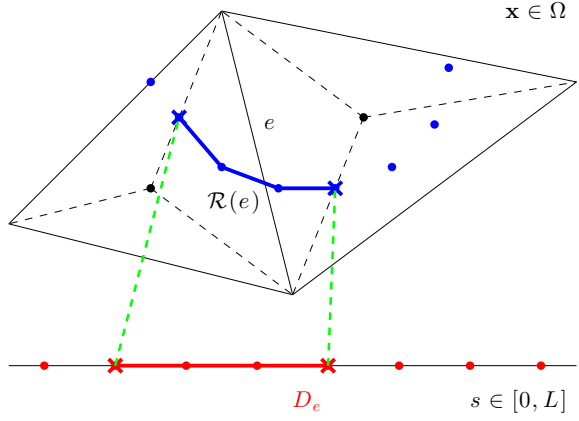}
\caption{Illustration of the preimage set $D_e$ in the Lagrangian coordinate system.}
\label{fig:var}
\end{figure}

We will start with a variational equation with local test functions for the continuous problem. Let $e \in \mathcal{F}_u^0$. Suppose $\bsym{v} \in H^1(\mathcal{R}(e))$. Let $D_e$ be the preimage set of $\mathcal{R}(e)$ under $\bsym{X}(\cdot, t)$, i.e.
\begin{equation}
D_e = \{s \in [0,L]: \bsym{X}(s,t) \in \mathcal{R}(e)\}.
\end{equation}
Figure \ref{fig:var} illustrates the preimage set $D_e$ in the Lagrangian coordinate system.

Without loss of generality assume $D_e$ is connected.
Similar to \cite{boffi03}, we have the following the variational equation:
\begin{definition}
Suppose $\bsym{X}(\cdot, t) \in W^{1,\infty}([0,L])$ for $t \in [0,T]$ and $\bsym{f} \in L^2([0,L] \times (0,T))$. Then for $e \in \mathcal{F}_u^0$, for $t \in (0,T)$, the force density $\bsym{F}(t)$ is defined as follows:
\begin{equation}
(\bsym{F}(t), \bsym{v})_{0,\mathcal{R}(e)} = \int_{D_e} \bsym{f}(s,t) \bsym{v}(\bsym{X}(s,t)) ds \quad \text{ for } \bsym{v} \in [H^1(\mathcal{R}(e))]^2.
\label{eq:IB_source3}
\end{equation}
\end{definition}

In particular, in our simple model, for $\bsym{v} \in [H^1(\mathcal{R}(e))]^2$, substituting \eqref{eq:sim_mod} into \eqref{eq:IB_source3} and using integration by parts over $D_e$, we have
\begin{equation}
\begin{split}
(\bsym{F}(t), \bsym{v})_{0,\Omega}
& = \sum_{e \in \mathcal{F}_u^0} \int_{D_e} \kappa \frac{\partial^2 \bsym{X}(s,t)}{\partial s^2} \bsym{v}(\bsym{X}(s,t)) ds\\
& = \sum_{e \in \mathcal{F}_u^0}  \left(\kappa \frac{\partial \bsym{X}(s,t)}{\partial s} \bsym{v}(\bsym{X}(s,t)) \Big\vert_{\partial D_e} - \int_{D_e} \kappa \frac{\partial \bsym{X}(s,t)}{\partial s} \frac{\partial \bsym{v}(\bsym{X}(s,t))}{\partial s} ds\right).
\end{split}
\label{eq:IB_source4}
\end{equation}
\revb{In our approach, we use \eqref{eq:IB_source4} for the force exerted on the fluid
by the immersed structure. For simplicity, we consider $X(s,t)$ as a piecewise linear function on the
partition, $s_0=0<s_1< \cdots <s_m=L$.}

Next, we consider a corresponding full discretization of \eqref{eq:IB_source4} using forward-Euler time-stepping as described in \eqref{eq:IB_source2}. Note that for $\bsym{X}_h \in S^h$, by construction we have
\begin{equation}
\frac{\partial \bsym{X}_h}{\partial s}(s) = \frac{\partial \bsym{X}_h}{\partial s}(s_{i-\frac{1}{2}}) \quad \mbox{ for } s \in J_i,\\
\label{eq:mid}
\end{equation}
where $s_{i-\frac{1}{2}} = \frac{1}{2}(s_{i-1} + s_i)$ for $i = 1, 2, \ldots, m$.
Using \eqref{eq:IB_source4}, the source term can be discretized by: for $\bsym{v} \in [U^h]^2$,
\begin{equation}
\begin{split}
(\bsym{F}^n_h, \bsym{v})_{0,\Omega}
& = \sum_{e \in \mathcal{F}_u^0} \left(\kappa \frac{\partial \bsym{X}^{n-1}_h}{\partial s} \bsym{v}(\bsym{X}^{n-1}_h) \Big\vert_{\partial D_e} - \int_{D_e} \kappa \frac{\partial \bsym{X}^{n-1}_h}{\partial s} \frac{\partial \bsym{v}(\bsym{X}^{n-1}_h)}{\partial s} ds\right)\\
& = \sum_{e \in \mathcal{F}_u^0} \left(\kappa \frac{\partial \bsym{X}^{n-1}_h}{\partial s} \bsym{v}(\bsym{X}^{n-1}_h) \Big\vert_{\partial D_e} - \sum_{i=1}^m \int_{D_e \cap J_i} \kappa \frac{\partial \bsym{X}^{n-1}_h}{\partial s} \frac{\partial \bsym{v}(\bsym{X}^{n-1}_h)}{\partial s} ds\right)\\
& = \sum_{e \in \mathcal{F}_u^0} \left(\kappa \frac{\partial \bsym{X}^{n-1}_h}{\partial s} \bsym{v}(\bsym{X}^{n-1}_h) \Big\vert_{\partial D_e} - \sum_{i=1}^m \kappa \frac{\partial \bsym{X}_h^{n-1}}{\partial s}(s_{i-\frac{1}{2}}) \int_{D_e \cap J_i} \frac{\partial \bsym{v}(\bsym{X}^{n-1}_h)}{\partial s} ds\right)\\
& = - \sum_{i=1}^m \kappa \frac{\partial \bsym{X}_h^{n-1}}{\partial s}(s_{i-\frac{1}{2}}) \Big(\bsym{v}(\bsym{X}^{n-1}_h(s_i)) - \bsym{v}(\bsym{X}^{n-1}_h(s_{i-1}))\Big)\\
& = \sum_{i=1}^m \kappa  \left(\frac{\partial \bsym{X}_h^{n-1}}{\partial s}(s_{i+\frac{1}{2}}) - \frac{\partial \bsym{X}_h^{n-1}}{\partial s}(s_{i-\frac{1}{2}}) \right)\bsym{v}(\bsym{X}^{n-1}_h(s_i)).\\
\end{split}
\label{eq:source_discrete}
\end{equation}
For the sake of simplifying notations, we use periodic indices, i.e. $s_{m+r} = s_r$. We remark that the variational equation \eqref{eq:IB_source3} for source term is local on $D_e$ in SDG-IBM and global on $\Omega$ in FE-IBM proposed by \cite{boffi03}.
Despite the difference in the variational equations, the resulting formula of the discrete source term in \eqref{eq:source_discrete} is identical to that of \revr{\cite{boffi07}}.

\subsection{Discretization of Euler-Lagrange equation}
\label{sec:EL}
Finally, we discuss the full discretization of Euler-Lagrange equation \eqref{eq:IB_EL2}. For $n = 1, 2, \ldots, K$, given $\bsym{X}^{n-1}_h \in S^h$ from the previous time step and a fixed-point solution $(\bsym{u}^n_h, p^n_h) \in [U^h]^2 \times P^h$ of \eqref{eq:IB_NS2}, we obtain the postprocessed velocity $\bsym{u}^{n,\star}_h$ from $\bsym{u}^n_h$ as discussed in Section \ref{sec:post}. The immersed boundary at time $t = t_n$ is then evolved by
\begin{equation}
\begin{split}
\bsym{X}^{n}_h(s_i) & = \bsym{X}^{n-1}_h(s_i) + \Delta t \, \bsym{u}^{n,\star}_h\left(\bsym{X}^{n-1}_h(s_i)\right) \; \mbox{ for } i = 0, 1, 2, \ldots, m,\\
\bsym{X}^{0}_h(s_i) & = \bsym{X}_0(s_i) \; \mbox{ for } i = 0, 1, 2, \ldots, m.
\end{split}
\label{eq:IB_EL3}
\end{equation}

\subsection{Summary of SDG-IBM}
\label{sec:summary}
The fully discrete SDG-IBM for numerically solving \eqref{eq:IB_NS}--\eqref{eq:IB_EL} is summarized as follows:
for $n = 1,2,\ldots,K$, given $\bsym{u}^{n-1}_h \in [U^h]^2$ and $\bsym{X}^{n-1}_h \in S^h$ from the previous time step,
\begin{enumerate}
\item let $\bsym{u}_{0,h}^n = \bsym{u}_h^{n-1}$ be the initial guess of the sequence of fixed-point iterations,
\item for $m = 1,2,\ldots$, given $\bsym{u}_{m-1,h}^n$ from the previous iteration,
\begin{enumerate}
\item obtain the postprocessed velocity $\bsym{V}_m^n = \bsym{u}_{m-1,h}^{n,\star}$ from $\bsym{u}_{m-1,h}^n$ by \eqref{eq:post1}--\eqref{eq:post4},
\item let $\alpha = \dfrac{\rho}{\Delta t}$, $\bsym{F} = \dfrac{\rho}{\Delta t} \bsym{u}_h^{n-1} + \bsym{F}^n$ and $\bsym{V} = \bsym{V}_m^n$ to obtain the linear system \eqref{eq:IB_NS3},
\item compute the discrete source term of $(\bsym{F}^n, \bsym{v})_{0,\Omega}$ for all $\bsym{v} \in [U^h]^2$ according to \eqref{eq:source_discrete},
\item formulate the system of linear equations \eqref{eq:NS_system} for the SDG method \eqref{eq:NS_discrete},
\item obtain the numerical solution $(\bsym{u}_{m,h}^n, p_{m,h}^n) \in [U^h]^2 \times P^h$,
\end{enumerate}
until a suitably specified stopping criterion is satisfied, and let $(\bsym{u}_h^n, p_h^n) \in [U^h]^2 \times P^h$ be the termination of the sequence of fixed-point iterations,
\item obtain the postprocessed velocity $\bsym{u}_h^{n,\star}$ from $\bsym{u}_h^n$ by \eqref{eq:post1}--\eqref{eq:post4},
\item obtain the new immersed boundary particle configuration $\bsym{X}^n_h \in S^h$ by \eqref{eq:IB_EL3}.
\end{enumerate}
We remark that despite the computation of the source term is placed under the inner iterations in the above procedure, the source term is independent of $m$ and needs to be computed only once for each time level.

\section{Stability analysis}
\label{sec:stab}
In this section, we will provide a stability analysis of SDG-IBM similar to \cite{boffi07}.
First, we introduce some tools which will facilitate our analysis.
The space $Q^h$ of piecewise polynomials on $\tau \in \mathcal{T}$ is defined by
\begin{equation}
Q^{h}
= \{ v \: : \: v|_{\tau} \in P^{k}(\tau); \; \tau\in\mathcal{T}\}.
\label{eq:Vh}
\end{equation}
We define the broken $H^1$ semi-norm $\vert \cdot \vert_{1,*}$ on $[Q^h]^2$ by
\begin{equation}
\vert \bsym{v} \vert_{1,*} = \left( \| \nabla_h \bsym{v} \|_{0, \Omega}^2 + \sum_{e \in \mathcal{F}} h_e^{-1} \| [\bsym{v}] \|_{0, e}^2 \right)^\frac{1}{2} \quad \text{ for } \bsym{v} \in [Q^h]^2.
\label{eq:h1semi}
\end{equation}
Note that the broken $H^1$ semi-norm $\vert \cdot \vert_{1,*}$ coincides with the energy norm $\| \cdot \|_h$ on $[U^h]^2$.

We begin with the following stability result:
\begin{lemma}
\label{stab1}
Let $(\bsym{u}_h, \bsym{w}_h, \bsym{z}_h, \widetilde{\bsym{w}}_h, \widetilde{\bsym{z}}_h, p_h) \in [U^h]^2 \times [W^h]^4 \times P^h$ be the solution of \eqref{eq:NS_discrete}. Then we have
\begin{equation}
\alpha \|\bsym{u}_h\|_{0, \Omega}^2 + \beta^2 \mu  \vert \bsym{u}_h \vert_{1, *}^2 \leq (\bsym{F}, \bsym{u}_h)_{0,\Omega},
\label{eq:stab10}
\end{equation}
where $\beta$ is the inf-sup constant $\beta_1$ in \eqref{eq:inf-sup-B}.
\begin{proof}
In \eqref{eq:NS_discrete}, we take test functions as follows:
\begin{equation}
\begin{split}
\bsym{v} & = \bsym{u}_h, \\
\bsym{\Psi}_1 & = -\bsym{w}_h, \\
\bsym{\Psi}_2 & = -\bsym{z}_h, \\
\bsym{\Psi}_3 & = -\frac{\rho}{2} \widehat{\bsym{w}}_h, \\
\bsym{\Psi}_4 & = -\frac{\rho}{2} \widehat{\bsym{z}}_h, \\
q & = -p_h,
\end{split}
\label{eq:stab11}
\end{equation}
where the definitions of $\widehat{\bsym{w}}_h$ and $\widehat{\bsym{z}}_h$ are given in \eqref{eq:newvar2}.
We then have
\begin{equation}
\begin{split}
\alpha (\bsym{u}_h,\bsym{u}_h)_{0,\Omega}  + \sqrt{\mu} B_h(\bsym{w}_h,u_{h,1})+ \sqrt{\mu} B_h(\bsym{z}_h,u_{h,2})+ & \\
\frac{\rho}{2}R_h\left(\widehat{\bsym{w}}_h,u_{h,1}\right) + \frac{\rho}{2}R_h\left(\widehat{\bsym{z}}_h,u_{h,2}\right)
+  b_h^*(p_h, \bsym{u}_h)& = (\bsym{F}, \bsym{u}_h)_{0,\Omega},\\
-\sqrt{\mu} B_h^*(u_{h,1},\bsym{w}_h) + \sqrt{\mu} (\widehat{\bsym{w}}_h,\bsym{w}_h)_{0,\Omega} & = 0,\\
-\sqrt{\mu} B_h^*(u_{h,2},\bsym{z}_h) + \sqrt{\mu} (\widehat{\bsym{z}}_h,\bsym{z}_h)_{0,\Omega} & = 0,\\
-\frac{\rho}{2} R_h^*(u_{h,1},\widehat{\bsym{w}}_h) + \frac{\rho}{2}(\widetilde{\bsym{w}}_h,\widehat{\bsym{w}}_h)_{0,\Omega} & = 0,\\
-\frac{\rho}{2} R_h^*(u_{h,2},\widehat{\bsym{z}}_h) +\frac{\rho}{2}(\widetilde{\bsym{z}}_h,\widehat{\bsym{z}}_h)_{0,\Omega} & = 0.\\
-b_h(\bsym{u}_h,p_h)&=0.
\end{split}
\label{eq:stab12}
\end{equation}
Summing up all the equations in \eqref{eq:stab12}, using the adjoint relations \eqref{eq:adjoint-B}, \eqref{eq:adjoint-b} and \eqref{eq:adjoint-r} and combining the terms, we have
\begin{equation}
\alpha \|\bsym{u}_h\|_{0, \Omega}^2 + \mu \|L_h\|_{0, \Omega}^2 = (\bsym{F}, \bsym{u}_h)_{0,\Omega}.
\label{eq:stab13}
\end{equation}
Next, by the first inf-sup condition of $U^h$ and $W^h$ in \eqref{eq:inf-sup-B}
\revr{ and then using \eqref{eq:adjoint-B}}, for all $v \in U^h$, we have
\begin{equation}
\| v \|_Z \leq \frac{1}{\beta} \sup_{\bsym{\Psi} \in W^h} \frac{B_h^*(v, \bsym{\Psi})}{\| \bsym{\Psi} \|_{X'}} \leq \frac{1}{\beta} \sup_{\bsym{\Psi} \in W^h} \frac{B_h^*(v, \bsym{\Psi})}{\| \bsym{\Psi} \|_{0,\Omega}},
\label{eq:stab14}
\end{equation}
for all $v \in U^h$. By the second equation of \eqref{eq:NS_discrete}, we have
\begin{equation}
\begin{split}
\| u_{h,1} \|_Z
& \leq \frac{1}{\beta} \sup_{\bsym{\Psi} \in W^h} \frac{B_h^*(u_{h,1}, \bsym{\Psi})}{\| \bsym{\Psi} \|_{0,\Omega}}\\
& = \frac{1}{\beta} \sup_{\bsym{\Psi} \in W^h} \frac{(\widehat{\bsym{w}}_h, \bsym{\Psi})_{0,\Omega}}{\| \bsym{\Psi} \|_{0,\Omega}},\\
& = \frac{1}{\beta} \|\widehat{\bsym{w}}_h \|_{0,\Omega}.
\end{split}
\label{eq:stab15a}
\end{equation}
Similarly, we have
\begin{equation}
\| u_{h,2} \|_Z
\leq \frac{1}{\beta} \|\widehat{\bsym{z}}_h \|_{0,\Omega}.
\label{eq:stab15b}
\end{equation}
Combining \eqref{eq:stab15a} and \eqref{eq:stab15b}, we obtain
\begin{equation}
\vert \bsym{u}_h \vert_{1, *} \leq \frac{1}{\beta} \| L_h \|_{0,\Omega}.
\label{eq:stab15}
\end{equation}
Substituting \eqref{eq:stab15} into \eqref{eq:stab13}, we have
\begin{equation}
\alpha \| \bsym{u}_h \|_{0, \Omega}^2 + \beta^2 \mu \vert \bsym{u}_h \vert_{1, *}^2 \leq (\bsym{F}, \bsym{u}_h)_{0,\Omega}.
\label{eq:stab16}
\end{equation}
\end{proof}
\end{lemma}
One important thing to note is that due to the skew-symmetric discretization of convection term, the convection velocity $\bsym{V}$ vanishes in the above estimate and the stability is therefore enhanced.

Now we are ready to present the following stability estimate:
\begin{theorem}
\label{stab2}
Let $(\bsym{u}_h^{n}, p_h^n, \bsym{X}_h^n) \in [U^h]^2 \times P^h \times S^h$ be the approximated solution of \eqref{eq:IB_NS}--\eqref{eq:IB_EL} at $t = t_n$ obtained by SDG-IBM discussed in Section \ref{sec:summary}. Then for $n = 1, 2, \ldots, K$, we have
\begin{equation}
\begin{split}
& \frac{\rho}{2\Delta t} \left( \|\bsym{u}_h^n\|_{0, \Omega}^2 - \|\bsym{u}_h^{n-1}\|_{0, \Omega}^2 \right) + \beta^2 \mu \vert \bsym{u}_h^n \vert_{1, *}^2 + \frac{\kappa}{2\Delta t} \left( \left\| \frac{\partial \bsym{X}_h^{n}}{\partial s} \right\|_{0,D}^2  - \left\| \frac{\partial \bsym{X}_h^{n-1}}{\partial s} \right\|_{0,D}^2  \right)\\
& \qquad \leq \left(\frac{C \kappa}{h_s h_x^\frac{1}{2}}  (L^{n-1})^\frac{3}{2} + \frac{C^\prime \kappa}{h_s}  L^{n-1}\right) \vert \bsym{u}_h^n - \bsym{u}_h^{n,\star} \vert_{1,*} +  \left(\frac{C^2 \kappa \Delta t}{h_s h_x} L^{n-1} + \frac{(C^\prime)^2 \kappa \Delta t}{h_s}\right) \vert \bsym{u}_h^{n,\star} \vert_{1,*}^2,
\end{split}
\label{eq:stab20}
\end{equation}
where $h_s$, $h_x$ and $L^{n-1}$ are defined as
\begin{equation}
\begin{split}
h_s & = \min_{1 \leq i \leq m} \vert s_i - s_{i-1} \vert, \\
h_x & = \min_{1 \leq i \leq m} \text{diam} (\hat{T}_i),\\
L^{n-1} & = \max_{1 \leq i \leq m} \vert \bsym{X}_h^{n-1} (s_i) -  \bsym{X}_h^{n-1} (s_{i-1}) \vert,
\end{split}
\end{equation}
and $\hat{T}_i$ is the union of all elements in $\mathcal{T}$ intersecting the segment joining $\bsym{X}_h^{n-1} (s_i)$ to $\bsym{X}_h^{n-1} (s_{i-1})$. All the constants appeared
in the above estimates are independent of discretization parameters, $h_x$, $h_s$, and $\Delta t$.
\begin{proof}
We recall Section~\ref{sec:summary}, $\alpha = \dfrac{\rho}{\Delta t}$, $\bsym{F} = \dfrac{\rho}{\Delta t} \bsym{u}_h^{n-1} + \bsym{F}^n$, and the discrete form $(\bsym{F}_h^n, \bsym{v})_{0,\Omega}$ of $(\bsym{F}^n,\bsym{v})_{0,\Omega}$ in \eqref{eq:source_discrete}.
We then have
\begin{equation}
\frac{\rho}{\Delta t} \|\bsym{u}_h^n\|_{0, \Omega}^2 + \beta^2 \mu \vert \bsym{u}_h^n \vert_{1, *}^2 \leq \frac{\rho}{\Delta t} (\bsym{u}_h^{n-1}, \bsym{u}_h^n)_{0, \Omega} + (\bsym{F}_h^n, \bsym{u}_h^n)_{0,\Omega},
\label{eq:stab21}
\end{equation}
since the above estimate holds true throughout the sequence of fixed point iterations.
By a direct calculation, we have
\begin{equation}
\| \bsym{u}_h^n - \bsym{u}_h^{n-1} \|_{0,\Omega}^2 = \| \bsym{u}_h^n \|_{0,\Omega}^2 - 2(\bsym{u}_h^{n-1}, \bsym{u}_h^n)_{0, \Omega} + \| \bsym{u}_h^{n-1} \|_{0,\Omega}^2.
\label{eq:stab22}
\end{equation}
Combining \eqref{eq:stab21} and \eqref{eq:stab22}, we obtain
\begin{equation}
\frac{\rho}{2\Delta t} \left( \|\bsym{u}_h^n\|_{0, \Omega}^2 - \|\bsym{u}_h^{n-1}\|_{0, \Omega}^2 + \| \bsym{u}_h^n - \bsym{u}_h^{n-1} \|_{0,\Omega}^2 \right) + \beta^2 \mu \vert \bsym{u}_h^n \vert_{1, *}^2 \leq (\bsym{F}_h^n, \bsym{u}_h^n)_{0,\Omega}.
\label{eq:stab23}
\end{equation}

We let $\bsym{v}_h^n = \bsym{u}_h^n - \bsym{u}_h^{n, \star}$. By \eqref{eq:source_discrete} and a rearrangement of indices, we have
\begin{equation}
\begin{split}
(\bsym{F}_h^n ,\bsym{u}_h^n)_{0,\Omega}
& = \sum_{i=1}^m \kappa  \left(\frac{\partial \bsym{X}_h^{n-1}}{\partial s}(s_{i+\frac{1}{2}}) - \frac{\partial \bsym{X}_h^{n-1}}{\partial s}(s_{i-\frac{1}{2}}) \right)\bsym{u}_h^n(\bsym{X}^{n-1}_h(s_i))\\
& = \sum_{i=1}^m \kappa  \frac{\partial \bsym{X}_h^{n-1}}{\partial s}(s_{i-\frac{1}{2}}) \Big(\bsym{u}_h^{n}(\bsym{X}^{n-1}_h(s_{i-1})) - \bsym{u}_h^{n}(\bsym{X}^{n-1}_h(s_i)) \Big) \\
& = \sum_{i=1}^m \kappa  \frac{\partial \bsym{X}_h^{n-1}}{\partial s}(s_{i-\frac{1}{2}})\Big( \bsym{v}_h^n(\bsym{X}^{n-1}_h(s_{i-1})) - \bsym{v}_h^n(\bsym{X}^{n-1}_h(s_i)) \Big) \\
& \qquad + \sum_{i=1}^m \kappa \frac{\partial \bsym{X}_h^{n-1}}{\partial s}(s_{i-\frac{1}{2}})\Big( \bsym{u}_h^{n,\star}(\bsym{X}^{n-1}_h(s_{i-1})) - \bsym{u}_h^{n,\star}(\bsym{X}^{n-1}_h(s_i)) \Big).
\end{split}
\label{eq:stab24}
\end{equation}

We first consider the first sum on the last equality of \eqref{eq:stab24}.
For $i = 1,2,\ldots,m$, we let $\gamma_i$ be the line segment connecting $\bsym{X}^{n-1}_h(s_{i-1})$ and $\bsym{X}^{n-1}_h(s_{i})$. Then we have
\begin{equation}
\Big\vert \bsym{v}_h^n(\bsym{X}^{n-1}_h(s_{i-1})) - \bsym{v}_h^n(\bsym{X}^{n-1}_h(s_i)) \Big\vert \leq \sum_l \int_{\gamma_{i,l}} \left\vert  \nabla \bsym{v}_h^n \right\vert + \sum_{e \in \mathcal{F}: e \cap \gamma_i \neq \phi} \left\vert [\bsym{v}_h^n] \vert_{e \cap \gamma_i} \right\vert,
\label{eq:stab25a0}
\end{equation}
where $\gamma_{i,l}$ are subsegments of $\gamma_i$ cut by the edges $e$ intersecting $\gamma_i$.
For the first term on the right hand side of \eqref{eq:stab25a0}, by an inverse inequality and then a trace inequality, we have
\begin{equation}
\begin{split}
\sum_l \int_{\gamma_{i,l}} \left\vert  \nabla \bsym{v}_h^n \right\vert
& \leq \Big\vert \bsym{X}^{n-1}_h(s_{i-1}) - \bsym{X}^{n-1}_h(s_i) \Big\vert^\frac{1}{2}
\sum_l \| \nabla \bsym{v}_h^n \|_{0,\gamma_{i,l}}  \\
& \leq C_0 (h_x)^{-\frac{1}{2}} (L^{n-1})^\frac{1}{2} \| \nabla \bsym{v}_h^n \|_{0,\hat{T}_i},
\end{split}
\label{eq:stab25a1}
\end{equation}
where $C_0$ is a constant depending on the number of $\gamma_{i,l}$'s but independent of $h$.
For the second term on the right hand side of \eqref{eq:stab25a0}, by the norm equivalence on the space of polynomials defined on the edges, we note that
\begin{equation}
\begin{split}
\sum_{e \in \mathcal{F}: e \cap \gamma_i \neq \phi} \left\vert [\bsym{v}_h^n] \vert_{e \cap \gamma_i} \right\vert
& \leq  \sum_{e \in \mathcal{F}: e \cap \gamma_i \neq \phi} \| [\bsym{v}_h^n] \|_{\infty, e} \\
& \leq C_0^\prime \sum_{e \in \mathcal{F}: e \cap \gamma_i \neq \phi} h_e^{-\frac{1}{2}} \| [\bsym{v}_h^n] \|_{0, e},
\end{split}
\label{eq:stab25a2}
\end{equation}
where $C_0^\prime$ is a constant independent of $h$.

Thus, we have
\begin{equation}
\begin{split}
& \; \left \vert \sum_{i=1}^m \kappa  \frac{\partial \bsym{X}_h^{n-1}}{\partial s}(s_{i-\frac{1}{2}})\Big( \bsym{v}_h^n(\bsym{X}^{n-1}_h(s_{i-1})) - \bsym{v}_h^n(\bsym{X}^{n-1}_h(s_i)) \Big)  \right \vert \\
\leq & \; \sum_{i=1}^m \kappa \left\vert \frac{\partial \bsym{X}_h^{n-1}}{\partial s}(s_{i-\frac{1}{2}}) \right\vert \Big\vert \bsym{v}_h^n(\bsym{X}^{n-1}_h(s_{i-1})) - \bsym{v}_h^n(\bsym{X}^{n-1}_h(s_i)) \Big\vert\\
\leq & \; \frac{\kappa}{h_s} L^{n-1} \sum_{i=1}^m \Big\vert \bsym{v}_h^n(\bsym{X}^{n-1}_h(s_{i-1})) - \bsym{v}_h^n(\bsym{X}^{n-1}_h(s_i)) \Big\vert\\
\leq & \; \frac{\kappa}{h_s} L^{n-1} \left(C_0 (h_x)^{-\frac{1}{2}} (L^{n-1})^\frac{1}{2} \sum_{i=1}^m \| \nabla \bsym{v}_h^n \|_{0,\hat{T}_i} + C_0^\prime \sum_{i=1}^m \sum_{e \in \mathcal{F}: e \cap \gamma_i \neq \phi} h_e^{-\frac{1}{2}} \| [\bsym{v}_h^n] \|_{0, e}\right) \\
\leq & \; \frac{\kappa}{h_s} L^{n-1} \left(C_0 C_1 (h_x)^{-\frac{1}{2}} (L^{n-1})^\frac{1}{2} \| \nabla \bsym{v}_h^n \|_{0, \Omega} + C_0^\prime C_1^\prime \sum_{e \in \mathcal{F}} h_e^{-\frac{1}{2}} \| [\bsym{v}_h^n] \|_{0, e}\right) \\
\leq & \; \left(\frac{C_0 C_1\kappa}{h_s h_x^\frac{1}{2}} (L^{n-1})^\frac{3}{2} + \frac{C_0^\prime C_1^\prime \kappa}{h_s} L^{n-1}\right)  \vert \bsym{u}_h^n - \bsym{u}_h^{n,\star} \vert_{1,*},
\end{split}
\label{eq:stab25a3}
\end{equation}
where $C_1$ is the maximum number of segments $\gamma_i$ intersected with an element $\tau \in \mathcal{T}$, and  $C_1^\prime$ is the maximum number of segments $\gamma_i$ intersected with an edge $e \in \mathcal{F}$.

For the second sum on the last equality of \eqref{eq:stab24}, by \eqref{eq:IB_EL3}, by a rearrangement of indices and a trick similar to \eqref{eq:stab22}, we have
\begin{equation}
\begin{split}
& \; \sum_{i=1}^m \kappa \frac{\partial \bsym{X}_h^{n-1}}{\partial s}(s_{i-\frac{1}{2}})\Big( \bsym{u}_h^{n,\star}(\bsym{X}^{n-1}_h(s_{i-1})) - \bsym{u}_h^{n,\star}(\bsym{X}^{n-1}_h(s_i)) \Big)\\
= & \; \sum_{i=1}^m \kappa
\left(\frac{\bsym{X}_h^{n-1}(s_{i}) - \bsym{X}_h^{n-1}(s_{i-1})}{s_{i} - s_{i-1}}\right)
\left(\frac{\bsym{X}_h^n(s_{i-1}) - \bsym{X}_h^{n-1}(s_{i-1})}{\Delta t} - \frac{\bsym{X}_h^n(s_i) - \bsym{X}_h^{n-1}(s_i)}{\Delta t}\right) \\
= & \; \sum_{i=1}^m \kappa
\left(\frac{\bsym{X}_h^{n-1}(s_{i}) - \bsym{X}_h^{n-1}(s_{i-1})}{s_{i} - s_{i-1}}\right)
\left(\frac{\bsym{X}_h^{n-1}(s_i)- \bsym{X}_h^{n-1}(s_{i-1})}{\Delta t} - \frac{\bsym{X}_h^n(s_i) - \bsym{X}_h^n(s_{i-1})}{\Delta t}\right) \\
= & \; \sum_{i=1}^m \frac{\kappa(s_i - s_{i-1})}{\Delta t}
\frac{\partial \bsym{X}_h^{n-1}}{\partial s} (s_{i-\frac{1}{2}})
\left(\frac{\partial \bsym{X}_h^{n-1}}{\partial s} (s_{i-\frac{1}{2}}) - \frac{\partial \bsym{X}_h^{n}}{\partial s} (s_{i-\frac{1}{2}})\right) \\
= & \;
 \; \frac{\kappa}{2\Delta t} \left(
\left\| \frac{\partial \bsym{X}_h^{n-1}}{\partial s} \right\|_{0,D}^2
- \left\| \frac{\partial \bsym{X}_h^{n}}{\partial s} \right\|_{0,D}^2
+ \left\| \frac{\partial \bsym{X}_h^{n-1}}{\partial s} - \frac{\partial \bsym{X}_h^{n}}{\partial s}  \right\|_{0,D}^2
\right).
\end{split}
\label{eq:stab25b0}
\end{equation}
Again, using a similar argument as \eqref{eq:stab25a1} and \eqref{eq:stab25a2} on $\bsym{u}_h^{n,\star}$, we obtain the following estimate:
\begin{equation}
\begin{split}
& \;
\left\| \frac{\partial \bsym{X}_h^{n-1}}{\partial s} - \frac{\partial \bsym{X}_h^{n}}{\partial s}  \right\|_{0,D}^2\\
= & \;
\sum_{i=1}^m  \frac{1}{s_i - s_{i-1}} \Big((\bsym{X}_h^{n-1}(s_i) - \bsym{X}_h^{n-1}(s_{i-1})) - (\bsym{X}_h^{n}(s_i) - \bsym{X}_h^{n}(s_{i-1}))\Big)^2\\
= & \;
\sum_{i=1}^m  \frac{1}{s_i - s_{i-1}} \Big((\bsym{X}_h^{n}(s_{i-1}) - \bsym{X}_h^{n-1}(s_{i-1})) - (\bsym{X}_h^{n}(s_{i}) - \bsym{X}_h^{n-1}(s_{i}))\Big)^2\\
= & \;
\sum_{i=1}^m \frac{\Delta t^2}{s_i - s_{i-1}} \Big( \bsym{u}_h^{n,\star}(\bsym{X}_h^{n-1}(s_{i-1})) - \bsym{u}_h^{n,\star}(\bsym{X}_h^{n-1}(s_{i})) \Big)^2\\
\leq & \;
\frac{2 \Delta t^2}{h_s}\left(\frac{(C_0 C_1)^2}{h_x} L^{n-1} + (C_0^\prime C_1^\prime)^2 \right)  \vert \bsym{u}_h^{n,\star} \vert_{1,*}^2.
\end{split}
\label{eq:stab25b1}
\end{equation}
Combining \eqref{eq:stab25b0} and \eqref{eq:stab25b1}, we have
\begin{equation}
\begin{split}
& \; \sum_{i=1}^m \kappa \frac{\partial \bsym{X}_h^{n-1}}{\partial s}(s_{i-\frac{1}{2}})\Big( \bsym{u}_h^{n,\star}(\bsym{X}^{n-1}_h(s_{i-1})) - \bsym{u}_h^{n,\star}(\bsym{X}^{n-1}_h(s_i)) \Big)\\
\leq & \;
 \; \frac{\kappa}{2\Delta t} \left(
\left\| \frac{\partial \bsym{X}_h^{n-1}}{\partial s} \right\|_{0,D}^2
- \left\| \frac{\partial \bsym{X}_h^{n}}{\partial s} \right\|_{0,D}^2
+ \frac{2 \Delta t^2}{h_s}\left(\frac{(C_0 C_1)^2}{h_x} L^{n-1} + (C_0^\prime C_1^\prime)^2 \right)  \vert  \bsym{u}_h^{n,\star} \vert_{1,*}^2 \right).
\end{split}
\label{eq:stab25b2}
\end{equation}
Finally, combining \eqref{eq:stab23}, \eqref{eq:stab24}, \eqref{eq:stab25a3} and \eqref{eq:stab25b2}, we obtain the desired result.
\end{proof}
\end{theorem}

Using Theorem \ref{stab2}, we can now establish the following CFL condition for SDG-IBM:
\begin{corollary}
\label{stab3}
For $n = 1, 2, \ldots, K$, define a CFL parameter $\eta^{n-1}$ by
\begin{equation}
\eta^{n-1} = \frac{\kappa \Delta t}{h_s} \left(1 + \frac{L^{n-1}}{h_x}\right).
\label{eq:CFLparam}
\end{equation}
Assume that $h_s = O(h)$ and $h_x = O(h)$,
and there exists a uniform constant $K_0$ such that
\begin{equation}
\beta^2 \mu - 2\max\{C, C^\prime\}^2 \eta^{j-1} \geq K_0 > 0 \text{ for } j = 1, 2, \ldots n,
\label{eq:CFLcond1}
\end{equation}
then the following energy property holds:
\begin{equation}
E^n \leq E^0 + R^n,
\end{equation}
where $E^n$ and $R^n$ are defined as
\begin{equation}
\begin{split}
E^n & = \frac{\rho}{2} \| \bsym{u}_h^n \|_{0, \Omega}^2 + \Delta t \sum_{j = 1}^n K_0 \vert \bsym{u}_h^j \vert_{1,*}^2 + \frac{\kappa}{2} \left\| \frac{\partial \bsym{X}_h^{n}}{\partial s} \right\|_{0,D}^2,\\
R^n & = \Delta t \sum_{j = 1}^n \left( \left(\frac{\widetilde{C} \kappa}{h^\frac{1}{2}}  (L^{j-1})^\frac{3}{2} + \widetilde{C}^\prime \kappa L^{j-1} \right) \| \bsym{u}^j \|_{2,\Omega} + 2\left(\widetilde{C}^2 \kappa \Delta t L^{j-1} + (\widetilde{C}^\prime)^2 \kappa h \Delta t \right) \| \bsym{u}^j \|_{2,\Omega}^2\right),
\end{split}
\end{equation}
and $\bsym{u}^j$ is the analytic solution of \eqref{eq:IB_NS2} at $t = t_j$.
All the constants in the above are independent of discretization parameters, $h_x$, $h_s$ and $\Delta t$.
\begin{proof}
By the triangle inequality and Cauchy-Schwarz inequality, we have
\begin{equation}
\vert \bsym{u}_h^{j,\star} \vert_{1,*}^2 \leq 2 \left( \vert \bsym{u}_h^{j} \vert_{1,*}^2 + \vert \bsym{u}_h^{j} - \bsym{u}_h^{j,\star} \vert_{1,*}^2 \right).
\label{eq:stab31}
\end{equation}
Applying this inequality to the right hand side of \eqref{eq:stab20}, rearranging the terms and using the assumption \eqref{eq:CFLcond1}, we obtain
\begin{equation}
\begin{split}
& \frac{\rho}{2\Delta t} \left( \|\bsym{u}_h^j\|_{0, \Omega}^2 - \|\bsym{u}_h^{j-1}\|_{0, \Omega}^2 \right) + K_0 \vert \bsym{u}_h^j \vert_{1, *}^2 + \frac{\kappa}{2\Delta t} \left( \left\| \frac{\partial \bsym{X}_h^{j}}{\partial s} \right\|_{0,D}^2  - \left\| \frac{\partial \bsym{X}_h^{j-1}}{\partial s} \right\|_{0,D}^2  \right)\\
& \qquad \leq \left(\frac{C \kappa}{h_s h_x^\frac{1}{2}}  (L^{n-1})^\frac{3}{2} + \frac{C^\prime \kappa}{h_s}  L^{n-1}\right) \vert \bsym{u}_h^j - \bsym{u}_h^{j,\star} \vert_{1,*} + 2 \left(\frac{C^2 \kappa \Delta t}{h_s h_x} L^{n-1} + \frac{(C^\prime)^2 \kappa \Delta t}{h_s}\right) \vert \bsym{u}_h^j - \bsym{u}_h^{j,\star} \vert_{1,*}^2.
\end{split}
\label{eq:stab32}
\end{equation}
Additionally, from \cite{sdg-hdg1} and \cite{sdg-ns2}, we have the following estimate
\begin{equation}
\vert \bsym{u}_h^j - \bsym{u}_h^{j,\star} \vert_{1,*} \leq C_u h \|\bsym{u}^j\|_{2,\Omega}.
\label{eq:stab33}
\end{equation}
By combining \eqref{eq:stab32} and \eqref{eq:stab33} and assuming $h_s = O(h)$ and $h_x = O(h)$, we have
\begin{equation}
\begin{split}
& \; \frac{1}{\Delta t} (E^j - E^{j-1}) \\
\leq & \left(\frac{C \kappa}{h^\frac{3}{2}}  (L^{n-1})^\frac{3}{2} + \frac{C^\prime \kappa}{h}  L^{n-1}\right) \left(C_u h \|\bsym{u}^j\|_{2,\Omega}\right) + 2\left(\frac{C^2 \kappa \Delta t}{h^2} L^{n-1} + \frac{(C^\prime)^2 \kappa \Delta t}{h}\right) \left(C_u h \|\bsym{u}^j\|_{2,\Omega}\right)^2 \\
\leq & \; \left(\frac{\widetilde{C} \kappa}{h^\frac{1}{2}}  (L^{n-1})^\frac{3}{2} + \widetilde{C}^\prime \kappa  L^{n-1}\right) \|\bsym{u}^j\|_{2,\Omega} + 2\left(\widetilde{C}^2 \kappa \Delta t L^{n-1} + (\widetilde{C}^\prime)^2 \kappa h \Delta t \right) \|\bsym{u}^j\|_{2,\Omega}^2.
\end{split}
\label{eq:stab34}
\end{equation}
Summing over $j = 1, 2, \ldots, n$, we obtain the desired result.
\end{proof}
\end{corollary}

We would like to make a few remarks here. From Corollary \ref{stab3}, the condition on $\Delta t$ for the stability is similar to that obtained in \cite{boffi07}. In \cite{boffi07}, linearized Navier-Stokes problem is considered, and we obtained the above stability result for the nonlinear problem by using the skew-symmetry property of the nonlinear term. This is an advantage of staggered DG formulation and obtained by using the splitting of the diffusion and convection term. The use of the post-processed velocity $\bsym{u}_h^{n,\star}$ in advancing the structure gives the additional term $R^n$ in the stability estimate. By using the post-processed velocity, we can observe a better mass-preserving property of the discrete problem.

We note that by assuming that $\bsym{X}^n(s)$ are uniformly Lipschitz for all time step $n$
we can bound $L^{n-1}$ by $C h$ for all $n$. The first term in $R^n$ can then be bounded
by $C h \Delta t \sum_{j=1}^n \|u^j\|_{2,\Omega}$ and the second term by
$C  h \Delta t \sum_{j=1}^n \|u^j\|_{2,\Omega}^2$.
We can then bound the term $R^n$ by $C h \sum_{j=1}^n \Delta t \| u^j \|_{2,\Omega}^2$,
which says that the term $R^n$ becomes harmless for a sufficiently small $h$.
\revb{The condition in \eqref{eq:CFLcond1} means that our scheme is more stable
for a model with larger $\mu$ and smaller $\kappa$.}

\section{Numerical results}
\label{sec:num}
In this section, we illustrate some numerical examples. We carry out numerical experiments to see the area conservation of the immersed boundary and the stability of the proposed method. In Sections \ref{sec:secexp1} and \ref{sec:secexp2}, we present numerical results of an ellipse and an L-shaped curve immersed in a static fluid. In Section \ref{sec:secexp3}, we perform an experiment to see an ellipse immersed in a rotating fluid. In Section \ref{sec:secexp4}, we examine the behaviour of a stretched curve immersed in a static fluid.
In Section~\ref{sec:secexp5}, we present stability of our method for a test example.

Polynomials with degree $k=1$ is used in the SDG spatial discretization. Throughout the experiments in the whole Section \ref{sec:num}, unless otherwise specified, the Lagrangian mesh defined in \eqref{eq:Lpart} is uniform. The physical quantities are set to be:
\begin{equation}
\rho = 1, \mu = 1, \kappa = 1 \text{ and } \Delta t = 0.01.
\end{equation}
We denote the number of divisions in $[0,1]$ in the Eulerian mesh by $N$, the number of divisions in $[0,L]$ in the Lagrangian mesh by $m$, and the number of divisions in $[0,T]$ in the temporal mesh by $K$, respectively.

\subsection{Ellipse immersed in a static fluid}
\label{sec:secexp1}
This experiment is to compare the area conversation of SDG-IBM with FE-IBM proposed in \cite{boffi04}.
The initial condition for the fluid motion is given by
\begin{equation}
\bsym{u}_0(x,y) = \bsym{0}, \quad (x,y) \in [0,1]^2.
\end{equation}
The initial configuration of the Lagrangian markers is given by
\begin{equation}
\bsym{X}_0(s) =
\begin{pmatrix} 0.2\cos(2\pi s) + 0.3\\ 0.1 \sin(2\pi s) + 0.3\end{pmatrix},
\quad s \in [0,1].
\end{equation}

Tests are performed with mesh sizes $N = 4, 8, 16, 32$ and $m = 64, 128, 256$ and $K = 200$. At $t = 2$, the area change is analyzed.

\begin{figure}[ht!]
\centering
\includegraphics[width=2.5in]{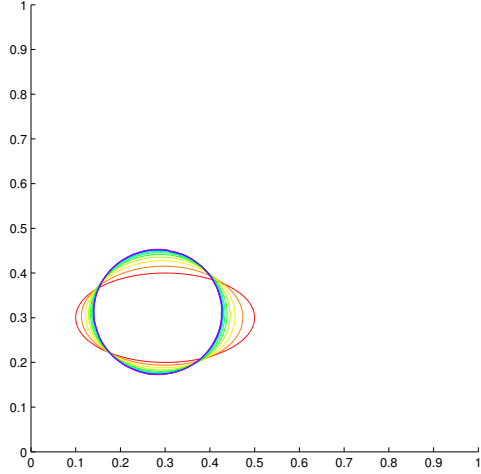}
\caption{Evolution of the immersed boundary throughout $t = 0$ and $t = 2$ in experiment \ref{sec:secexp1}.}
\label{fig:exp1}
\end{figure}

Table \ref{tab:area1} records the area change of the immersed boundary in experiment \ref{sec:secexp1}, and shows that the area conservation of SDG-IBM is very outstanding. With $N = 32$ and $m = 256$, the area loss is 0.07\%, significantly less than 2.3\% of FE-IBM in \cite{boffi04}. 
\revr{We note that in our numerical experiments we calculated the area change of the immersed boundary by comparing the area enclosed by the $m$-sided polygons with vertices $\mathbf{X}_h^K(s_i)$ and $\mathbf{X}_h^0(s_i)$ respectively.} 
Figure \ref{fig:exp1} shows the evolution of the immersed boundary throughout $t = 0$ and $t = 2$ with $N = 32$ and $m = 256$. It can be seen that the Lagrangian markers tend to the equilibrium configuration, which is a circle in shape.
\begin{table}[ht]
\begin{center}
\begin{tabular}{c c|c c c}
& & \multicolumn{3}{c}{$m$} \\
& & 64 & 128 & 256 \\ \hline
\multirow{4}{*}{$N$}
& 4 & 0.2314 & -0.4009 & -4.1318 \\
& 8 & -0.1507 & -0.1466 & -0.1105 \\
& 16 & -0.4359 & -0.1286 & 0.0511 \\
& 32 & -1.3745 & -0.2429 & -0.0763 \\
\end{tabular}
\caption{Percentage of area change in experiment \ref{sec:secexp1}.}
\label{tab:area1}
\end{center}
\end{table}

\subsection{L-shaped curve immersed in a static fluid} 
\label{sec:secexp2}

We consider an experiment which has identical set-ups as experiment \ref{sec:secexp1}, except the initial configuration of the Lagrangian markers is replaced by an L-shaped closed curve.

\begin{figure}[ht!]
\centering
\includegraphics[width=2.5in]{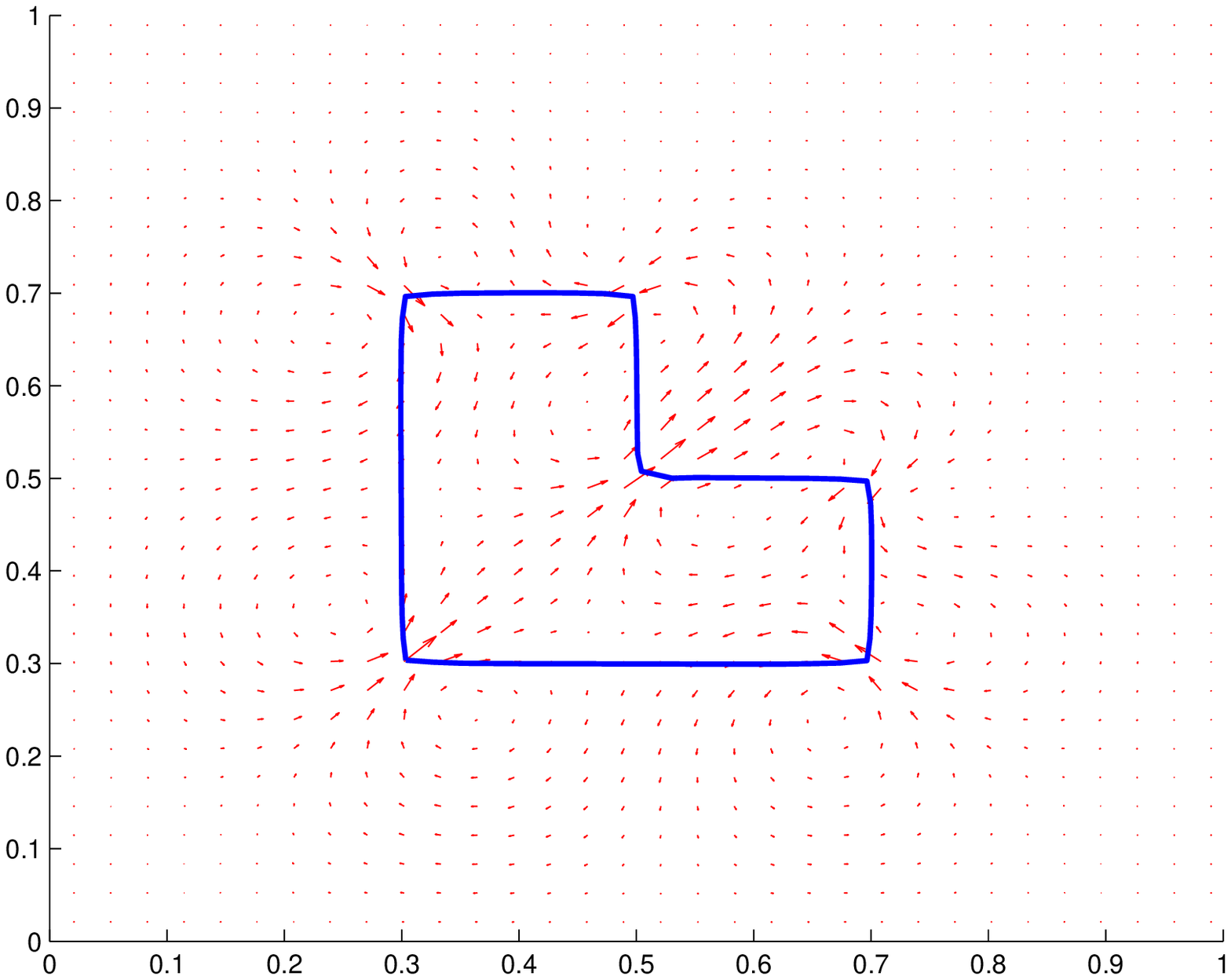}
\caption{Profile of the fluid flow and configuration of the Lagrangian markers at $t = \Delta t$ in experiment \ref{sec:secexp2}.}
\label{fig:exp2a}
\end{figure}

\begin{figure}[ht!]
\centering
\includegraphics[width=2.5in]{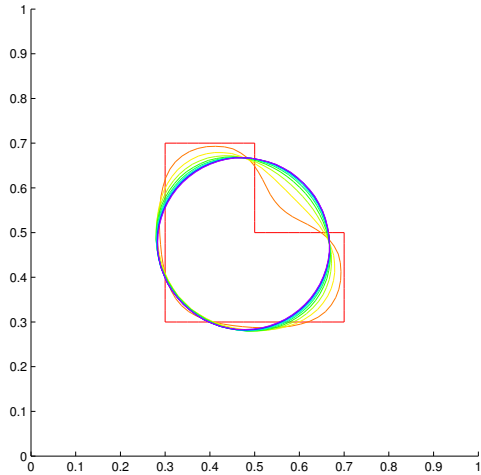}
\caption{Evolution of the immersed boundary throughout $t = 0$ and $t = 2$ in experiment \ref{sec:secexp2}.}
\label{fig:exp2b}
\end{figure}

Table \ref{tab:area2} records the area change of the immersed boundary in experiment \ref{sec:secexp2}. Figure \ref{fig:exp2a} shows the profile of the fluid flow and the configuration of the Lagrangian markers at $t = \Delta t$ with $N = 32$ and $m = 256$. It can be observed that the fluid flow out of the immersed boundary at the inner corner and flow into the immersed boundary at the other corners. The flow substantially pushes the inner corner out. Figure \ref{fig:exp2b} shows the evolution of the immersed boundary throughout $t = 0$ and $t = 2$ with $N = 32$ and $m = 256$. Again, the Lagrangian markers tend to the circular equilibrium configuration.

\begin{table}[ht]
\begin{center}
\begin{tabular}{c c|c c c}
& & \multicolumn{3}{c}{$m$} \\
& & 64 & 128 & 256 \\ \hline
\multirow{4}{*}{$N$}
& 4 & -6.5987 & 4.0340 & -55.1772 \\
& 8 & -0.2046 & 0.1075 & 0.8988 \\
& 16 & -0.6720 & 0.1189 & 0.2321 \\
& 32 & -3.5787 & -0.2989 & -0.0429 \\
\end{tabular}
\caption{Percentage of area change in experiment \ref{sec:secexp2}.}
\label{tab:area2}
\end{center}
\end{table}

\subsection{Ellipse immersed in a rotating fluid} 
\label{sec:secexp3}

\begin{figure}[ht!]
\centering
\includegraphics[width=3.0in]{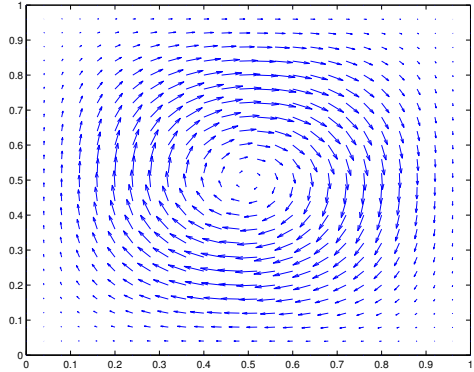}
\caption{Vector plot for the velocity field $\bsym{v}$ in experiment \ref{sec:secexp3}.}
\label{fig:exp3a}
\end{figure}

We consider a model with the immersed boundary driven by a rotating fluid. In addition to the elastic force acting on the fluid by the immersed boundary, an external force for maintaining the Navier-Stokes flow of the rotating velocity field
\begin{equation}
\bsym{v}(x,y) =
\begin{pmatrix} -0.4(1-\cos(2 \pi x)) \sin(2 \pi y) \\  0.4 \sin(2 \pi x) (1 -  \cos(2 \pi y)) \end{pmatrix}, \quad (x,y) \in [0,1]^2
\end{equation}
is added to the source term.
Figure \ref{fig:exp3a} shows a vector plot for the velocity field $\bsym{v}$ on $[0,1]^2$.

The initial condition for the fluid motion is given by
\begin{equation}
\bsym{u}_0(x,y) = \bsym{v}(x,y), \quad (x,y) \in [0,1]^2.
\end{equation}
The initial configuration of the Lagrangian markers is given by
\begin{equation}
\bsym{X}_0(s) =
\begin{pmatrix} 0.2\cos(2\pi s) + 0.4\\ 0.1 \sin(2\pi s) + 0.5\end{pmatrix},
\quad s \in [0,1].
\end{equation}

Tests are performed with mesh sizes $N = 4, 8, 16, 32$ and $m = 64, 128, 256$ and $K = 200$. At $t = 2$, the area change is analyzed.

Table \ref{tab:area3} records the area change of the immersed boundary in experiment \ref{sec:secexp3}. \revb{The result in Table~\ref{tab:area3} is less satisfactory when decreasing $N$ and $m$. 
We note that we have used a uniform time step size $\Delta t=1/100$.
Since the accuracy of time discretization and the stability of the scheme also affects the
area conservation, we test the same model problem with decreasing $\Delta t$ and report the area change
in Table~\ref{tab:area3dt}. For a fixed $m$ and $N$, we can observe the area change decreases
when decreasing $\Delta t$, which confirms our assertion. 
In Tables~\ref{tab:area3dt2} and \ref{tab:area3dt3}, the ratios $h_x/h_s$ and $\Delta t/h_s$ are fixed respectively, 
and the reduction of area change is similar to that in Table~\ref{tab:area3dt}. 
This shows that the time discretization is accounted for 
the relatively poor area conservation in this experiment.
}

\begin{figure}[ht!]
\centering
\includegraphics[width=2.5in]{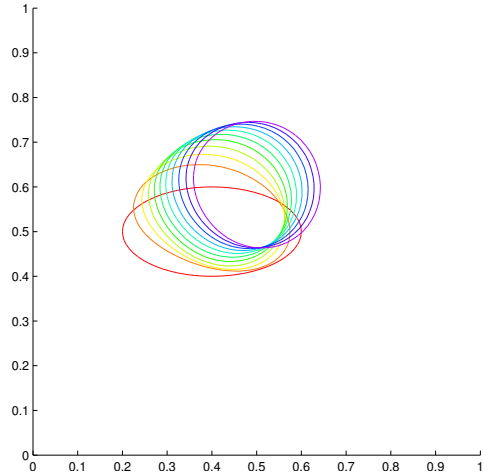}
\caption{Evolution of the immersed boundary throughout $t = 0$ and $t = 2$ in experiment \ref{sec:secexp3}.}
\label{fig:exp3b}
\end{figure}

Figure \ref{fig:exp3b} shows the evolution of the immersed boundary throughout $t = 0$ and $t = 2$ with $N = 32$ and $m = 256$. It can be seen that the Lagrangian markers are driven by the rotating velocity field, and they tend to the circular equilibrium configuration simultaneously.

\begin{table}[ht]
\begin{center}
\begin{tabular}{c c|c c c}
& & \multicolumn{3}{c}{$m$} \\
& & 64 & 128 & 256 \\ \hline
\multirow{4}{*}{$N$}
& 4 & 2.4434 & 1.5623 & 1.6037 \\
& 8 & 1.8129 & 1.7919 & 1.8131 \\
& 16 & 1.8121 & 1.8432 & 1.8081 \\
& 32 & 0.2928 & 1.8271 & 1.7833 \\
\end{tabular}
\caption{Percentage of area change in experiment \ref{sec:secexp3}.}
\label{tab:area3}
\end{center}
\end{table}

\begin{table}[ht]
\begin{center}
\begin{tabular}{c c|c  c  c}
& & \multicolumn{3}{c}{$m$} \\
& & 64 & 128 & 256 \\ \hline
\multirow{3}{*}{$\Delta t$}
& 1/100 & 1.8121 & 1.8432 & 1.8081 \\
& 1/200 & 0.9013 & 0.9724 & 0.9205 \\
& 1/400 & 0.4249 & 0.4975 & 0.4481 \\
\end{tabular}
\caption{\revb{Percentage of area change in experiment \ref{sec:secexp3}: varying $\Delta t$ and $m$ with $N = 16$.}}
\label{tab:area3dt}
\end{center}
\end{table}

\begin{table}[ht]
\begin{center}
\begin{tabular}{c c|c  c  c  c}
& & \multicolumn{4}{c}{$m$} \\
& & 32 & 64 & 128 & 256 \\ \hline
\multirow{3}{*}{$\Delta t$}
& 1/100 & 1.6815 & 1.8129 & 1.8432 & 1.7833 \\
& 1/200 & 0.7912 & 0.9055 & 0.9724 & 0.8913 \\
& 1/400 & 0.3859 & 0.4286 & 0.4975 & 0.4644 \\
\end{tabular}
\caption{\revb{Percentage of area change in experiment \ref{sec:secexp3}: varying $\Delta t$ and $m$ with $N=m/8$.}}
\label{tab:area3dt2}
\end{center}
\end{table}

\begin{table}[ht]
\begin{center}
\begin{tabular}{c c|c c  c c}
& & \multicolumn{3}{c}{$\Delta t$} \\
& & 1/100 & 1/200 & 1/400 \\ \hline
\multirow{4}{*}{$N$}
& 4 & 2.4434 & 0.7953 & 0.4085 \\
& 8 & 1.8129 & 0.8878 & 0.4583 \\
& 16 & 1.8121 & 0.9724 & 0.4481 \\
& 32 & 0.2928 & 0.9187 & 0.4644 \\
\end{tabular}
\caption{\revb{Percentage of area change in experiment \ref{sec:secexp3}: varying $\Delta t$ and $N$ with $m=0.64/\Delta t$.}}
\label{tab:area3dt3}
\end{center}
\end{table}

\subsection{Stretched immersed boundary} 
\label{sec:secexp4}

\begin{figure}[ht!]
\centering
\includegraphics[width=2.5in]{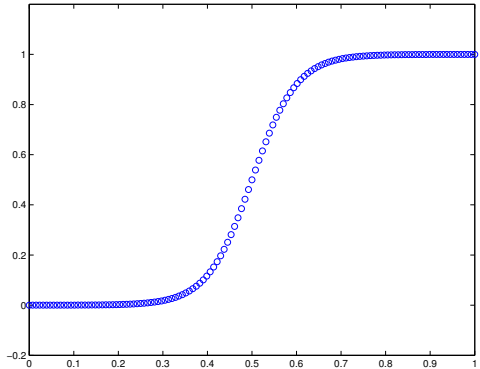}
\caption{Graph of the linearly rescaled sigmoid curve $\tilde{s} = \tilde{G}(s)$ in experiment \ref{sec:secexp4}.}
\label{fig:exp4a}
\end{figure}

We consider a model with the immersed boundary is initially stretched, i.e. the initial configuration of the Lagrangian markers is a non-uniformly spaced circle. We define a sigmoid function $G: [0,1] \to (0,1)$ by
\begin{equation}
G(s) = \frac{1}{1 + e^{-10+20s}}, \quad s \in [0,1].
\end{equation}
Let $\widetilde{G}$ be obtained by linearly rescaling the range of $G$ onto $[0,1]$. More precisely, $\widetilde{G}$ is defined as
\begin{equation}
\widetilde{G}(s) = \frac{G(s) - G(0)}{G(1) - G(0)}, \quad s \in [0,1].
\end{equation}
Figure \ref{fig:exp4a} shows the graph of the linearly rescaled sigmoid function $\tilde{s} = \widetilde{G}(s)$.

The initial condition for the fluid motion is given by
\begin{equation}
\bsym{u}_0(x,y) = \bsym{0}, \quad (x,y) \in [0,1]^2.
\end{equation}
The initial configuration of the Lagrangian markers is given by
\begin{equation}
\bsym{X}_0(s) =
\begin{pmatrix} 0.2\cos(2\pi \tilde{G}(s)) + 0.5\\ 0.2 \sin(2\pi \tilde{G}(s)) + 0.5\end{pmatrix},
\quad s \in [0,1].
\end{equation}

In the non-uniform parametrization, some markers are farther away from their neighbours. The longer distance between a particle and its neighbouring particles has a higher tension and models a stretched portion of the curve. In this experiment, the immersed boundary is stretched at an interval around $s = 0.5$.

Tests are performed with mesh sizes $N = 4, 8, 16, 32$ and $m = 64, 128, 256$ and $K = 200$. At $t = 2$, the area change is analyzed.

\begin{figure}[ht!]
\centering
\includegraphics[width=2.5in]{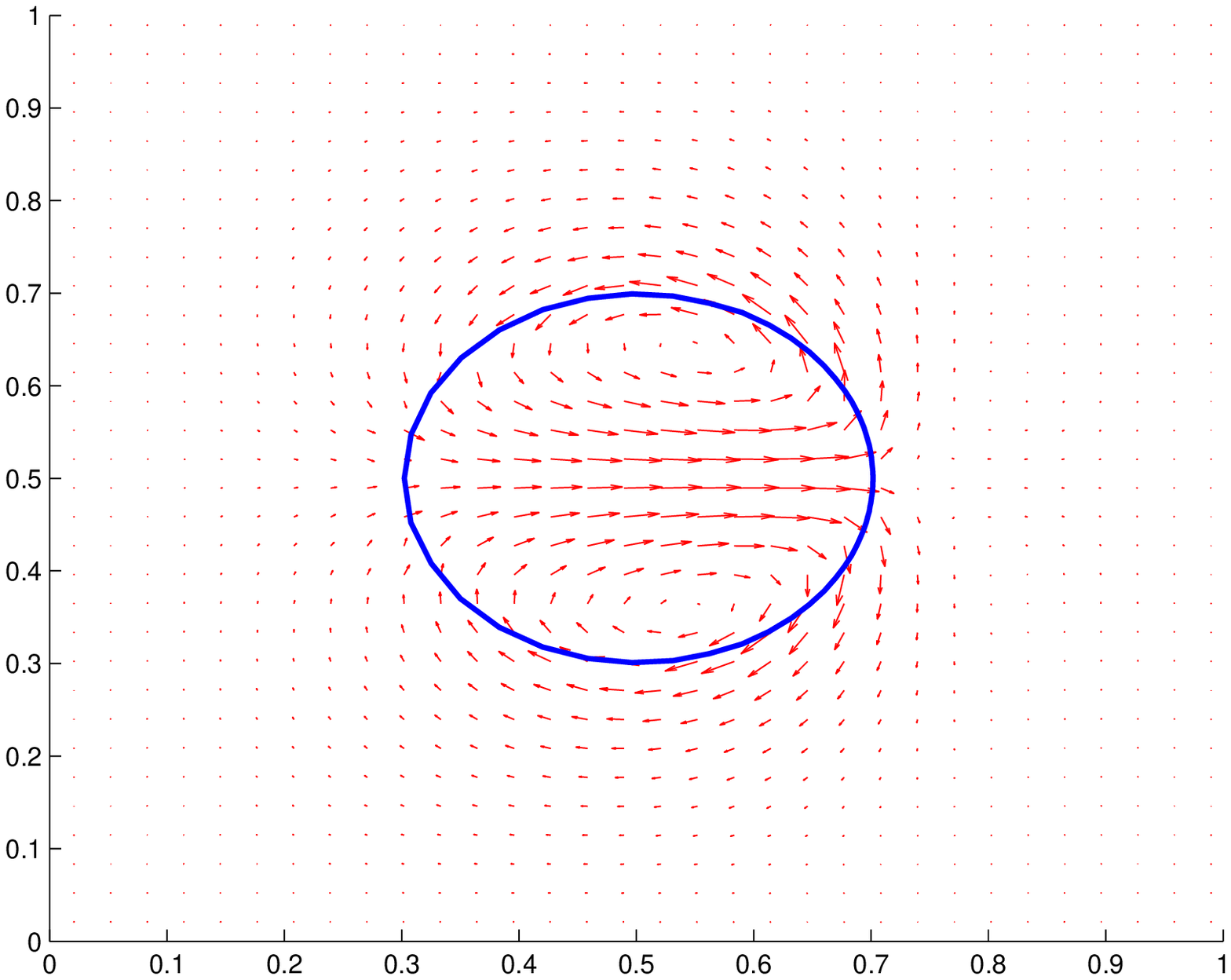}
\caption{Profile of the fluid flow and configuration of the Lagrangian markers at $t = \Delta t$ in experiment \ref{sec:secexp4}.}
\label{fig:exp4b}
\end{figure}

\begin{figure}[ht!]
\centering
\includegraphics[width=2.5in]{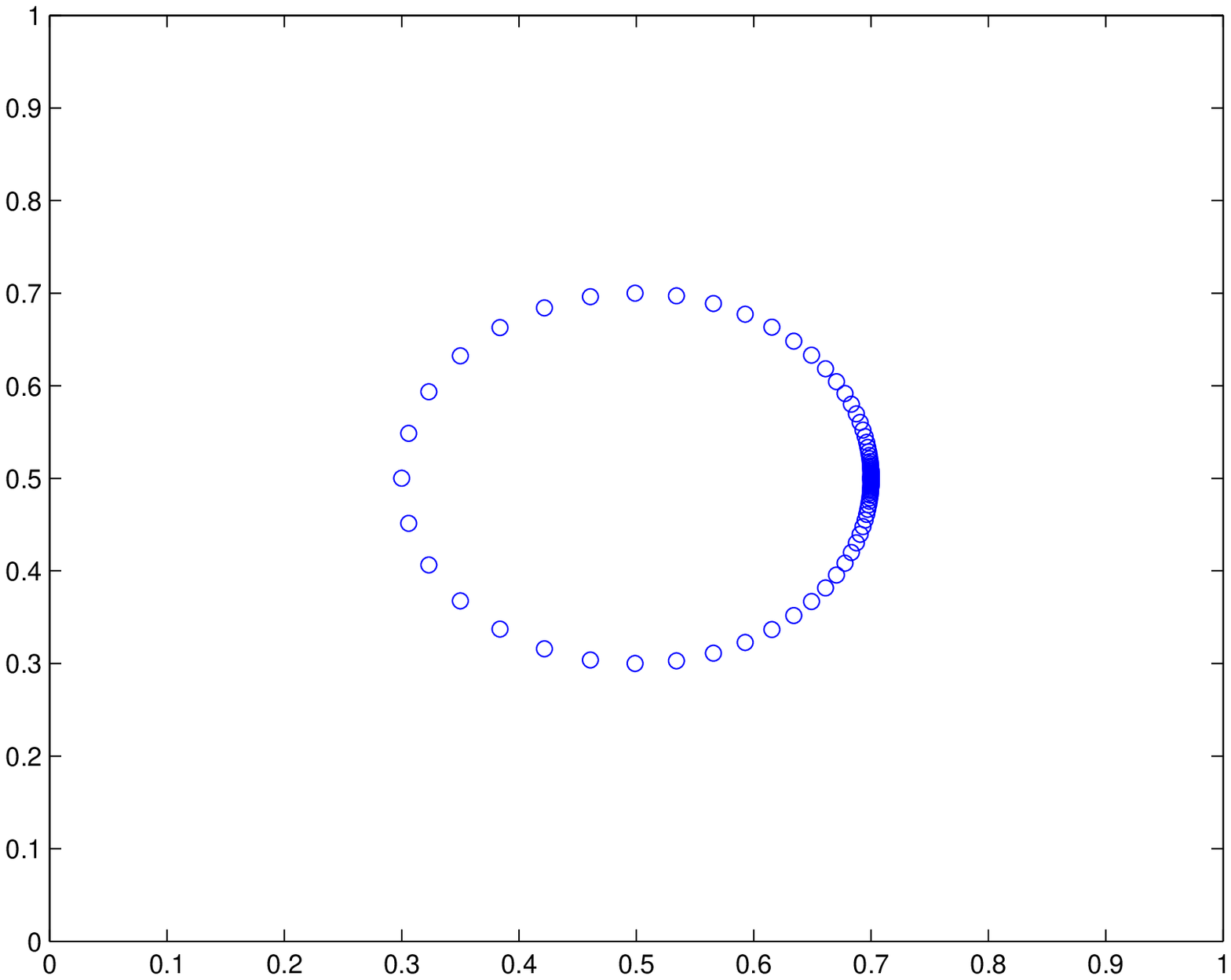}
\caption{Configuration of the Lagrangian markers $\bsym{X}(s,t)$ at $t = 0$ in experiment \ref{sec:secexp4}.}
\label{fig:exp4c}
\end{figure}

\begin{figure}[ht!]
\centering
\includegraphics[width=2.5in]{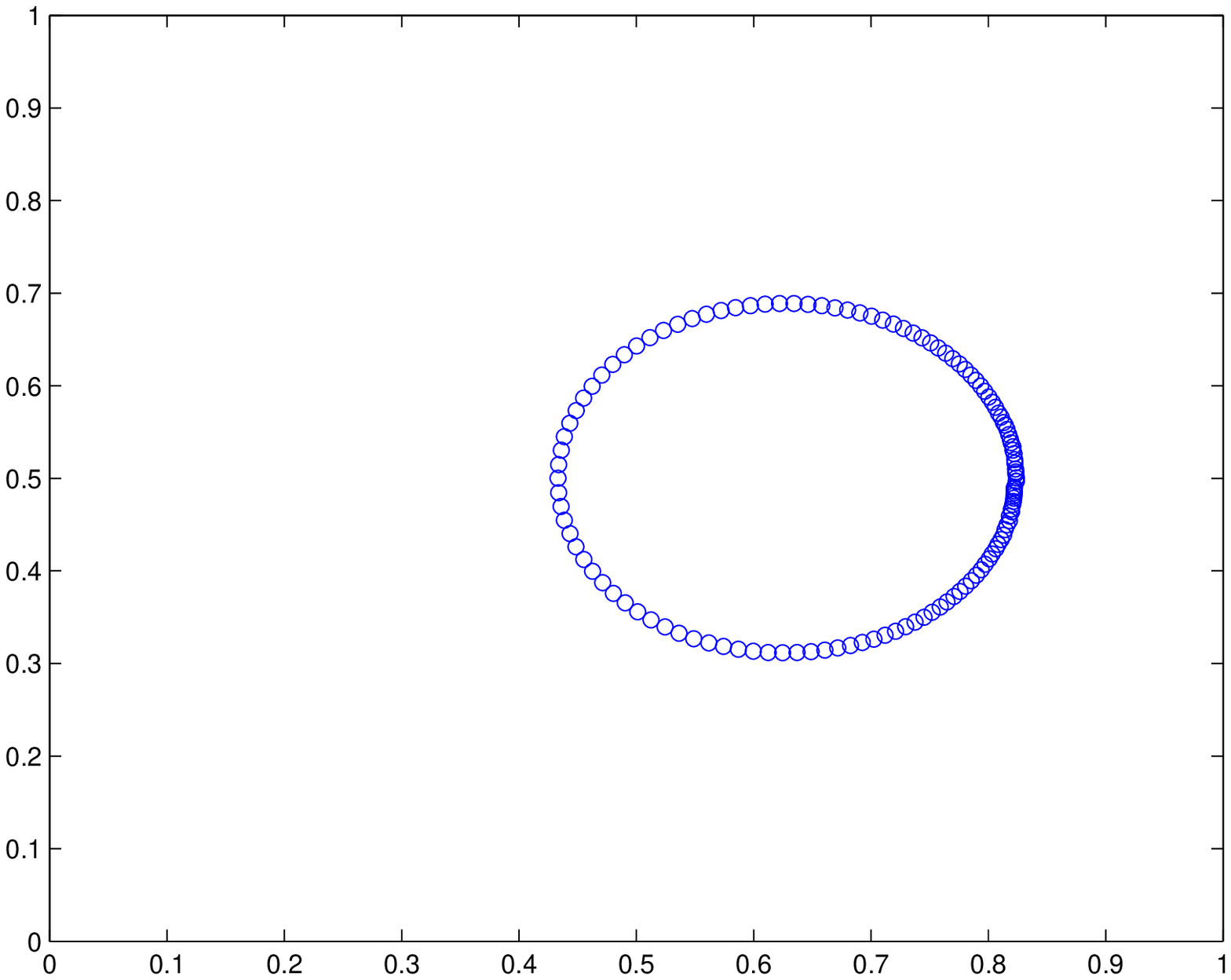}
\caption{Configuration of the Lagrangian markers $\bsym{X}(s,t)$ at $t = 1$ in experiment \ref{sec:secexp4}.}
\label{fig:exp4d}
\end{figure}

\begin{figure}[ht!]
\centering
\includegraphics[width=2.5in]{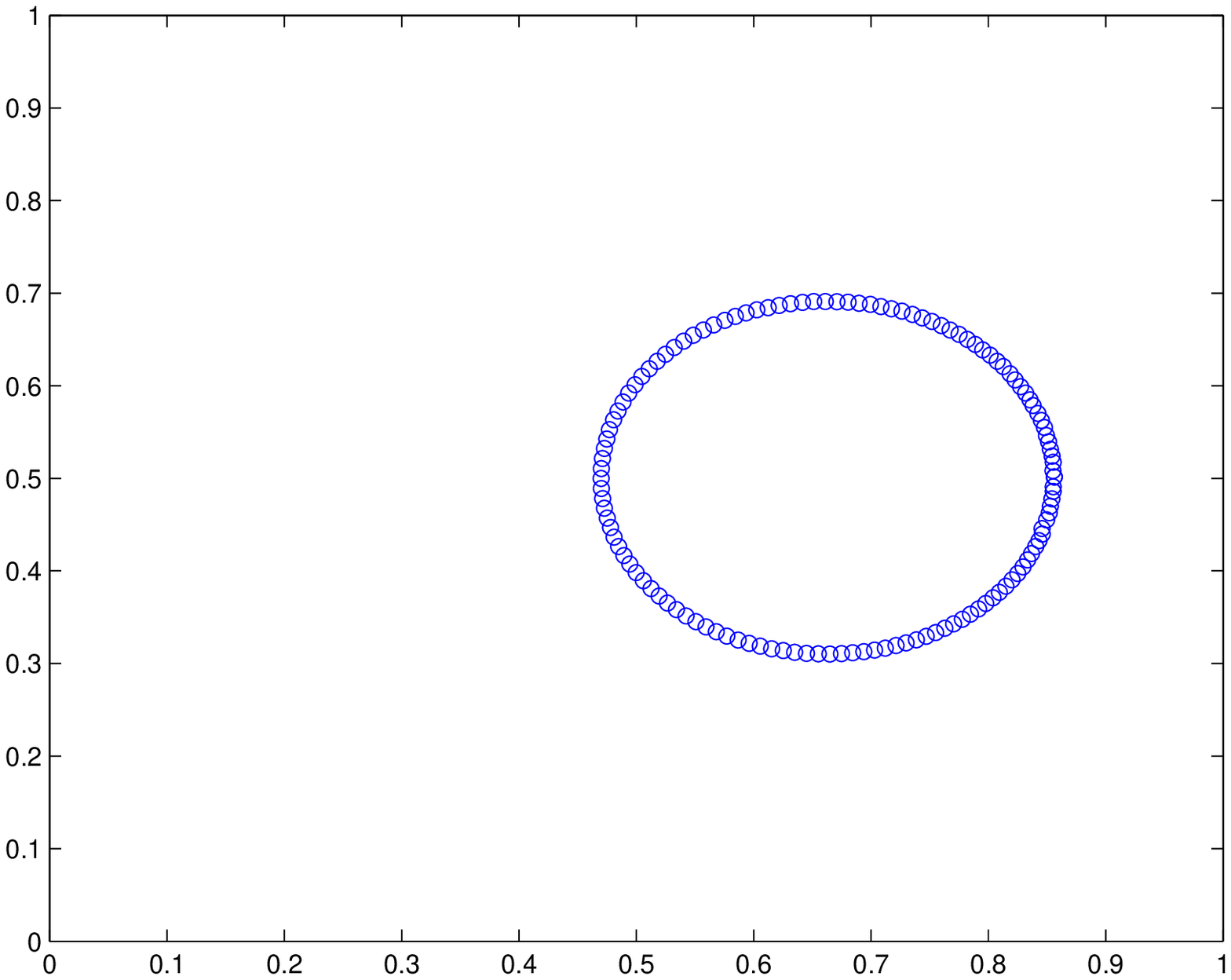}
\caption{Configuration of the Lagrangian markers $\bsym{X}(s,t)$ at $t = 2$ in experiment \ref{sec:secexp4}.}
\label{fig:exp4e}
\end{figure}

\begin{figure}[ht!]
\centering
\includegraphics[width=2.5in]{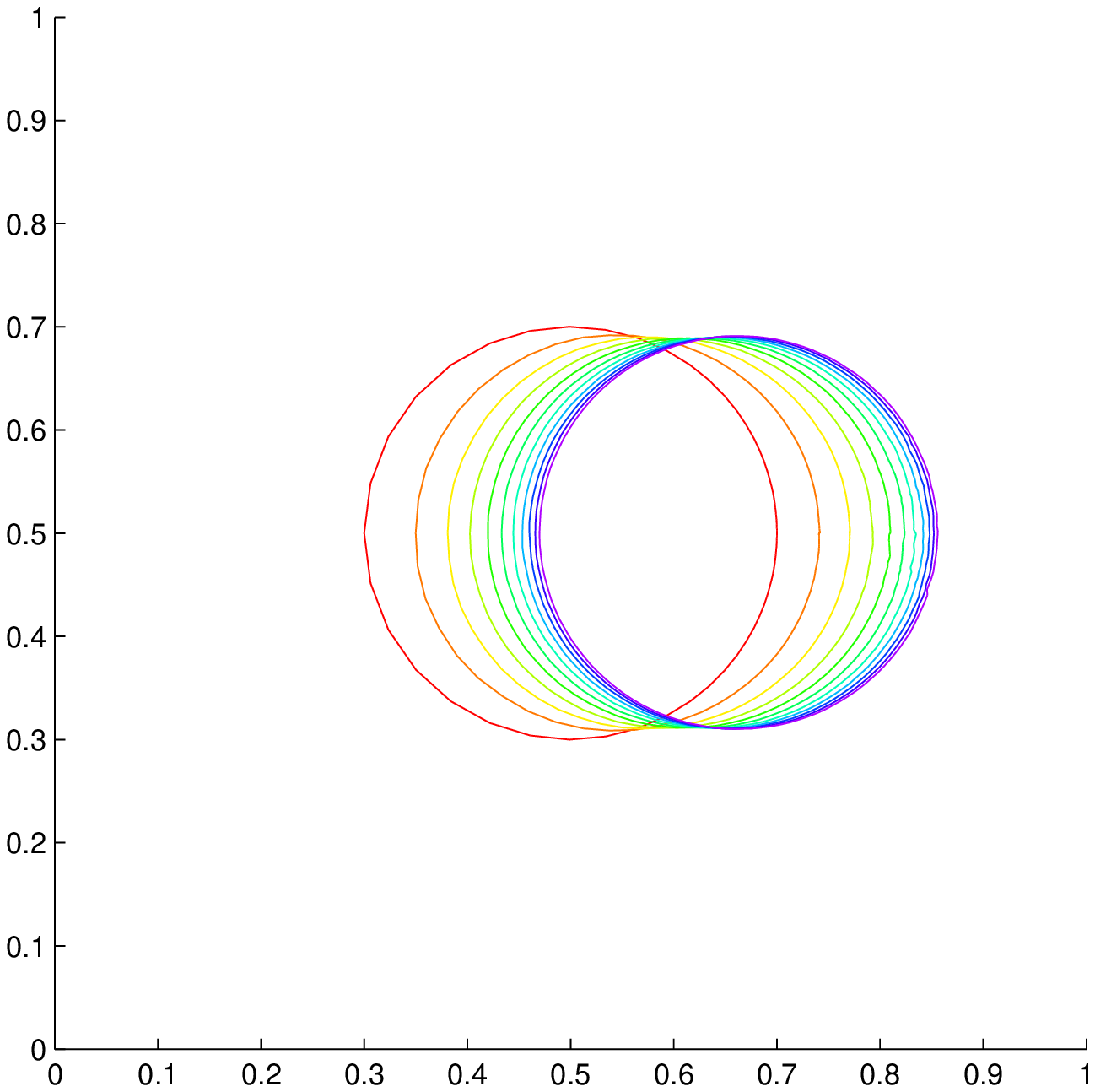}
\caption{Evolution of the immersed boundary throughout $t = 0$ and $t = 2$ in experiment \ref{sec:secexp4}.}
\label{fig:exp4f}
\end{figure}

Table \ref{tab:area4} records the area change of the immersed boundary in experiment \ref{sec:secexp4}. It can be observed that the area conservation depends heavily on a balance in the number of divisions $N$ in Eulerian mesh and $m$ in Lagrangian mesh. Figure \ref{fig:exp4b} shows the profile of the fluid flow and the configuration of the Lagrangian markers at $t = \Delta t$ with $N = 32$ and $m = 256$. It can be observed that the fluid flows into the immersed boundary at the scretched portion. The flow substantially pushes the immersed boundary in the direction away from the stretched portion. Figures \ref{fig:exp4c}--\ref{fig:exp4e} show the configurations of the immersed boundary $\bsym{X}(s,t)$ at different time $t$ with $N = 32$ and $m = 256$. Figure \ref{fig:exp4f} shows the evolution of the immersed boundary throughout $t = 0$ and $t = 2$ with $N = 32$ and $m = 256$. It can be seen that the Lagrangian markers tend to the circular equilibrium configuration and become evenly spaced.

\begin{table}[ht]
\begin{center}
\begin{tabular}{c c|c c c}
& & \multicolumn{3}{c}{$m$} \\
& & 64 & 128 & 256 \\ \hline
\multirow{4}{*}{$N$}
& 4 & 2.9441 & 0.7611 & 21.3877 \\
& 8 & -0.1162 & 1.6579 & 2.5126 \\
& 16 & -15.8600 & -0.3102 & 1.6067 \\
& 32 & -30.6409 & -7.7084 & 0.6536 \\
\end{tabular}
\caption{Percentage of area change in experiment \ref{sec:secexp4}.}
\label{tab:area4}
\end{center}
\end{table}

\subsection{Numerical stability}
\label{sec:secexp5}
The last numerical experiment is devoted to inspecting the numerical stability of SDG-IBM. According to the results in Section \ref{sec:stab}, if $\eta^n$ is sufficiently small, then the method would be stable and the energy would not blow up.

We consider a model proposed in \cite{boffi07}. A balloon with radius $R$ is inflated and placed at rest in the middle of a square domain $[0,1]^2$ filled with fluid. The initial condition for the fluid motion is given by
\begin{equation}
\bsym{u}_0(x,y) = \bsym{0}, \quad (x,y) \in [0,1]^2.
\end{equation}
The initial configuration of the Lagrangian markers is given by
\begin{equation}
\bsym{X}_0(s) =
\begin{pmatrix} R\cos(s/R) + 0.5\\ R \sin(s/R) + 0.5\end{pmatrix},
\quad s \in [0,2 \pi R].
\end{equation}
In this experiment, we set $R = 0.4$.

Tests are performed with mesh sizes $N = 32$ and $m = 128$ and $K = 120, 300, 600$. The elasticity is set to be $\kappa = 1, 2, 4$. Throughout $t = 0$ to $t = 3$, the quantities $E^n$ and $\eta^n$ are analyzed. The parameters are chosen in order to compare our method with FE-IBM in \cite[Fig. 3]{boffi07}.

\begin{figure}[ht!]
\centering
\begin{subfigure}[b]{0.3\textwidth}
\includegraphics[width=\textwidth]{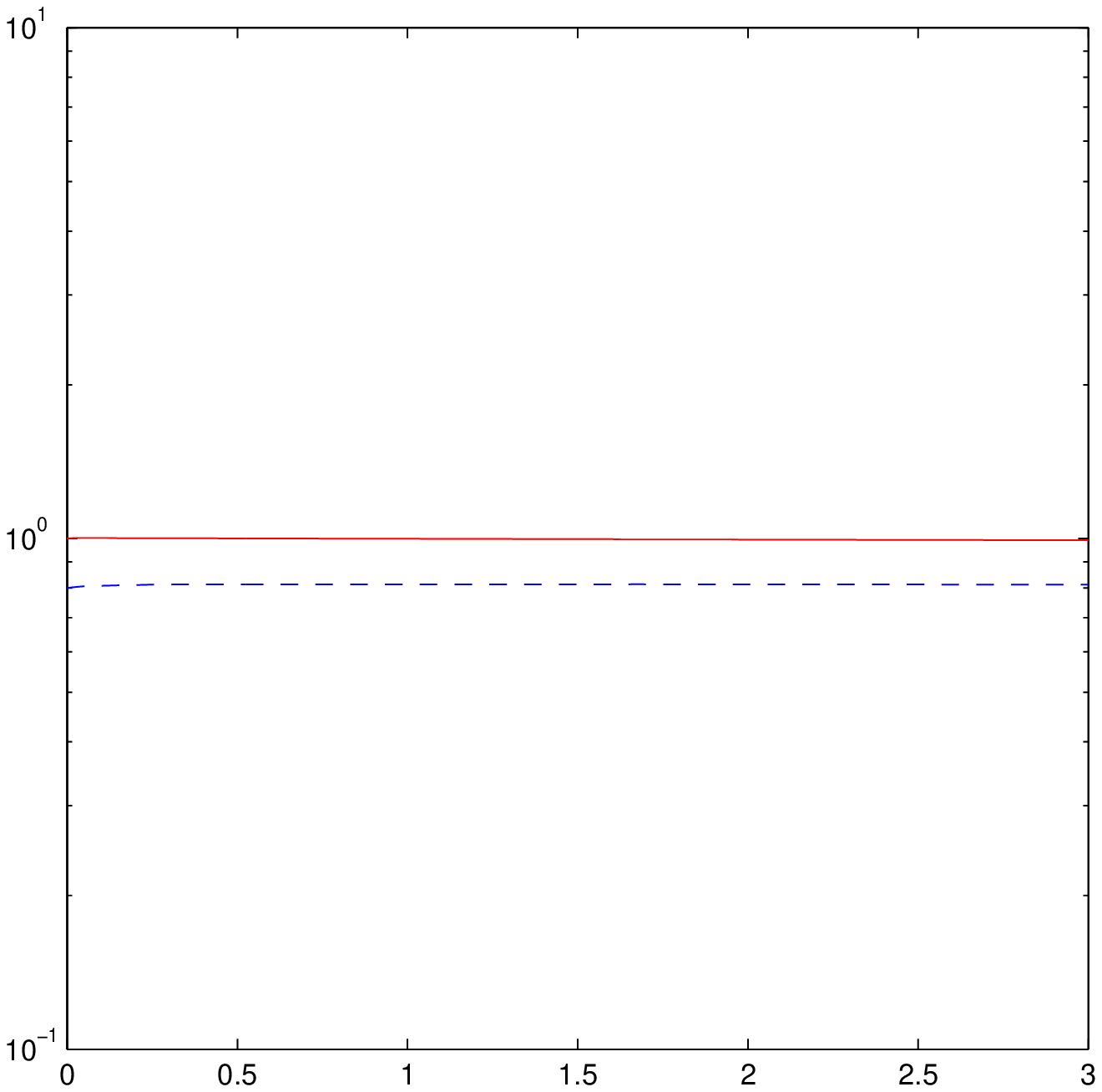}
\caption{$K = 120, \kappa = 1$}
\end{subfigure}
\begin{subfigure}[b]{0.3\textwidth}
\includegraphics[width=\textwidth]{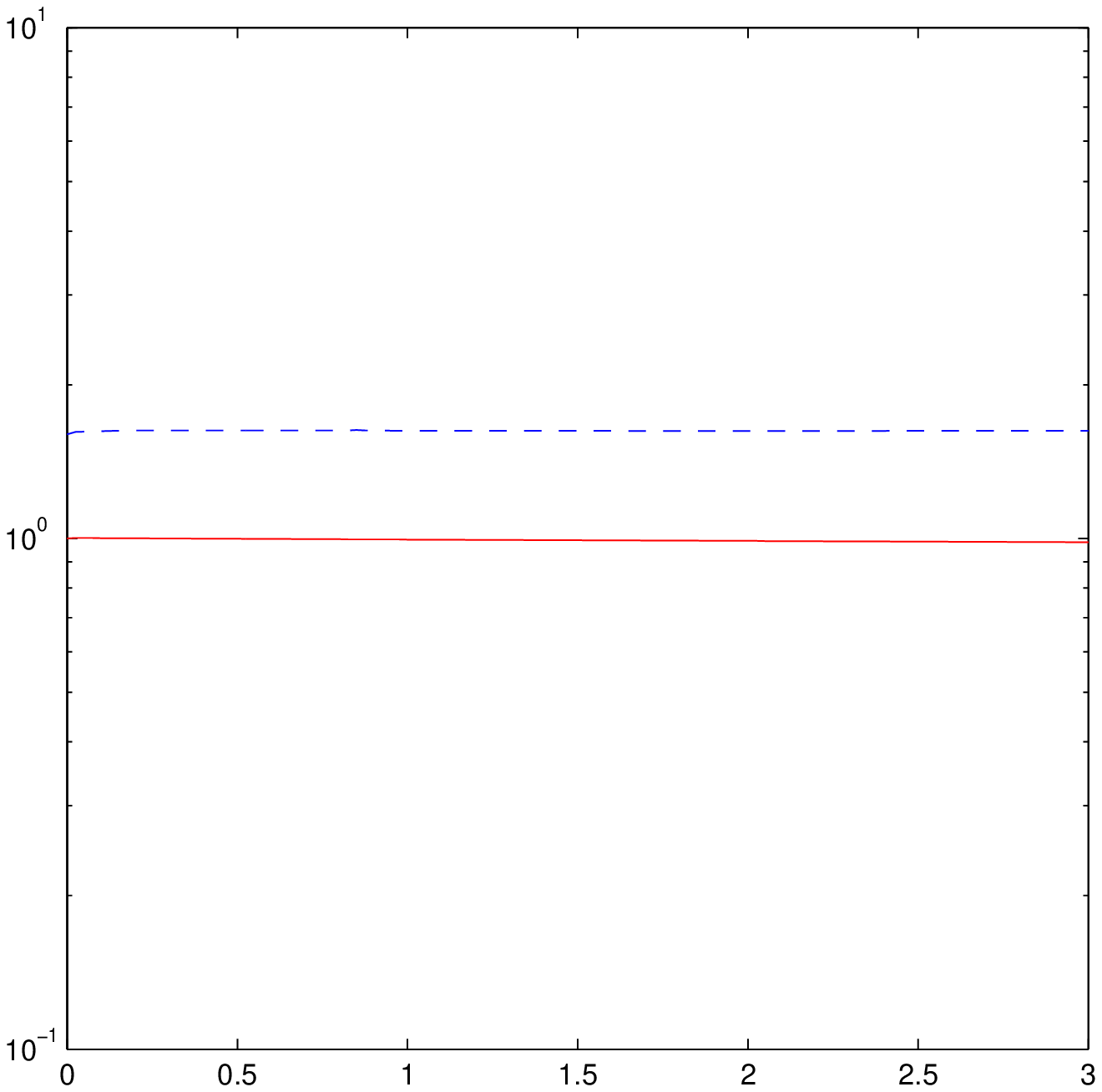}
\caption{$K = 120, \kappa = 2$}
\end{subfigure}
\begin{subfigure}[b]{0.3\textwidth}
\includegraphics[width=\textwidth]{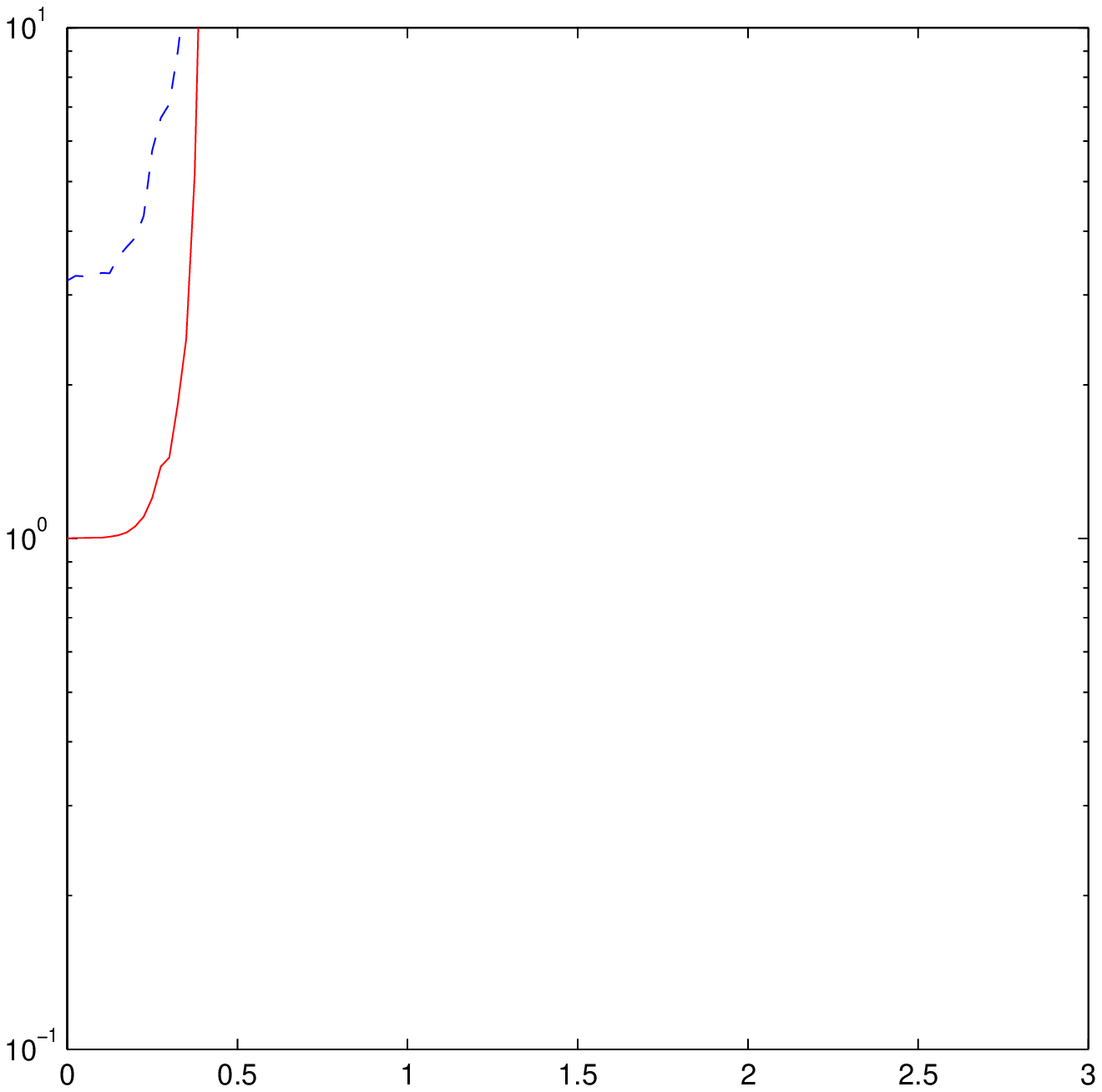}
\caption{$K = 120, \kappa = 4$}
\end{subfigure}
\begin{subfigure}[b]{0.3\textwidth}
\includegraphics[width=\textwidth]{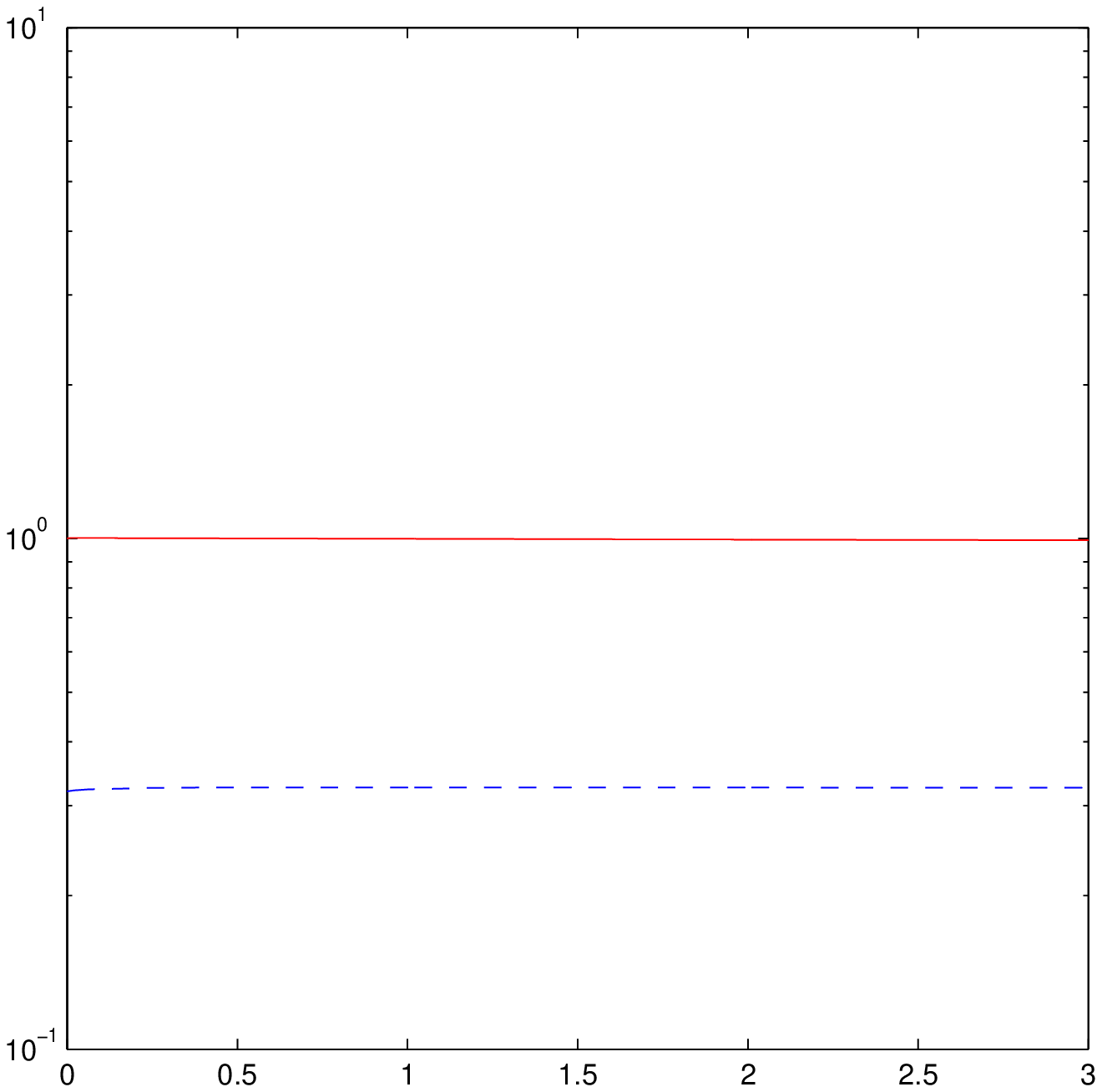}
\caption{$K = 300, \kappa = 1$}
\end{subfigure}
\begin{subfigure}[b]{0.3\textwidth}
\includegraphics[width=\textwidth]{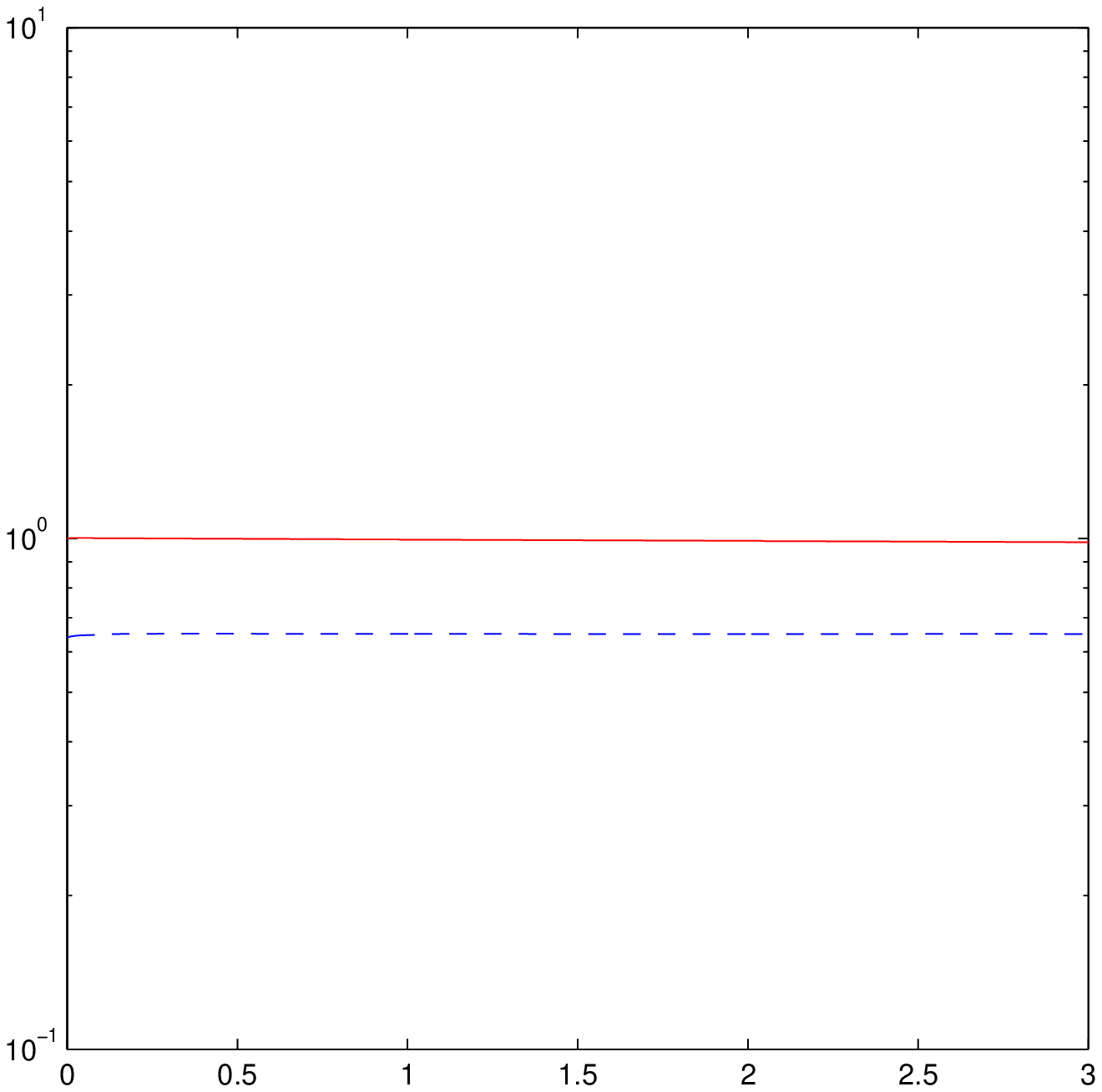}
\caption{$K = 300, \kappa = 2$}
\end{subfigure}
\begin{subfigure}[b]{0.3\textwidth}
\includegraphics[width=\textwidth]{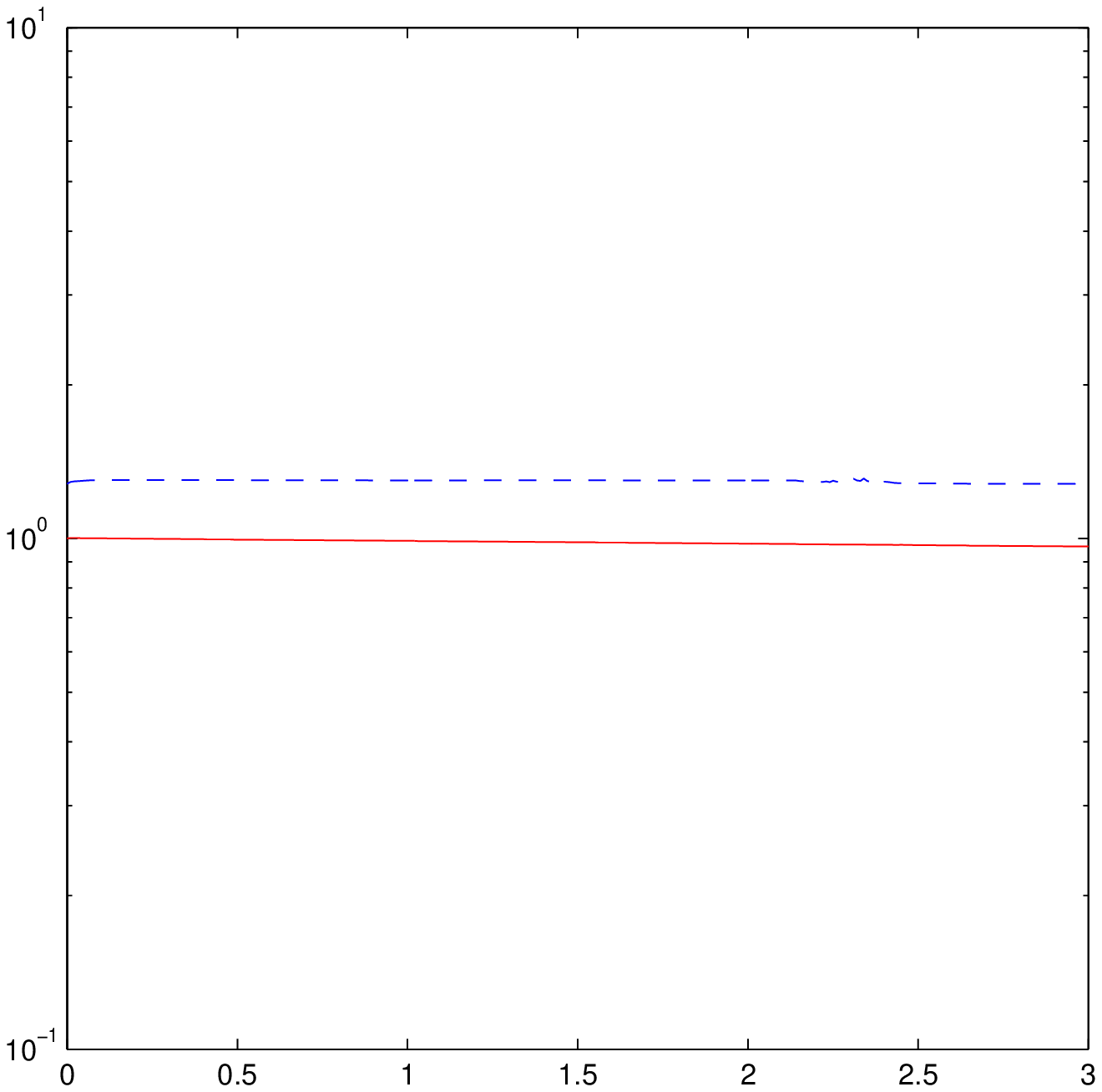}
\caption{$K = 300, \kappa = 4$}
\end{subfigure}
\begin{subfigure}[b]{0.3\textwidth}
\includegraphics[width=\textwidth]{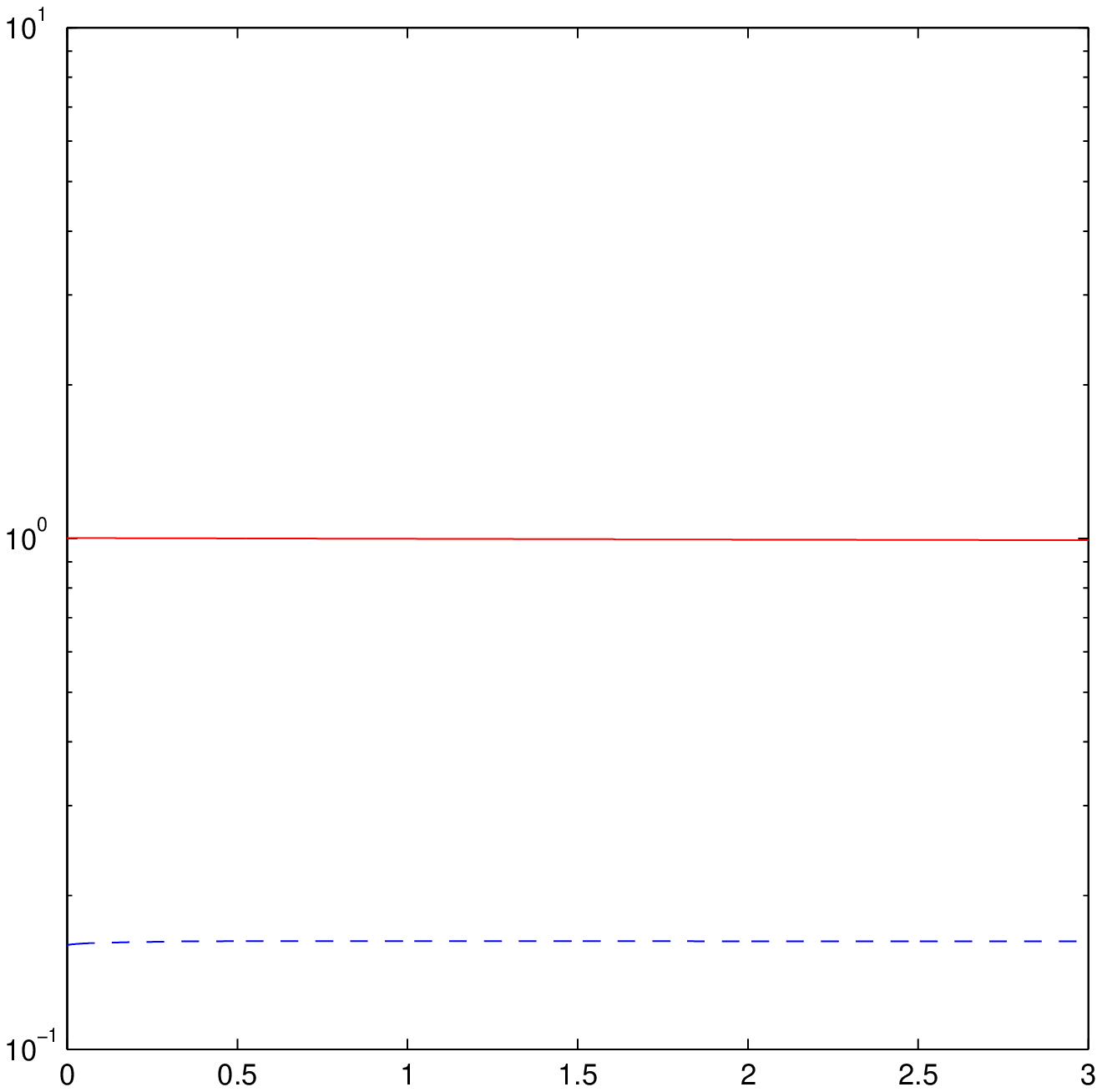}
\caption{$K = 600, \kappa = 1$}
\end{subfigure}
\begin{subfigure}[b]{0.3\textwidth}
\includegraphics[width=\textwidth]{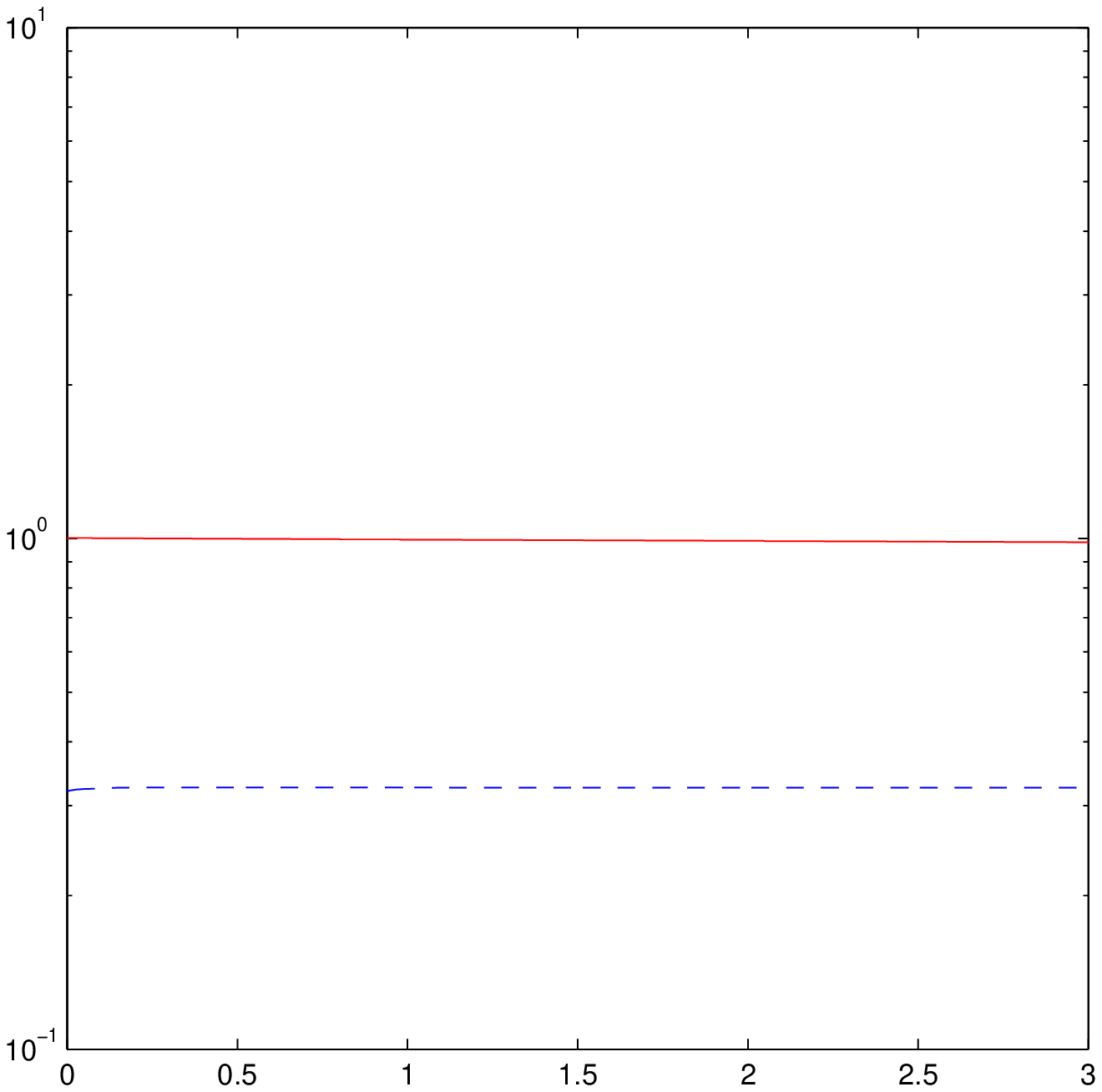}
\caption{$K = 600, \kappa = 2$}
\end{subfigure}
\begin{subfigure}[b]{0.3\textwidth}
\includegraphics[width=\textwidth]{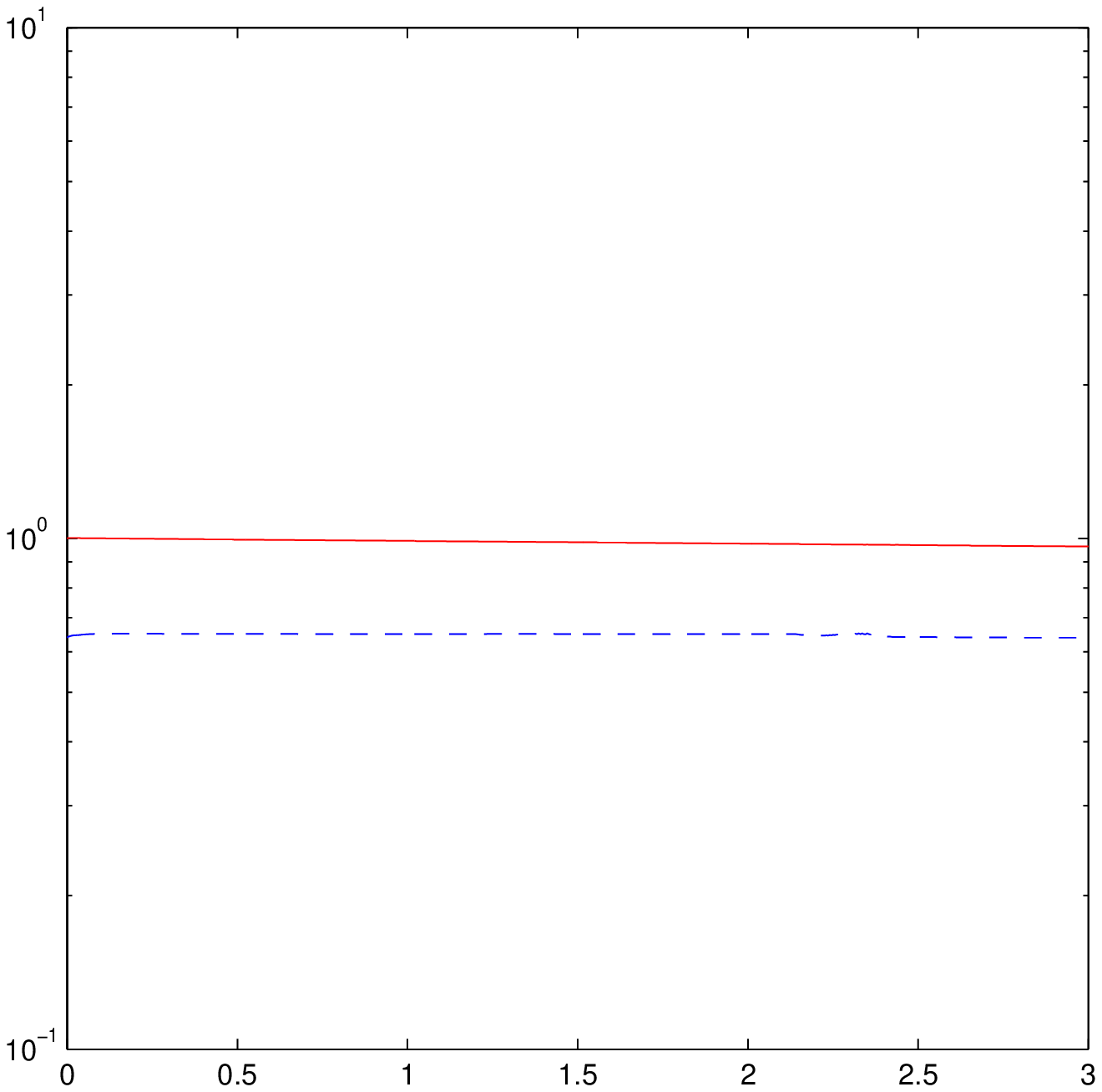}
\caption{$K = 600, \kappa = 4$}
\end{subfigure}
\caption{\revb{History of $E^n$ (red solid line) and $\eta^n$ (blue dashed line) throughout $t = 0$ and $t = 3$ with $m=128$ and $N=32$ in experiment \ref{sec:secexp5}.}}
\label{fig:exp5}
\end{figure}

\begin{figure}[ht!]
\centering
\begin{subfigure}[b]{0.3\textwidth}
\includegraphics[width=\textwidth]{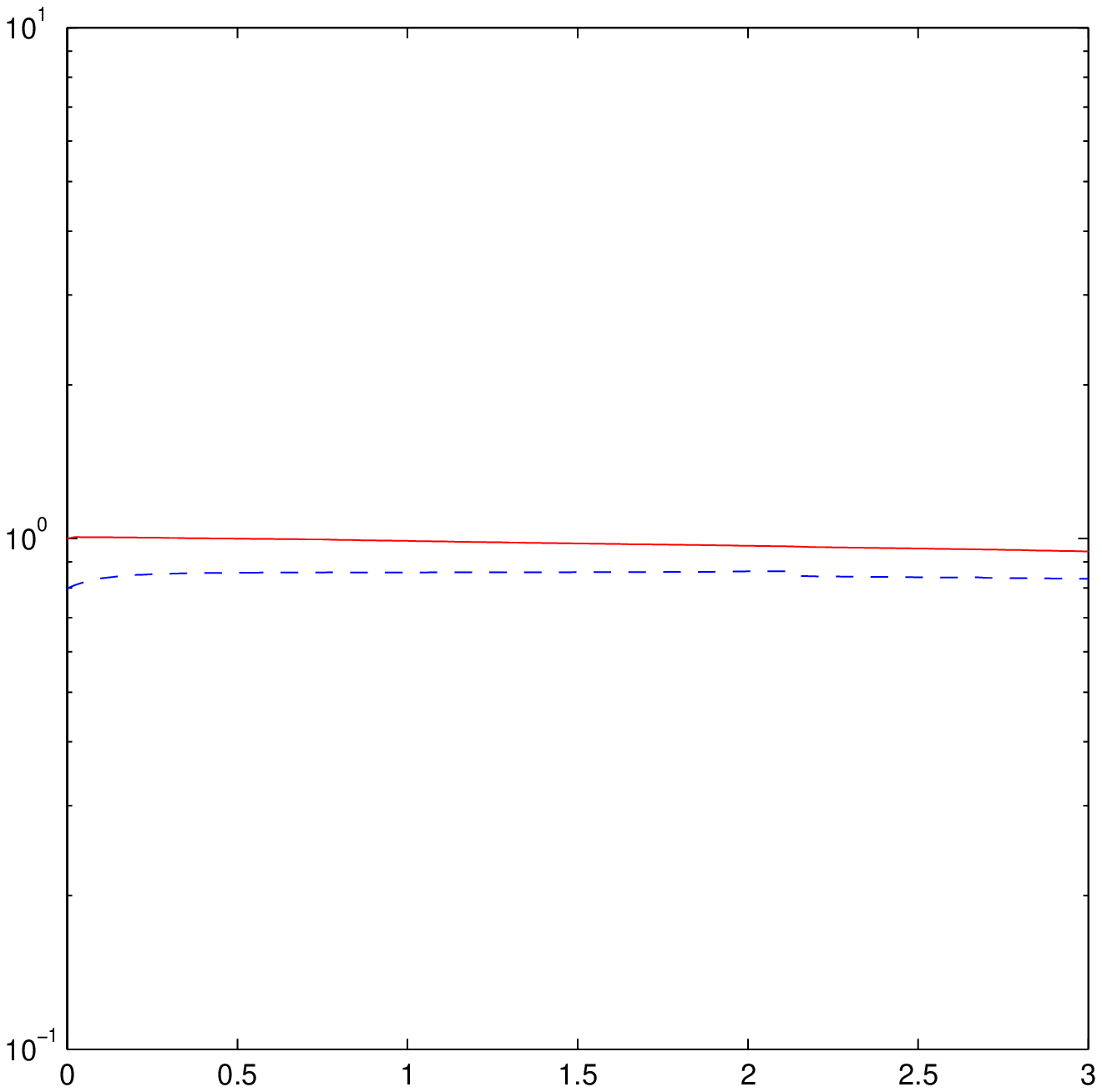}
\caption{$K = 120, m = 32, N = 8$}
\end{subfigure}
\begin{subfigure}[b]{0.3\textwidth}
\includegraphics[width=\textwidth]{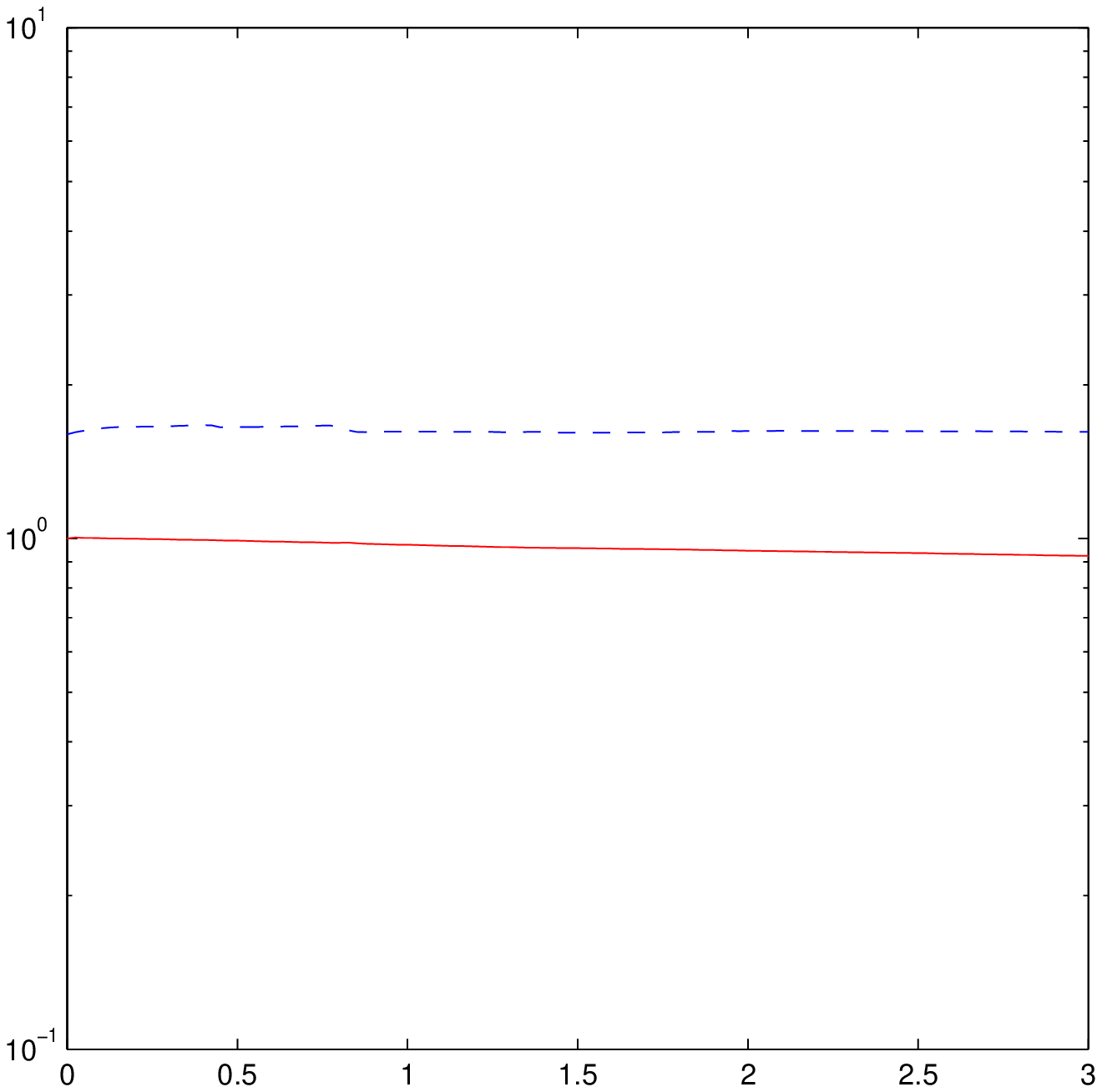}
\caption{$K = 120, m = 64, N = 16$}
\end{subfigure}
\begin{subfigure}[b]{0.3\textwidth}
\includegraphics[width=\textwidth]{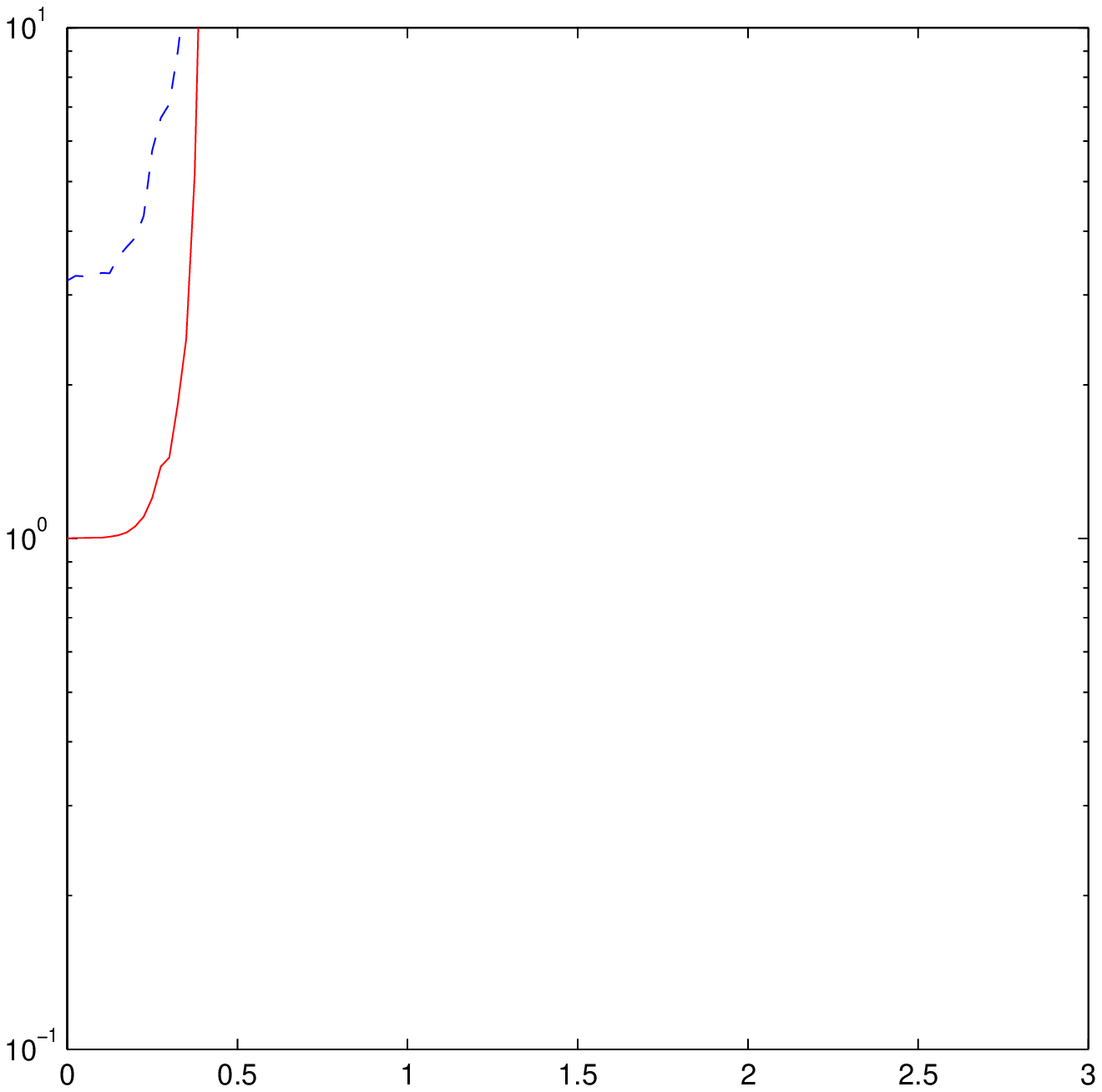}
\caption{$K = 120, m = 128, N = 32$}
\end{subfigure}
\begin{subfigure}[b]{0.3\textwidth}
\includegraphics[width=\textwidth]{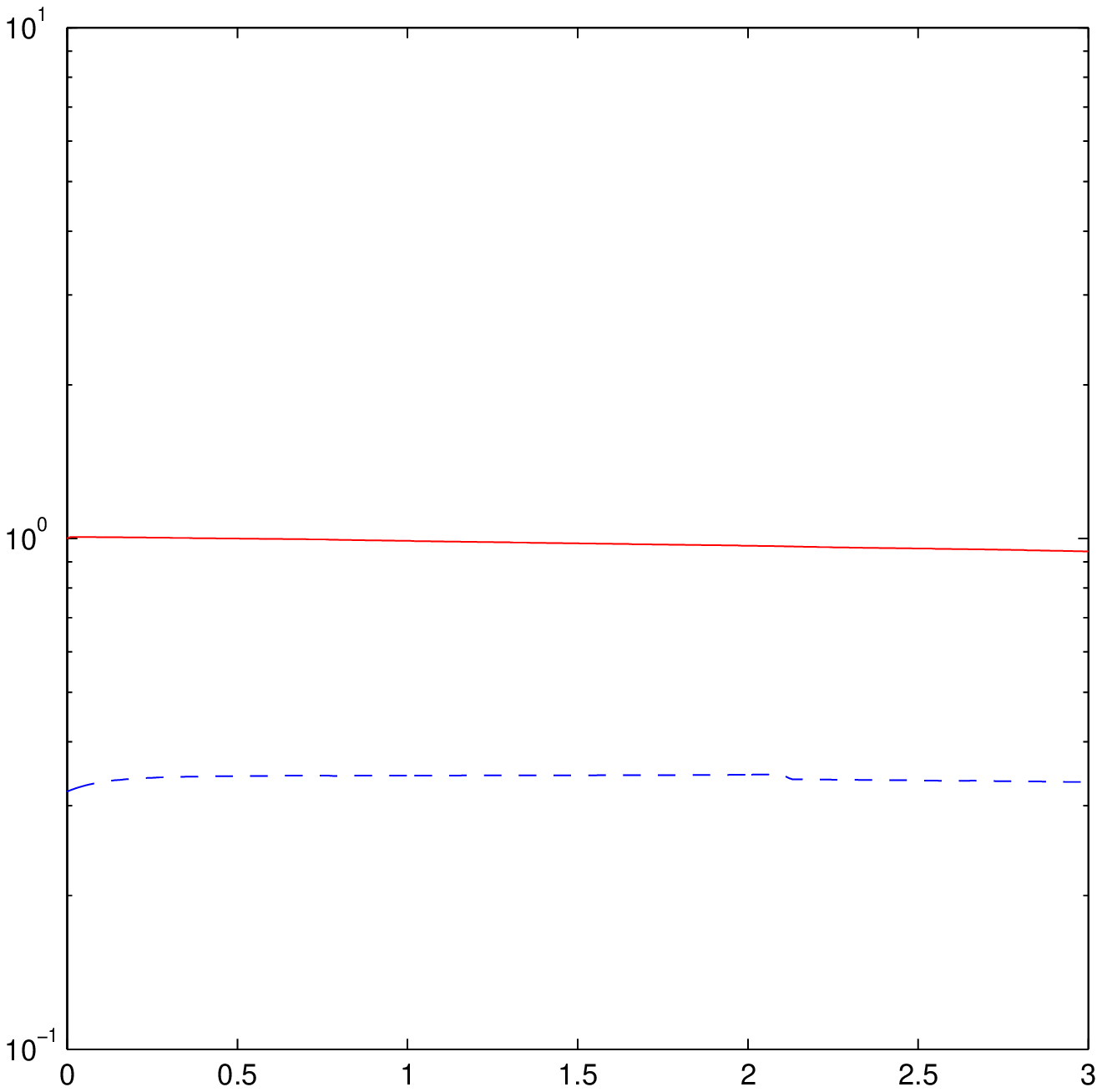}
\caption{$K = 300, m = 32, N = 8$}
\end{subfigure}
\begin{subfigure}[b]{0.3\textwidth}
\includegraphics[width=\textwidth]{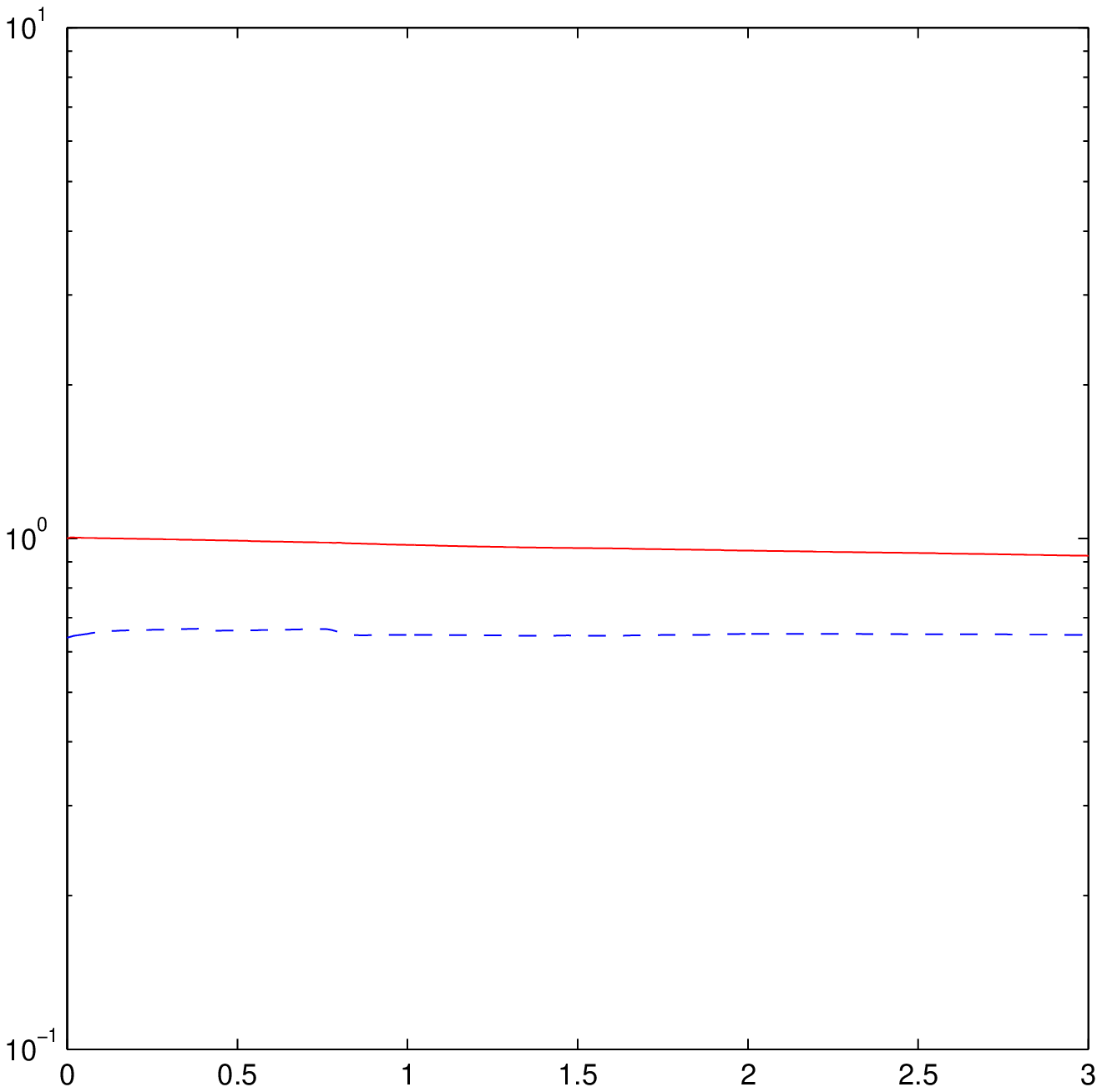}
\caption{$K = 300, m = 64, N = 16$}
\end{subfigure}
\begin{subfigure}[b]{0.3\textwidth}
\includegraphics[width=\textwidth]{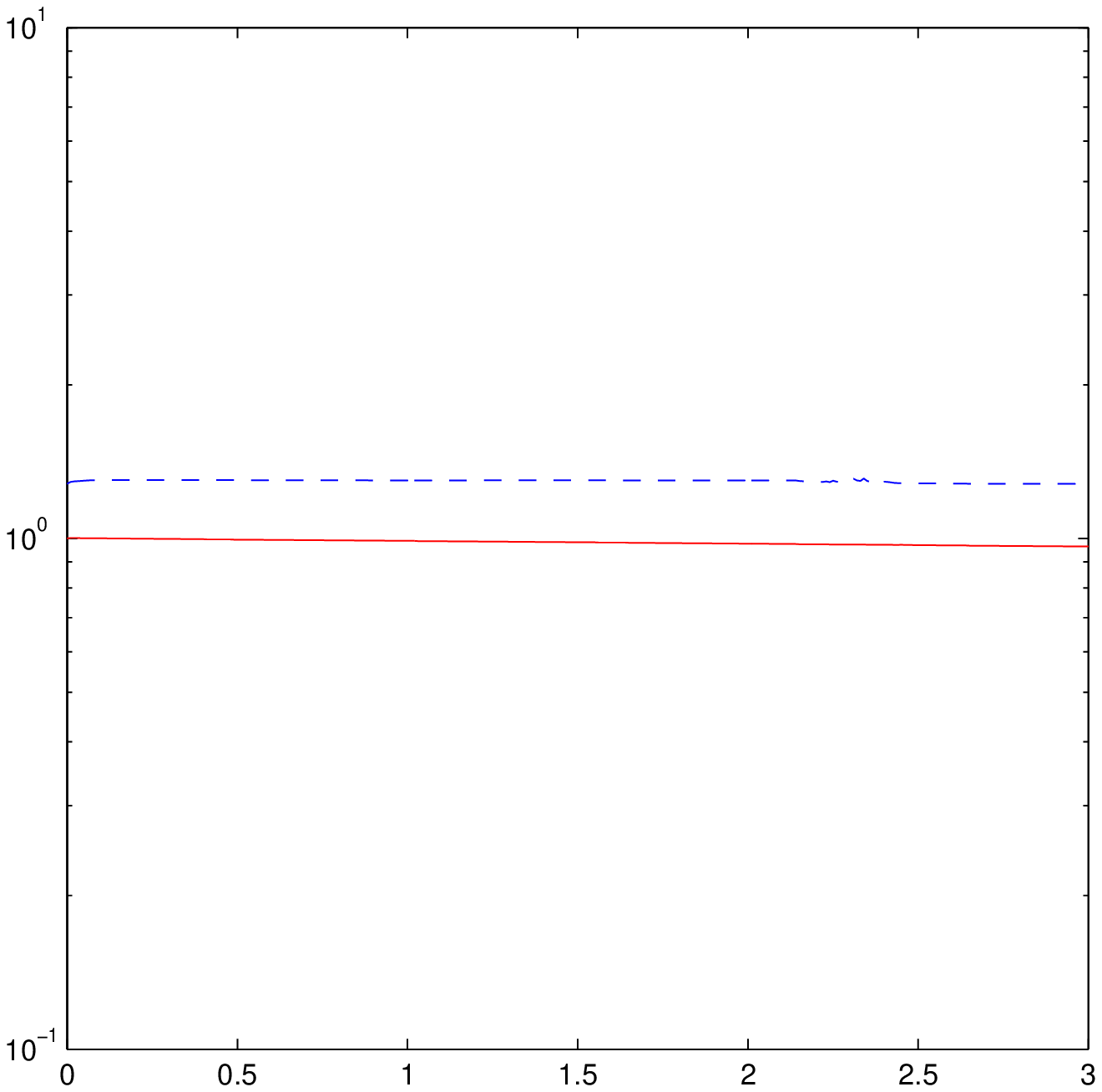}
\caption{$K = 300, m = 128, N = 32$}
\end{subfigure}
\begin{subfigure}[b]{0.3\textwidth}
\includegraphics[width=\textwidth]{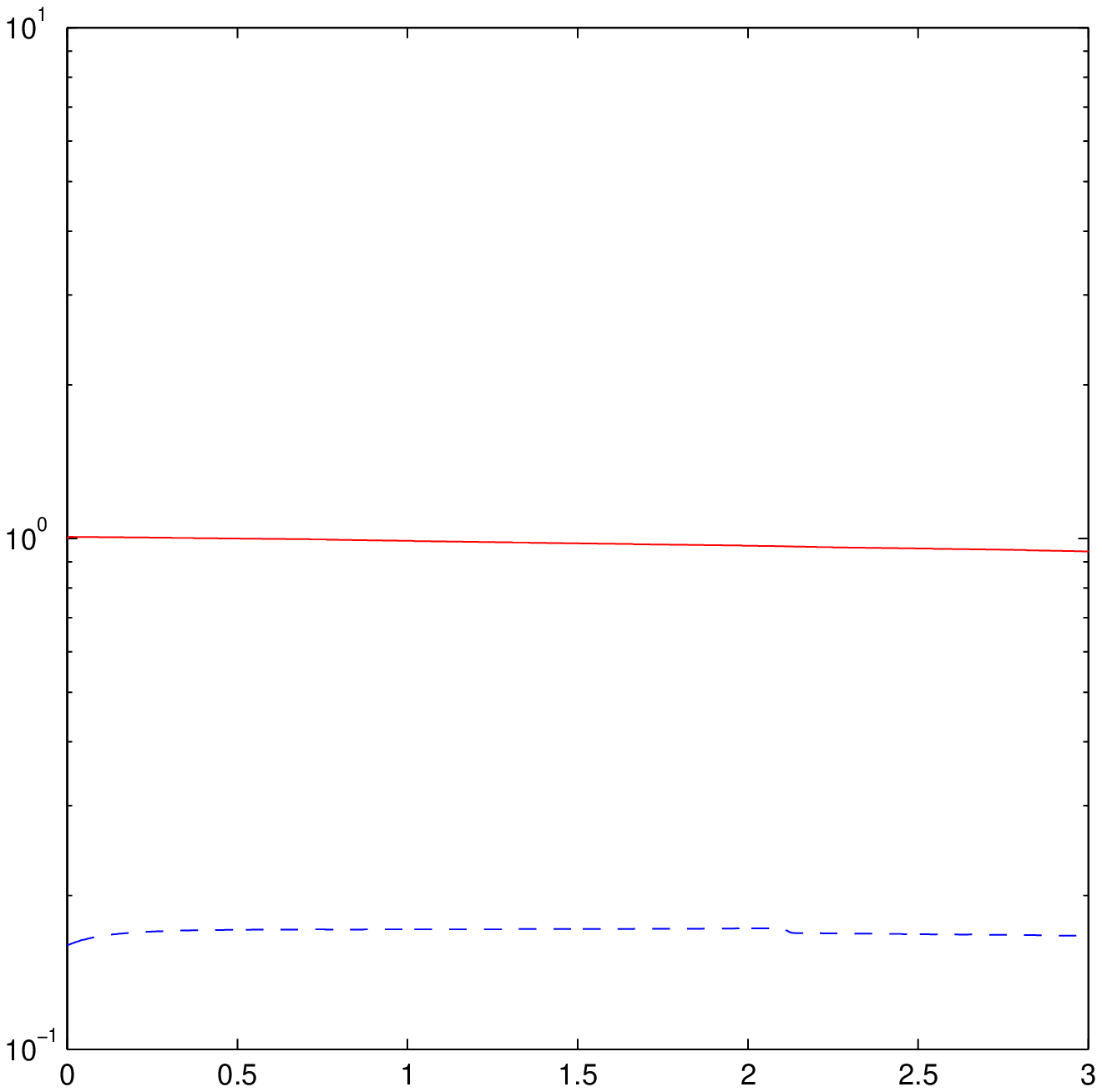}
\caption{$K = 600, m = 32, N = 8$}
\end{subfigure}
\begin{subfigure}[b]{0.3\textwidth}
\includegraphics[width=\textwidth]{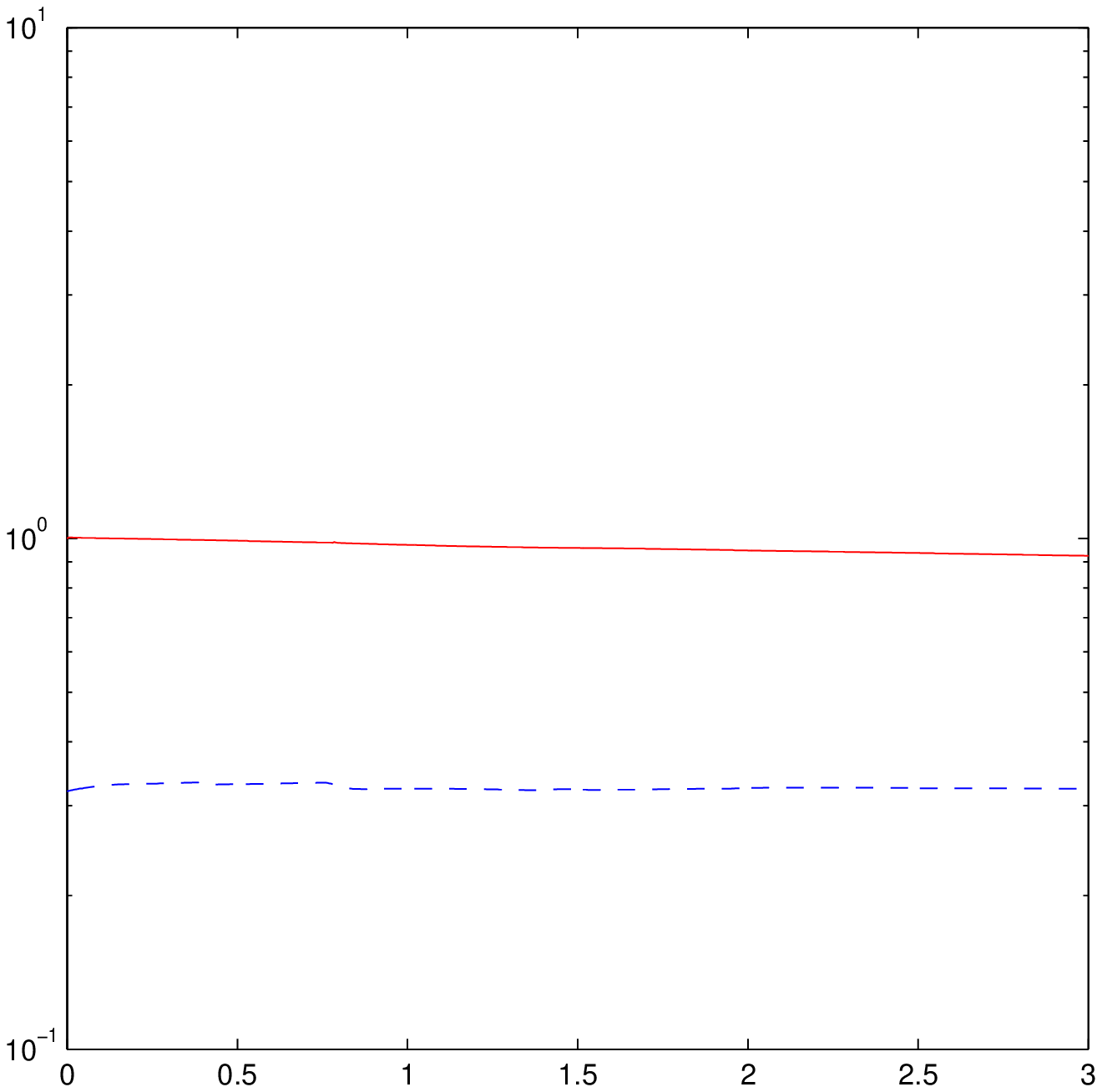}
\caption{$K = 600, m = 64, N = 16$}
\end{subfigure}
\begin{subfigure}[b]{0.3\textwidth}
\includegraphics[width=\textwidth]{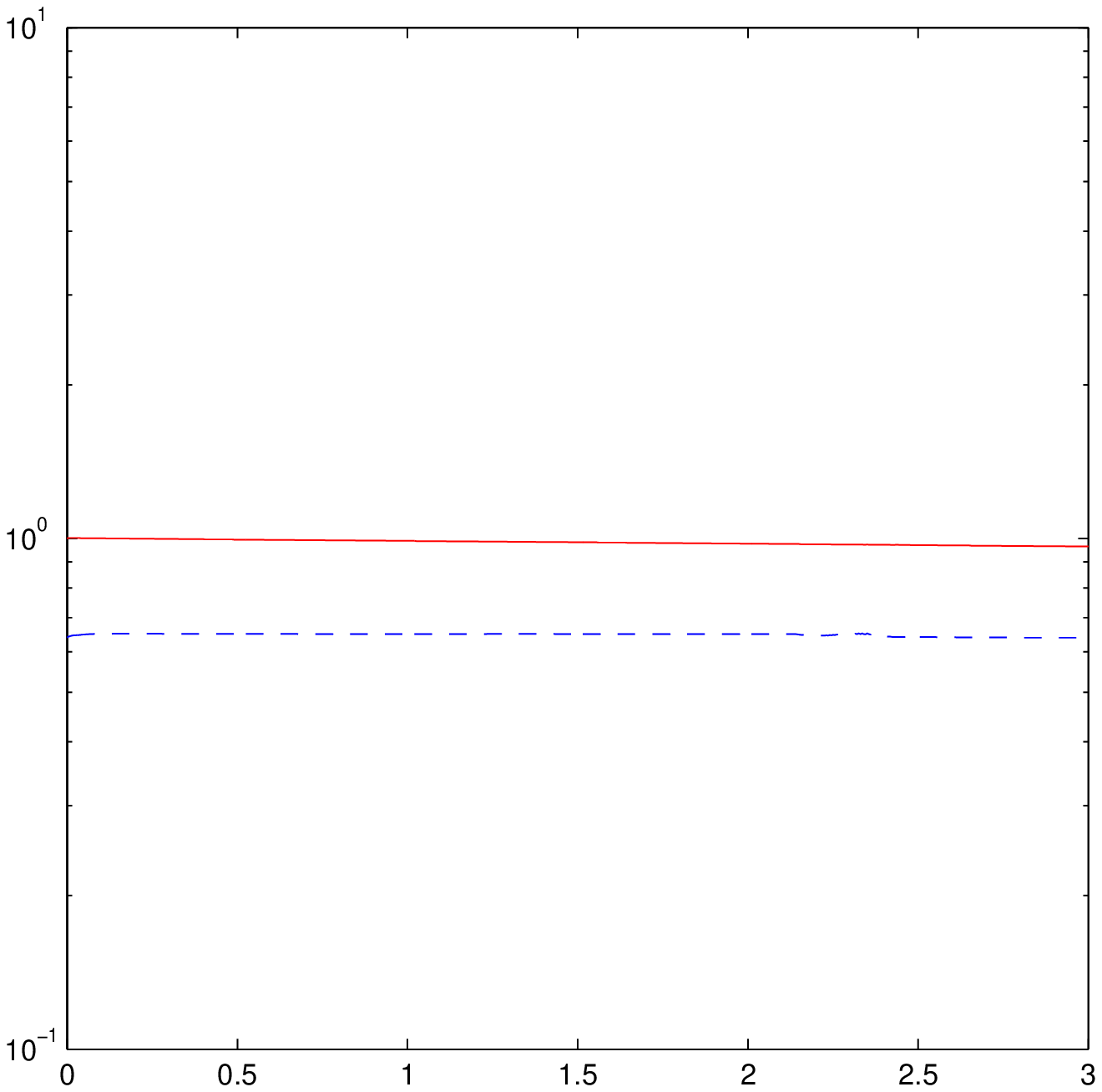}
\caption{$K = 600, m = 128, N = 32$}
\end{subfigure}
\caption{\revb{History of $E^n$ (red solid line) and $\eta^n$ (blue dashed line) throughout $t = 0$ and $t = 3$ with $\kappa=4$ in experiment \ref{sec:secexp5}.}}
\label{fig:exp5b}
\end{figure}

\begin{figure}[ht!]
\centering
\begin{subfigure}[b]{0.3\textwidth}
\includegraphics[width=\textwidth]{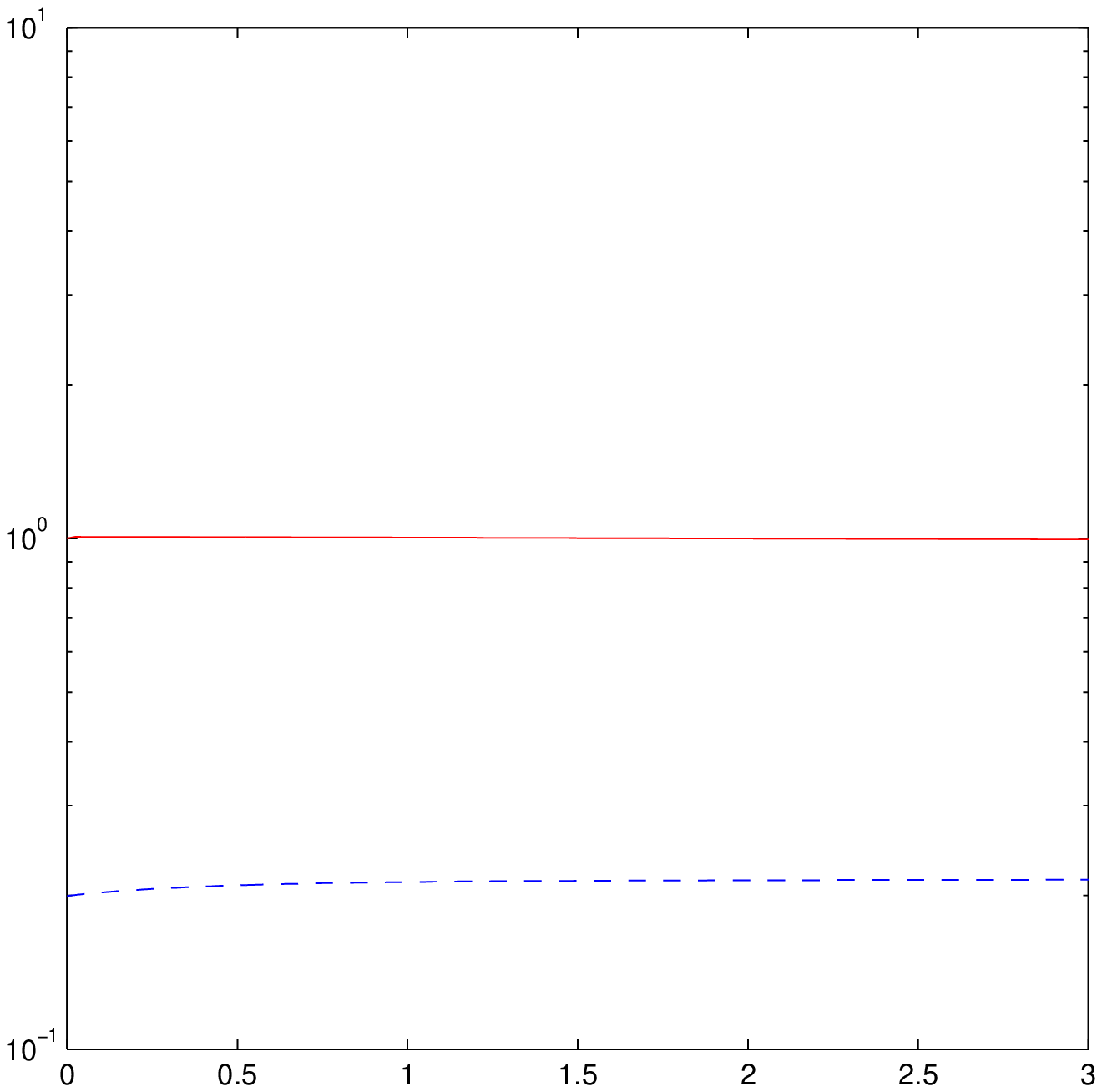}
\caption{$\kappa=1, m = 32, N = 8$}
\end{subfigure}
\begin{subfigure}[b]{0.3\textwidth}
\includegraphics[width=\textwidth]{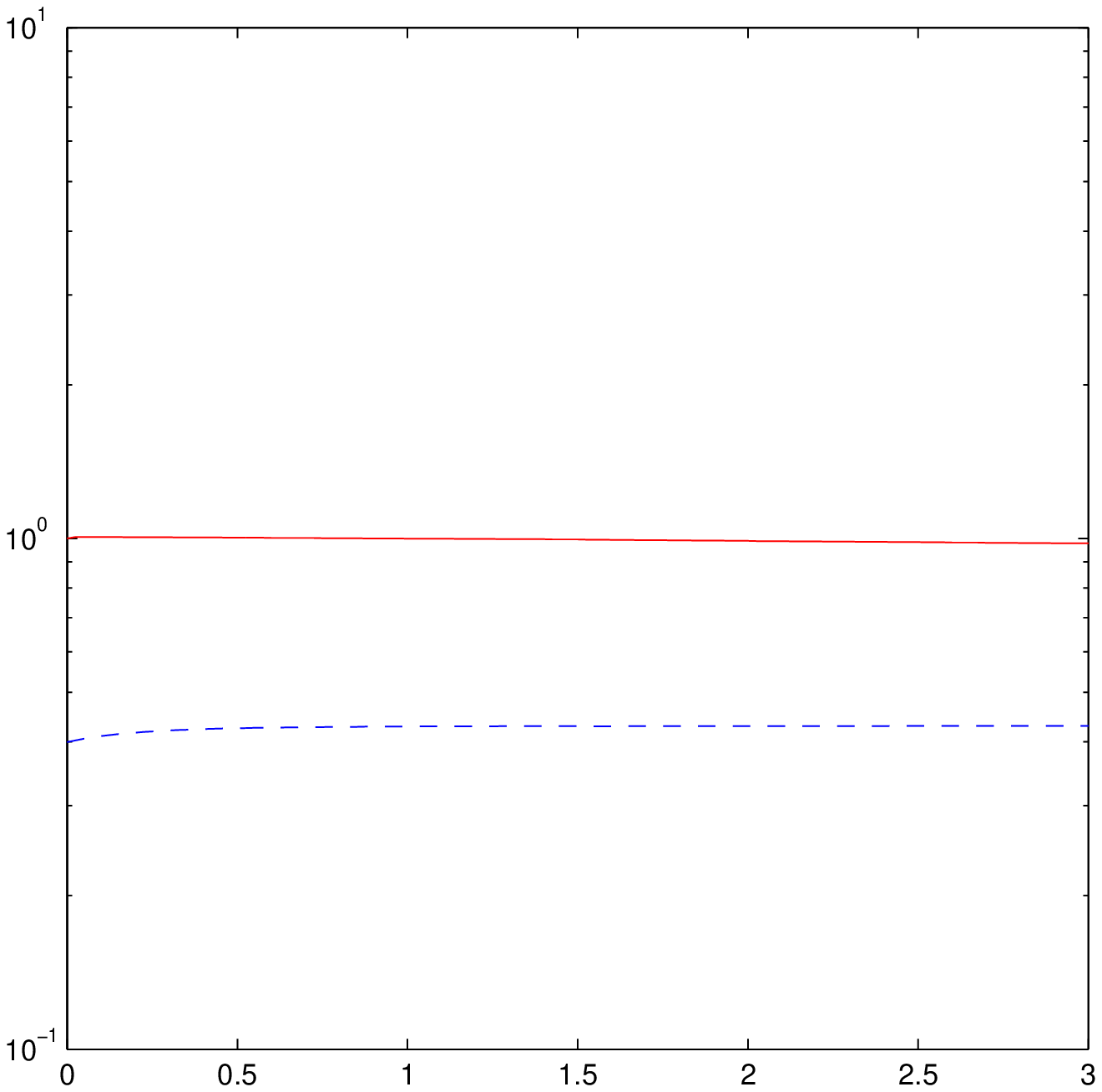}
\caption{$\kappa=2, m = 32, N = 8$}
\end{subfigure}
\begin{subfigure}[b]{0.3\textwidth}
\includegraphics[width=\textwidth]{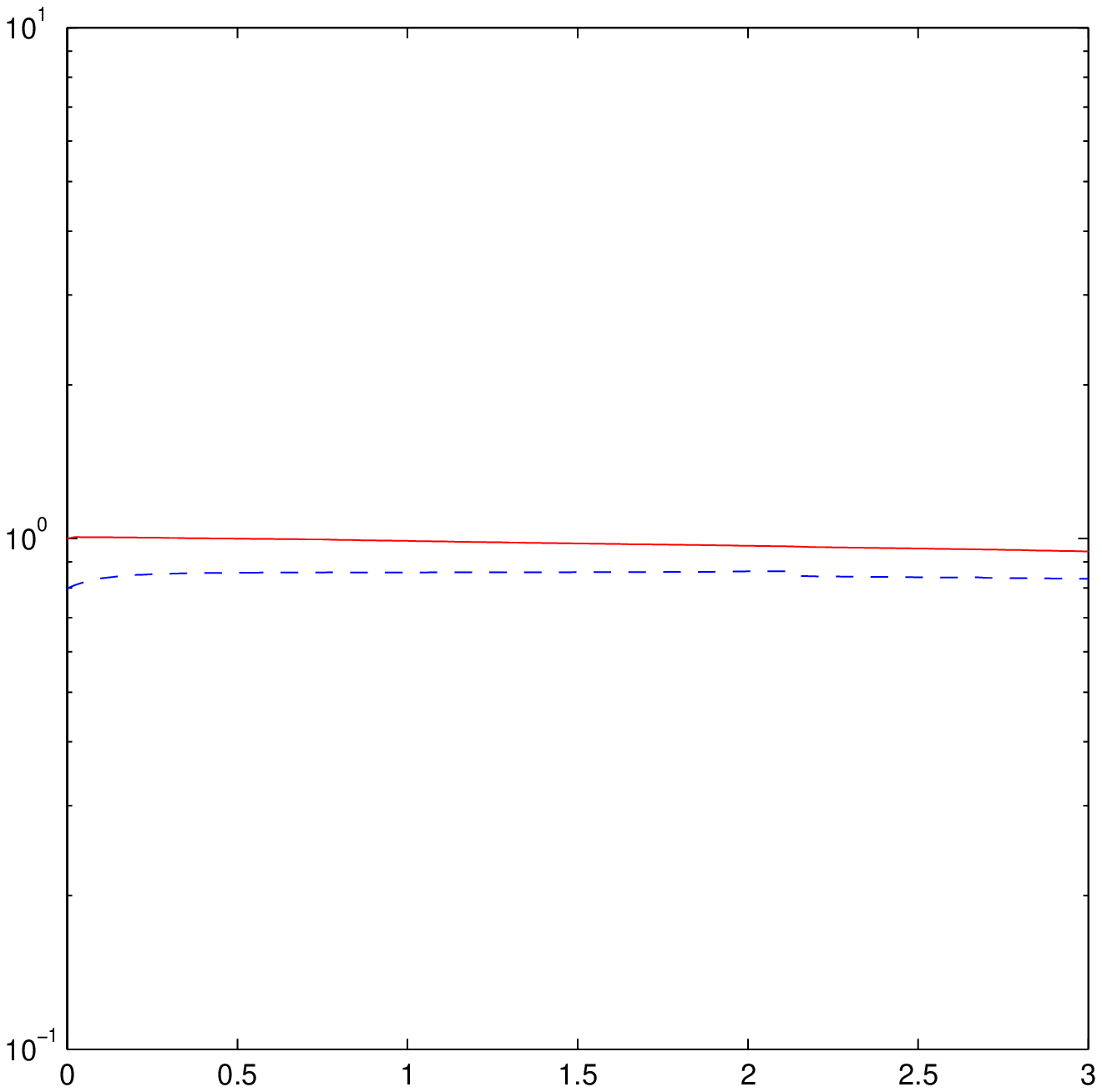}
\caption{$\kappa=4, m = 32, N = 8$}
\end{subfigure}
\begin{subfigure}[b]{0.3\textwidth}
\includegraphics[width=\textwidth]{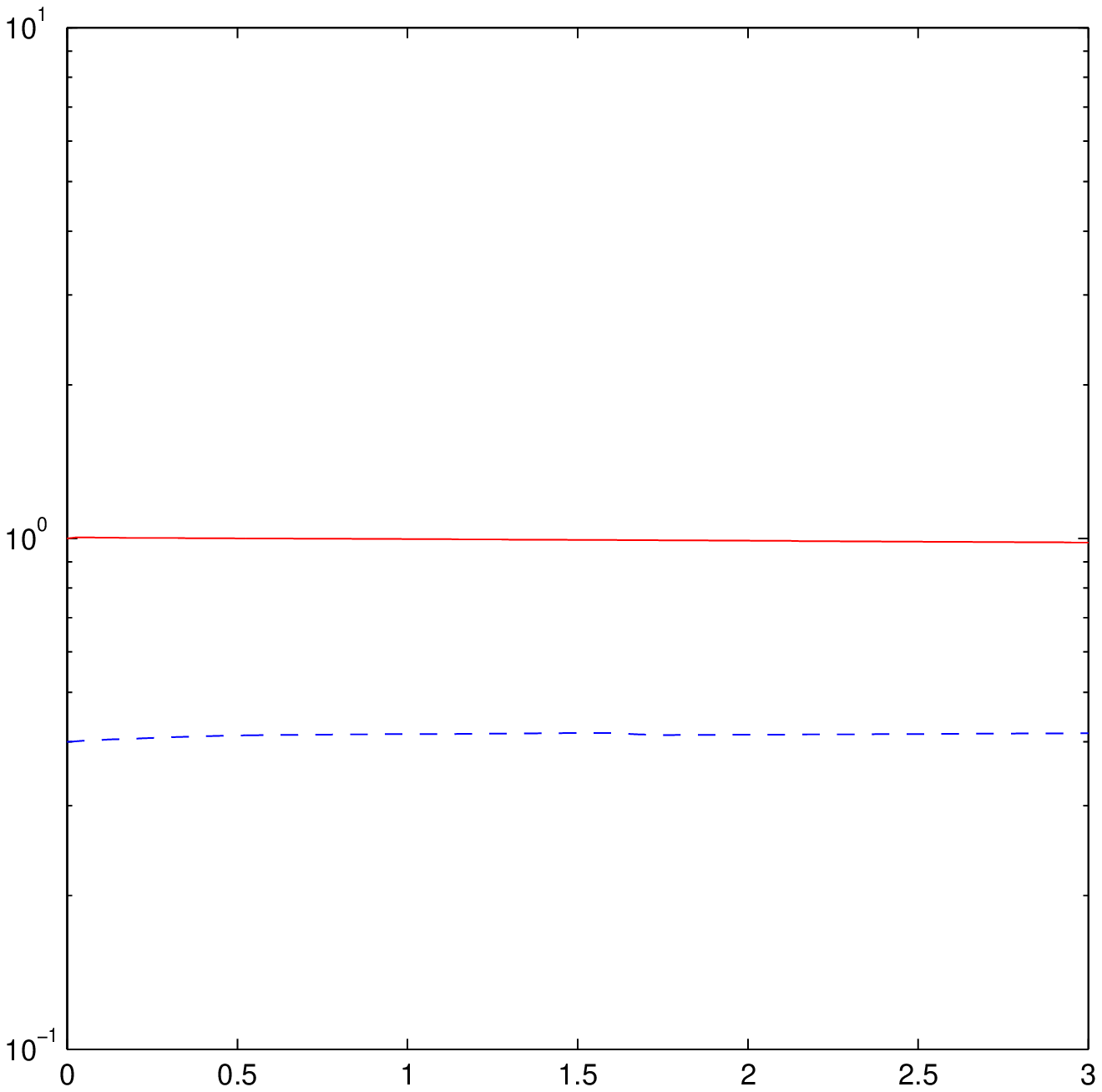}
\caption{$\kappa=1, m = 64, N = 16$}
\end{subfigure}
\begin{subfigure}[b]{0.3\textwidth}
\includegraphics[width=\textwidth]{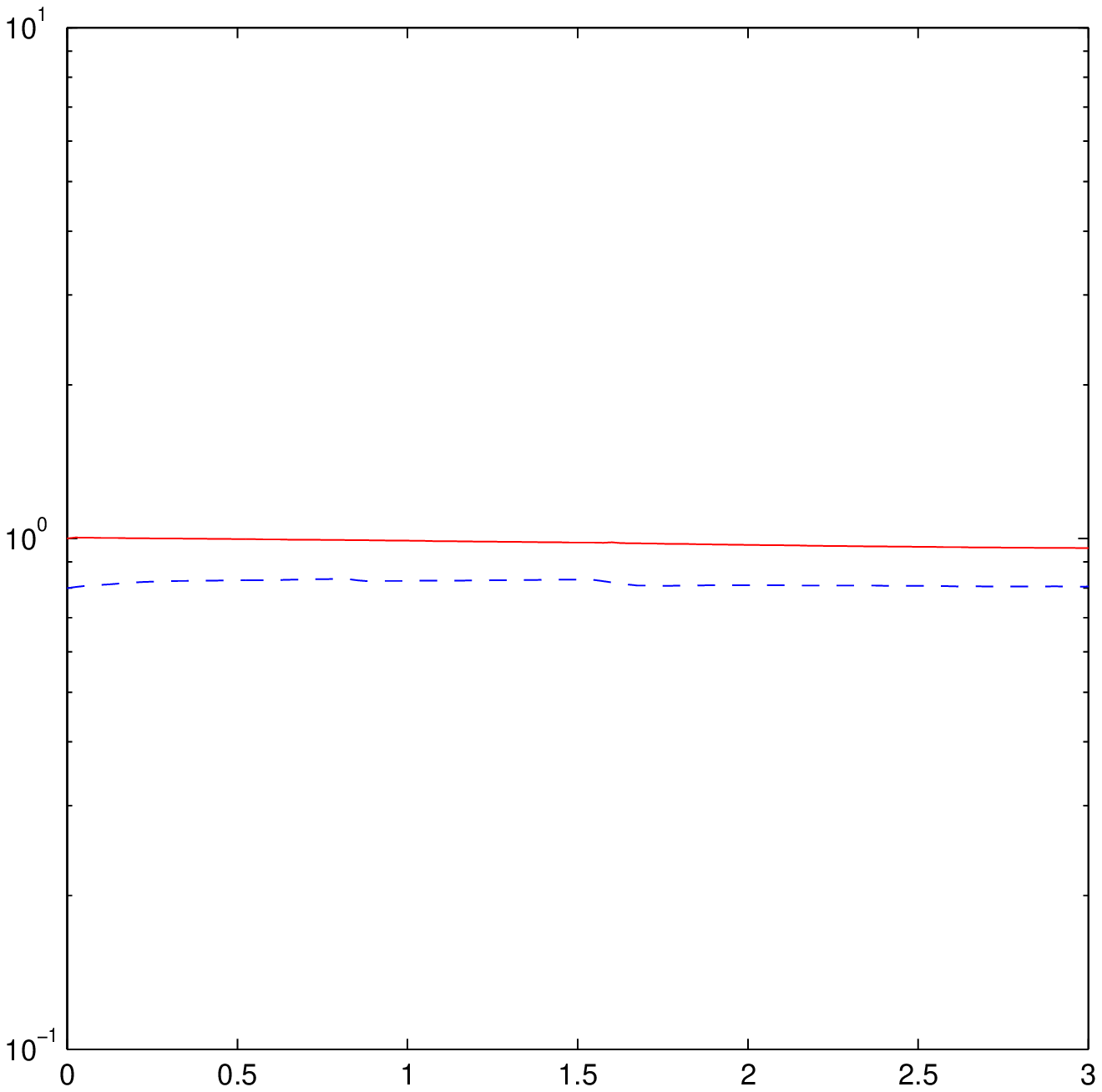}
\caption{$\kappa=2, m = 64, N = 16$}
\end{subfigure}
\begin{subfigure}[b]{0.3\textwidth}
\includegraphics[width=\textwidth]{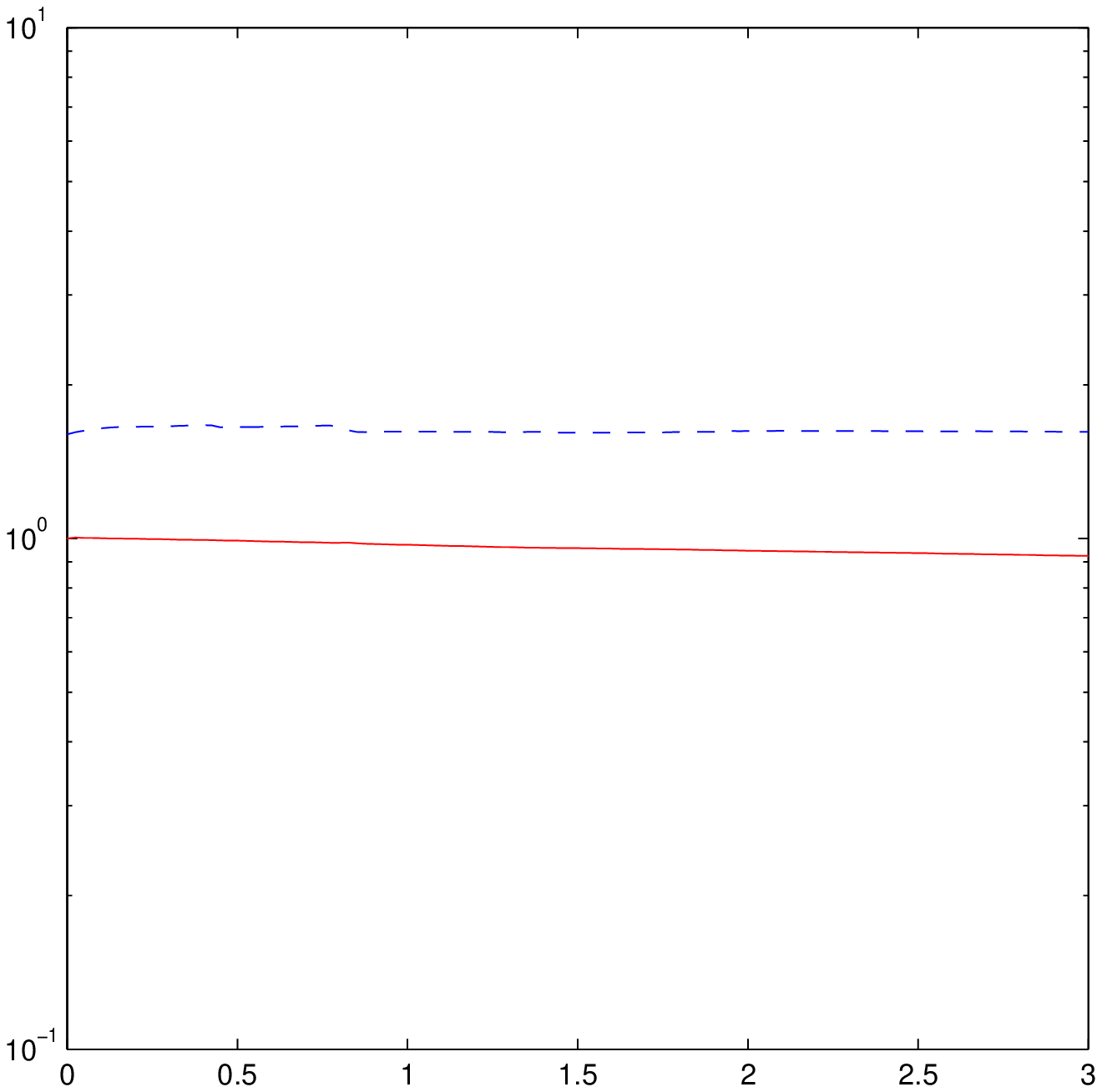}
\caption{$\kappa=4, m = 64, N = 16$}
\end{subfigure}
\begin{subfigure}[b]{0.3\textwidth}
\includegraphics[width=\textwidth]{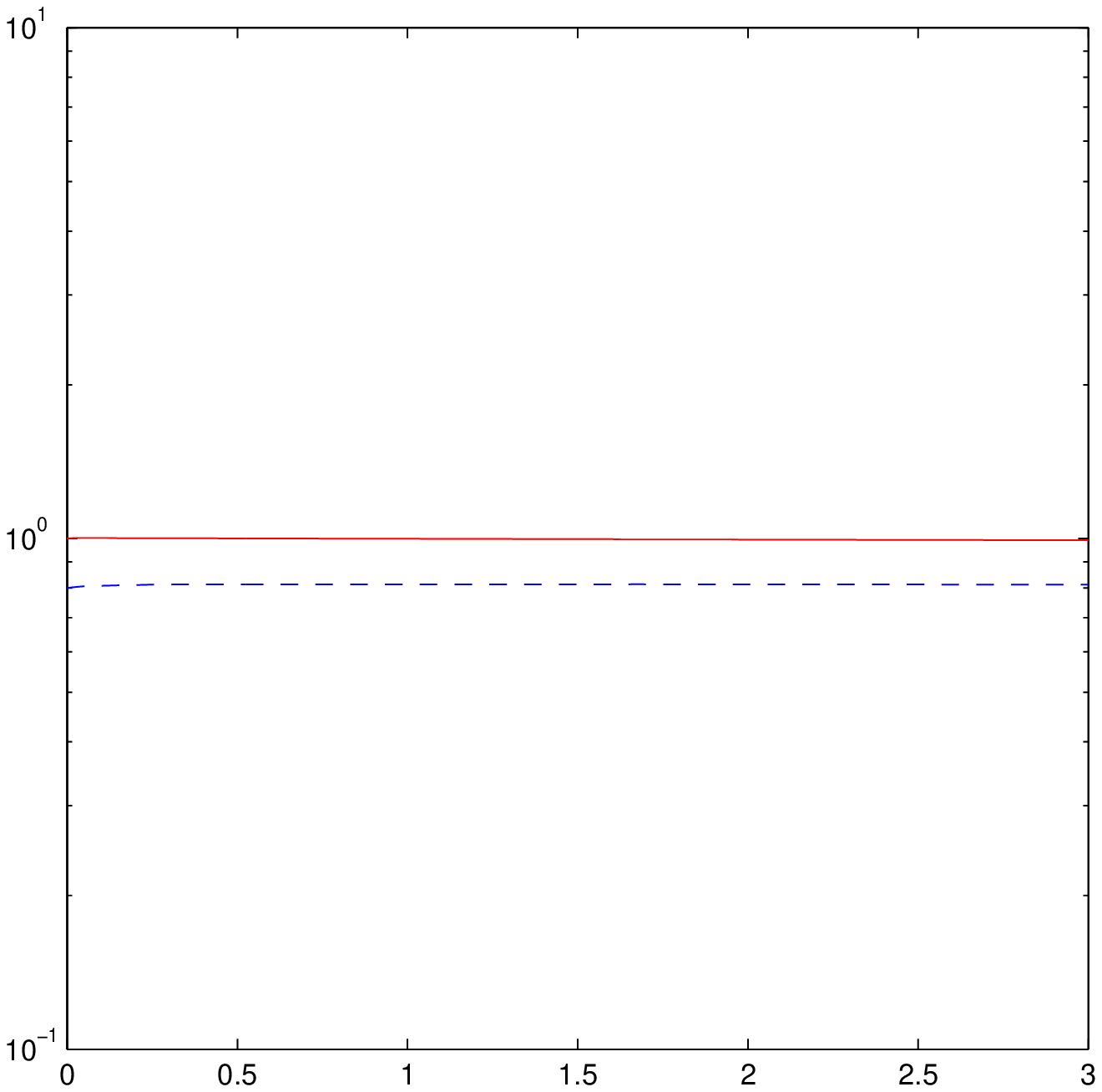}
\caption{$\kappa=1, m = 128, N = 32$}
\end{subfigure}
\begin{subfigure}[b]{0.3\textwidth}
\includegraphics[width=\textwidth]{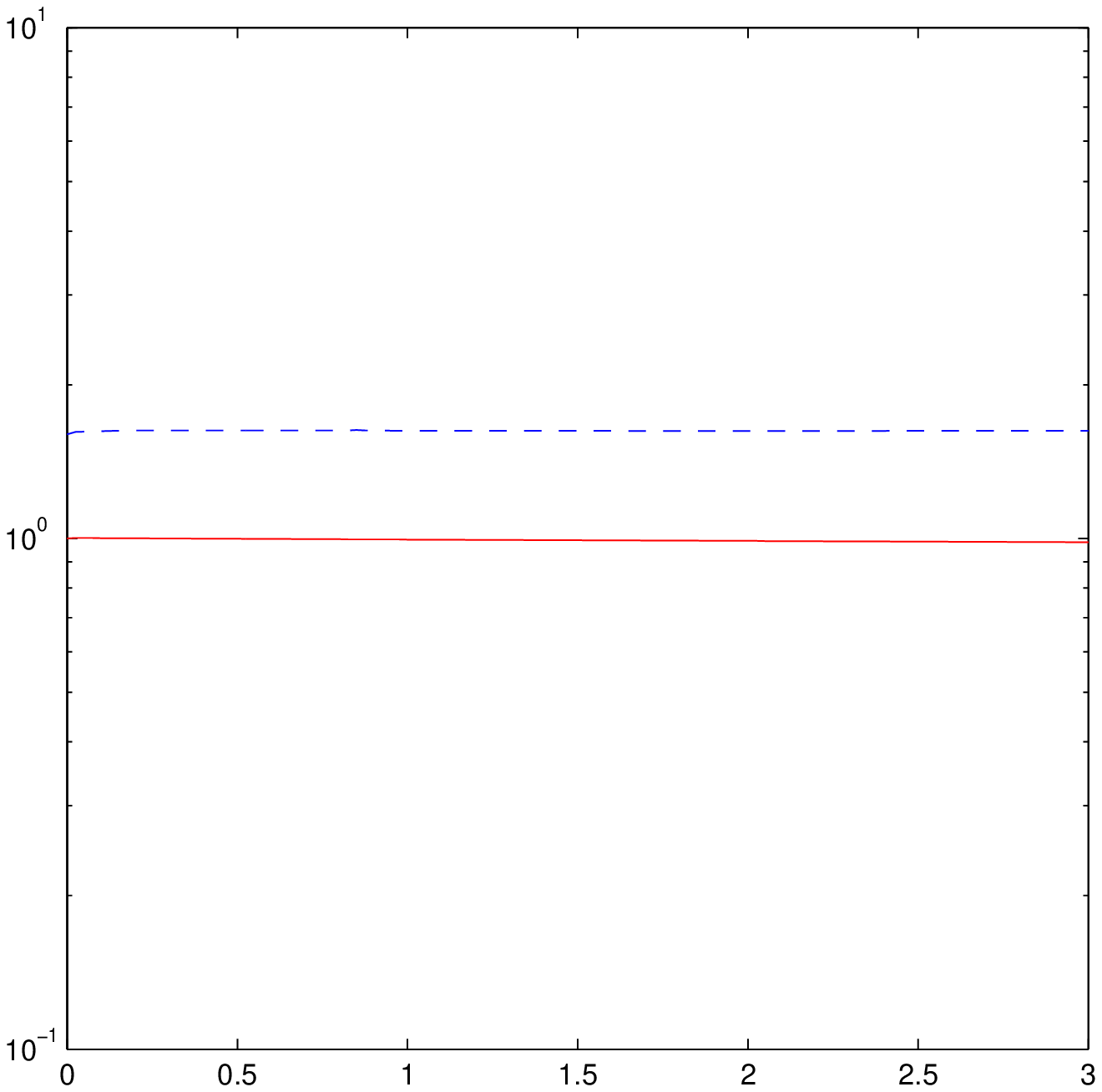}
\caption{$\kappa=2, m = 128, N = 32$}
\end{subfigure}
\begin{subfigure}[b]{0.3\textwidth}
\includegraphics[width=\textwidth]{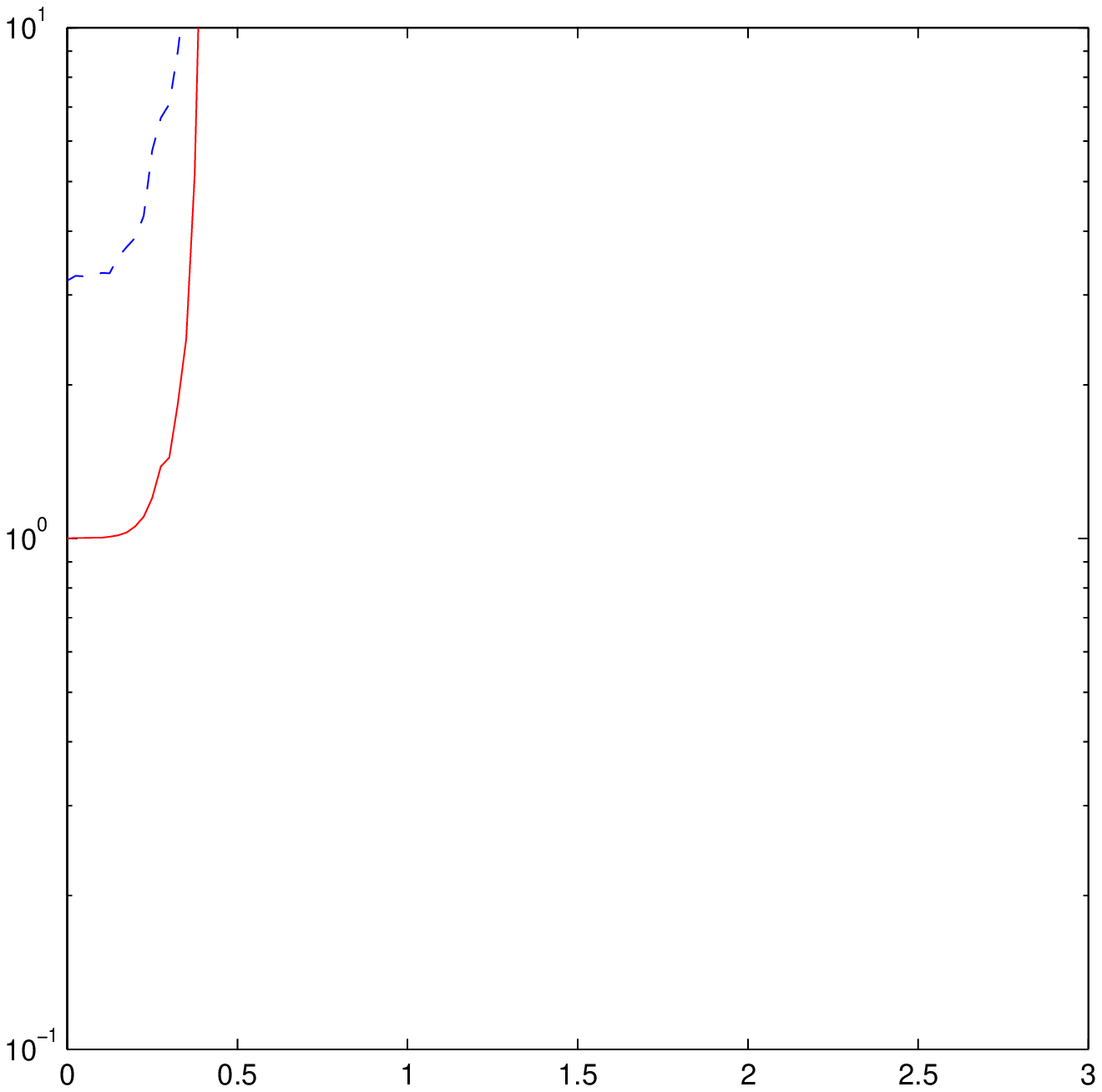}
\caption{$\kappa=4, m = 128, N = 32$}
\end{subfigure}
\caption{\revb{History of $E^n$ (red solid line) and $\eta^n$ (blue dashed line) throughout $t = 0$ and $t = 3$ with $K = 120$ in experiment \ref{sec:secexp5}.}}
\label{fig:exp5c}
\end{figure}

Figure \ref{fig:exp5} records the history of $E^n$ and $\eta^n$ throughout $t = 0$ and $t = 3$. It can be observed that our method is stable with the combinations of $K = 120, \kappa = 2$ and $K = 300, \kappa = 4$, in which FE-IBM is unstable. This shows our method provides good energy stability. \revb{Figures \ref{fig:exp5b} and \ref{fig:exp5c} records the same quantities with varying mesh sizes $m$ and $N$. Our analysis in Corollary \ref{stab3} suggests that the ratio $m/N$ should be fixed and the ratio $m \kappa / K$ should be sufficiently small for stability. From figures \ref{fig:exp5b} and \ref{fig:exp5c}, it can be observed that for a fixed ratio $m = 4N$, it is sufficient to ensure $m \kappa / K$ does not exceed a threshold of around   $32/15$ in order to achieve stability.}

\section{Conclusion}

In this paper, we develop a new staggered discontinuous Galerkin immersed boundary method.
We use the so-called BE/FE scheme for temporal discretization in order to avoid implicit coupling of nonlinear equations. \revb{Stability of our scheme is thus subject to the CFL type time-step restriction.}
We discuss our staggered discontinuous Galerkin scheme for solving the incompressible Navier-Stokes equations, and also a variational way of treating the fluid-structure interaction which suits our method.
The novel splitting of the convection term and the diffusion term realizes the possibility of involving the convection term without loss of energy stability.
Another important feature of our method is the improvement in volume conservation through the use of pointwise divergence-free post-processed velocity in driving the Lagrangian markers of the immersed boundary.
From the numerical experiments, we see that the exact divergence-free velocity field provides excellent volume conservation properties for the immersed boundary, the robustness of our method in treating immersed curves of different shapes, and also the energy stability of the nonlinear fluid model.
\revb{For a stretched immersed boundary model, we observe that the area conservation heavily depends on
a balance in the number of divisions $N$ in Eulerian grid and $m$ in Lagrangian grid.}




\end{document}